\documentclass[preprint,12pt]{elsarticle}
\usepackage[hmargin=0.8in,vmargin=0.9in]{geometry}
\usepackage{fancyhdr}
\usepackage{graphicx}
\usepackage{subfigure}
\usepackage{amssymb}
\usepackage{amsmath}
\usepackage{amsfonts}
\usepackage{theorem}
\usepackage{mathrsfs}
\usepackage{mathtools}
\usepackage{makecell}
\usepackage{multirow}
\usepackage{empheq}
\usepackage{algorithm}
\usepackage{algpseudocode}
\usepackage{bm}
\usepackage{color}
\algnewcommand\algorithmicforeach{\textbf{for each}}
\algdef{S}[FOR]{ForEach}[1]{\algorithmicforeach\ #1\ \algorithmicdo}
\usepackage{algpseudocode}
\usepackage{varwidth}
\usepackage{setspace}
\usepackage[labelfont=bf]{caption}
\usepackage{epstopdf}
\usepackage{float}
\usepackage{cite}

\makeatletter
\def\ps@pprintTitle{%
   \let\@oddhead\@empty
   \let\@evenhead\@empty
   \let\@oddfoot\@empty
   \let\@evenfoot\@oddfoot
}
\makeatother

\usepackage{booktabs} 
\usepackage{array} 
\usepackage{paralist} 
\usepackage{verbatim} 
\usepackage{subfigure} 
\usepackage{natbib}
\usepackage{hyperref}
\usepackage{xstring}

\usepackage{hyperref}
\hypersetup{
	colorlinks   = true,
	citecolor    = blue,
	linkcolor=blue  
}

\usepackage{xstring}
\makeatletter

{\theorembodyfont{\rmfamily}\newtheorem{remark}{Remark}[section]}

\allowdisplaybreaks
\numberwithin{equation}{section}
\numberwithin{figure}{section}
\numberwithin{table}{section}

\usepackage{ae,aecompl}

\def\Chi {{\Large{\mbox{$\chi$}}}}

\newcommand{\mF}{\bm{F}}
\newcommand{\mS}{\bm{S}}

\newcommand{\mG}{\bm{G}}

\newcommand{\mH}{\bm{H}}

\newcommand{\mU}{\bm{U}}

\newcommand{\Ups}{\bm{\Upsilon}}

\newcommand\eref[1]{(\ref{#1})}%
\newcommand\tref[1]{Table \ref{#1}}%
\newcommand\fref[1]{Fig.~\ref{#1}}%
\makeatletter
\newcommand{\refcheckize}[1]{%
	\expandafter\let\csname @@\string#1\endcsname#1%
	\expandafter\DeclareRobustCommand\csname relax\string#1\endcsname[1]{%
		\csname @@\string#1\endcsname{##1}\wrtusdrf{##1}}%
	\expandafter\let\expandafter#1\csname relax\string#1\endcsname
}
\makeatother

\newcommand*\xbarsh[1]{%
	\hbox{%
		\vbox{%
			\hrule height 0.5pt 
			\kern0.4ex
			\hbox{%
				\kern-0.22em
				\ensuremath{#1}%
				\kern-0.05em
			}%
		}%
	}%
}

\newcommand*\xbar[1]{%
	\hbox{%
		\vbox{%
			\hrule height 0.5pt 
			\kern0.4ex
			\hbox{%
				\kern-0.05em
				\ensuremath{#1}%
				\kern-0.05em
			}%
		}%
	}%
}

\begin{document}

\begin{frontmatter}

\author{Yekaterina Epshteyn\textsuperscript{1}\fnref{label2} and Thuong Nguyen\fnref{label3}}
\fntext[label2]{Department of Mathematics, The University of Utah, Salt Lake City, UT 84112, USA;
	{\tt epshteyn@math.utah.edu} }
\fntext[label3]{Department of Mathematics, The University of Utah, Salt Lake City, UT 84112, USA;
	{\tt tnguyen@math.utah.edu} }

\title{Adaptive Central-Upwind Scheme on Triangular Grids for the
	Saint-Venant System}


\begin{abstract}
In this work we develop a robust adaptive well-balanced and
positivity-preserving central-upwind scheme on unstructured triangular
grids for shallow water equations. The numerical method is an extension of the scheme from  [{\sc Liu {\em et al.}},J. of
Comp. Phys, 374 (2018), pp. 213 - 236]. As a part of the adaptive
central-upwind algorithm, we obtain  local a posteriori error
estimator for the efficient mesh refinement strategy.  The accuracy, high-resolution and
efficiency of new adaptive central-upwind scheme are demonstrated on a number
of challenging tests for shallow water models.
\end{abstract}

\end{frontmatter}

	\section{Introduction}\label{sec1}
	We consider the two-dimensional (2-D) Saint-Venant system of shallow water equations,
	
	\begin{subequations}\label{eq:swe}
		\begin{align}
		&h_t+(hu)_x+(hv)_y=0,\label{eq:swe1}\\
		&(hu)_t+\Big(hu^2+\frac{g}{2}h^2\Big)_x+(huv)_y=-ghB_x,\label{eq:swe2}\\
		&(hv)_t+(huv)_x+\Big(hv^2+\frac{g}{2}h^2\Big)_y=-ghB_y,\label{eq:swe3}
		\end{align}
	\end{subequations}
	
	where $t$ is the time, $x$ and $y$ are horizontal spatial coordinates ($(x,y)\in\Omega$), $h(x,y,t)$ is the water height, $u(x,y,t)$ and
	$v(x,y,t)$ are the $x$- and $y$-components of the flow velocity, $B(x,y)$ is the bottom topography, and $g$ is the constant gravitational
	acceleration. The system \eref{eq:swe1}--\eref{eq:swe3} was originally proposed in \citep{Sai}, but it is still widely used to model water flow
	in rivers, lakes and coastal areas, to name a few
        examples. The Saint-Venant system  \eref{eq:swe1}--\eref{eq:swe3}  is an example
        of the hyperbolic system of balance/conservation laws. The
        design of robust and accurate numerical algorithms for the computation
of its solutions is important and challenging problem that has been
extensively studied in the recent years.
\par An accurate numerical scheme for shallow water equations
\eref{eq:swe1}--\eref{eq:swe3} should preserve the physical properties
of the flow. For example, {\bf i)} the numerical method should be positivity preserving,
that is, the water height $h$ should be nonnegative at all times. The
positivity preserving property ensures a robust
performance of the algorithm on dry ($h$ is zero) or almost dry ($h$ is near zero)
states;  {\bf ii)}  in addition,  the numerical method for system
\eref{eq:swe1}--\eref{eq:swe3}  should be  well-balanced,
the method should exactly preserve the ``lake-at-rest'' solution,
$h+B\equiv const, u\equiv 0, v\equiv 0.$ This property diminishes the
appearance of unphysical oscillations of magnitude proportional to the
grid size. In the past decade, several well-balanced \citep{MR2562498,MR2799526,
  MR2799646, ABBLP,GHS,Jin,JW,KLab,KP1,SDBL,LeV,LeVbook,NPPN,PS,Rus1,Rus2,XS1,XS2,MR2799526,MR3640671,DLGeorge,MR3738315, MR2804645, AIPBEKP, LAEK} and positivity preserving
\citep{ABBLP,KLab,KP1,PS,MR3640671, MR2513368, MR2799646,MR3738315,
  MR2804645, AIPBEKP, LAEK} schemes (non exhaustive lists) for shallow
water models have been
proposed, but only few satisfy both major properties  {\bf i)}  and
 {\bf ii)} simultaneously. 
\par The traditional numerical methods for system
\eref{eq:swe1}--\eref{eq:swe3} consider very fine
fixed meshes to reconstruct delicate features of the
solution. However, this can lead to high computational cost, as well
as to a poor accuracy of small scale characteristics of the
problem. Therefore, the main goal of this work is to design adaptive
numerical algorithm for shallow water equations. 
 In this work, we extend numerical method in \citep{LAEK} to adaptive well-balanced and
 positivity-preserving central-upwind
finite volume method on unstructured triangular grids. The central
Nessyahu-Tadmor schemes, their generalization into higher resolution
central schemes and semidiscrete central-upwind schemes are a family
of efficient and accurate Godunov-type Riemann problem-free
projection-evolution finite volume methods for hyperbolic
problems. They were originally developed in
\citep{MR1047564,MR1756766,KNP}. The main advantages of these
numerical algorithms are the high-resolution, efficiency and their simplicity. The class of
central-upwind methods has been successfully used for problems in
science and engineering, including, for geophysical flow problems
and related models,
e.g. \citep{MR1756766,MR1763829,KNP,KLab,SDBL,KP1,KM,CEHK,MR3229988,MR3187922,KurLiu,MR3071176,KP05,AIPBEKP,MR2804645,MR3440159,LAEK}. There
is some very recent effort on the design of adaptive well-balanced and
positivity-preserving central-upwind schemes on quad-tree grids for
shallow water models \citep{MR3315267,GHAZIZADEH2020104633}, but no research has been
done for the development of such adaptive schemes on unstructured
triangular grids.
\par This paper is organized as follows. In Section \ref{sect2}, we briefly
review well-balanced positivity-preserving central-upwind scheme
on unstructured triangular grids \citep{LAEK} which serves as the
underlying discretization for the developed adaptive algorithm.  We
give summary of the adaptive central-upwind method in Section \ref{sect3a0}. We
discuss adaptive mesh refinement strategy in
Section \ref{sect3a}. In Section \ref{sect3b},  we present adaptive
second-order strong stability preserving Runge-Kutta method, employed
as a part of time evolution for the adaptive central-upwind scheme. We
derive local a posteriori error estimator in Section \ref{sect3c} which is
used as a robust indicator for the adaptive mesh refinement in our
work. Finally, in Section \ref{sect4}, we illustrate the high accuracy and
efficiency of the developed adaptive central-upwind scheme on a number
of challenging tests for shallow water models.

	\section{Semi-Discrete Central-Upwind Scheme \!--\! an
		Overview}\label{sect2}
	In this work, we employ the central-upwind scheme
        discussed in this section as the underlying
	discretization for the developed adaptive central-upwind
        algorithm,  Section \ref{sect3}.  Therefore, in this section, we will briefly review a 
	semi-discrete second-order well-balanced positivity preserving
	central-upwind scheme on unstructured triangular grids for the Saint-Venant system of shallow water
	equations \citep{MR2804645,LAEK}. 
\par In the first work \citep{MR2804645}, a new second-order semi-discrete central-upwind scheme was developed for computing the solutions of the system \eref{eq:swe1}--\eref{eq:swe3} on
	unstructured triangular grids.
	The key ideas in the
	development of the scheme in \citep{MR2804645} were: {\bf 1)} Change of variables  from $(h,hu,hv)^T$ to variables $(w:=h+B,hu,hv)^T$. This change of variables simplifies the
	construction of the well-balanced scheme since in the ``lake-at-rest''
	steady-state, it is the equilibrium variable, the water
        surface $w\equiv h+B$ (but not the conservative variable, the
        water height $h$) that
	has to stay constant; {\bf 2)}
	Replacement of the bottom topography function $B$ with its continuous piecewise linear approximation;
	{\bf 3)} Design of the special positivity preserving correction of the  piecewise linear reconstruction for the water surface $w$;
	{\bf 4)} Development of a special well-balanced finite-volume-type quadrature
	for the discretization of the cell averages of the geometric source
	term. The developed scheme in \citep{MR2804645}, enforced the positivity of
	the water height $h$,  and preserved the ``lake-at-rest'' steady state
	in the case of fully submerged bottom topography. In the recent work \citep{LAEK}, we further improved the
        well-balanced property of the scheme from
	\citep{MR2804645},  and extended the scheme to accurate and stable simulations
	of shallow water models with dry or near dry states (e.g., waves arriving or leaving the shore).
We will briefly review below the central-upwind scheme from \citep{LAEK}.
	
	First, we rewrite the system \eref{eq:swe1}--\eref{eq:swe3} in the following equivalent form,
	\begin{equation}\label{eqs2}
	\mU_t+\mF(\mU,B)_x+\mG(\mU,B)_y=\mS(\mU,B),
	\end{equation}
	where the variables $\bm{U}$ and the fluxes $\bm{F}$ and $\bm{G}$ are
	\begin{equation*}
	\mU=\left(\begin{array}{c}w\\[1.8ex]hu\\[1.8ex]hv\end{array}\right),\quad
	\mF=\left(\begin{array}{c}hu\\\dfrac{(hu)^2}{w-B}+\dfrac{g}{2}(w-B)^2\\[1.5ex]\dfrac{(hu)(hv)}{w-B}\end{array}\right),\quad
	\mG=\left(\begin{array}{c}hv\\\dfrac{(hu)(hv)}{w-B}\\[1.5ex]\dfrac{(hv)^2}{w-B}+\dfrac{g}{2}(w-B)^2\end{array}\right),
	\end{equation*}
	and the source term $\mS$ is
	\begin{equation*}
	\mS=\left(\begin{array}{c}0\\-g(w-B)B_x\\[0.5ex]-g(w-B)B_y\end{array}\right).
	\end{equation*}
	\begin{figure}[ht!]
	\centering
\includegraphics[width=8cm]{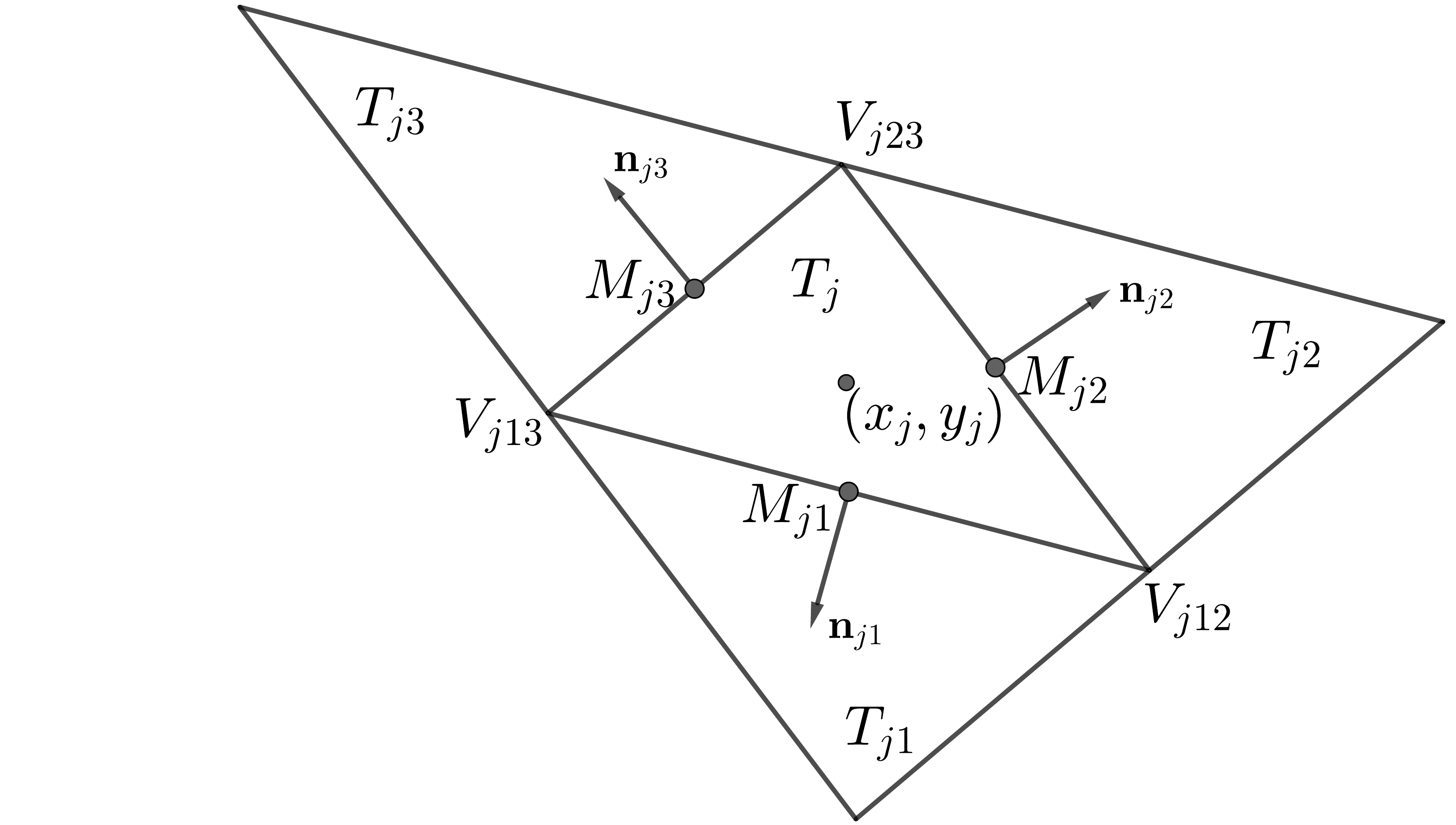}\\
		\vspace*{5mm}
	\caption{{ A typical triangular cell with three neighbors.\label{Fig1}}}
	\vspace*{2mm}
	\end{figure}

	As illustrated in Figure \ref{Fig1},  we denote,
	
	\noindent
	${\mathcal {T}}:=\bigcup_j T_j$ is an unstructured triangulation of the computational domain $\Omega$;
	
	\noindent
	$T_j\in \mathcal{T}$ is a triangular cell of size $|T_j|$ with the barycenter $(x_j,y_j)$;
	
	\noindent
	$V_{j\kappa}=(\widetilde{x}_{j\kappa},\widetilde{y}_{j\kappa}),~\kappa=12,23,31$ are the three vertices of $T_j$;
	
	\noindent
	$T_{jk},~k=1,2,3$ are the neighboring triangles that share a common side with $T_j$;
	
	\noindent
	$\ell_{jk}$ is the length of the common side of $T_j$ and
        $T_{jk}$,  and $M_{jk}$ is its midpoint;
	
	\noindent
	$\bm{n}_{jk}:=(\cos(\theta_{jk}),\sin(\theta_{jk}))^\top$ is the outer unit normal to the $k$th side of $T_j$.
	
	Next, in order to develop the positivity-preserving and well-balanced
	scheme, the bottom topography $B$ is replaced with its continuous piecewise linear approximation $\widetilde{B}$ given by
	\begin{equation*}
	\left|\begin{array}{ccc}x-\widetilde{x}_{j12}&y-\widetilde{y}_{j12}&{\widetilde B}(x,y)-\widehat{B}_{j12}\\[0.5ex]
	\widetilde{x}_{j23}-\widetilde{x}_{j12}&\widetilde{y}_{j23}-\widetilde{y}_{j12}&\widehat{B}_{j23}-\widehat{B}_{j12}\\[0.5ex]
	\widetilde{x}_{j13}-\widetilde{x}_{j12}&\widetilde{y}_{j13}-\widetilde{y}_{j12}&\widehat{B}_{j13}-\widehat{B}_{j12}\end{array}\right|=0,
	\quad(x,y)\in T_j,
	\end{equation*}
	where, in the case of continuous bottom topography,
	$\widehat{B}_{j\kappa}:=B(V_{j\kappa})~\kappa=12,23,31$. Then, denote:
	$$
	B_{jk}:=\widetilde{B}(M_{jk}),\quad B_j:=\widetilde{B}(x_j,y_j)=\frac{1}{3}(\widehat{B}_{j12}+\widehat{B}_{j23}+\widehat{B}_{j13}).
	$$
	
	At time $t$, define by $\xbar\mU_j(t)$ the approximation of the cell averages of the solution,
	\begin{equation*}
	\xbar\mU_j(t)\approx\frac{1}{|T_j|}\iint\limits_{T_j}\mU(x,y,t)\,dxdy.
	\end{equation*}
	Then, it can be shown (see \citep{LAEK, MR2804645}), that the semi-discrete
	second-order central-upwind scheme for the Saint-Venant system
	(\ref{eqs2})  on triangular grid is given by the following system of ODEs,
	\begin{equation}
	\begin{aligned}
	\frac{d\,\xbar{\mU}_j}{dt}=&-\frac{1}{|T_j|}\big[\mH_{j1}+\mH_{j2}+\mH_{j3}\big]+\xbar{\mS}_j,
	\label{10eq8}
	\end{aligned}
	\end{equation}
	where the numerical fluxes through the edges of the triangular
        cell $T_j$ are
	\begin{equation}
	\begin{aligned}
	\mH_{jk}=&\frac{\ell_{jk}\cos(\theta_{jk})}{a_{jk}^{\rm in}+a_{jk}^{\rm out}}
	\Big[{a_{jk}^{\rm in}\mF(\mU_{jk}(M_{jk}),B_{jk})+a_{jk}^{\rm out}\mF(\mU_j(M_{jk}),B_{jk})}\Big]\\
	+&\frac{\ell_{jk}\sin(\theta_{jk})}{a_{jk}^{\rm in}+a_{jk}^{\rm out}}
	\Big[{a_{jk}^{\rm in}\mG(\mU_{jk}(M_{jk}),B_{jk})+a_{jk}^{\rm out}\mG(\mU_j(M_{jk}),B_{jk})}\Big]\\
	-&\ell_{jk}\frac{a_{jk}^{\rm in}a_{jk}^{\rm out}}{a_{jk}^{\rm in}+a_{jk}^{\rm out}}\big[\mU_{jk}(M_{jk})-\mU_j(M_{jk})\big],\quad k=1,2,3.
	\label{flux}
	\end{aligned}
	\end{equation}
	
	Here, $\bm{U}_{j}(M_{jk})$ and $\bm{U}_{jk}(M_{jk})$ are the
        reconstructed point values of $\mU$ at the middle points of
        the edges $M_{jk}$. To obtain these values \citep{LAEK}, first, 
	a piecewise linear reconstruction of the variables
        $\Ups:=(w,u,v)^\top$ is computed as,
	\begin{equation}
	\widetilde\Ups(x,y)=\sum\limits_j\Ups_j(x,y)\Chi_{T_j},\quad\Ups_j(x,y):=\Ups_j+(\widehat{\Ups}_x)_j(x-x_j)+(\widehat{\Ups}_y)_j(y-y_j),
	\label{pwlinear}
	\end{equation}
	where $\Chi_{T_j}$ is the characteristic function of the cell $T_j$, $\Ups_j$ are the point values of $\Ups$ at the cell centers and
	$(\widehat{\Ups}_x)_j$ and $(\widehat{\Ups}_y)_j$ are the limited partial derivative. After that, the second and third components of the point values $\bm{U}_{j}(M_{jk})$ and $\bm{U}_{jk}(M_{jk})$ are obtained from,
$\Ups_j(M_{jk})$ and $\Ups_{jk}(M_{jk})$,
$$
\begin{aligned}
&(hu)_j(M_{jk})=(w_j(M_{jk})-B_{jk})\,u_j(M_{jk}),&&(hu)_{jk}(M_{jk})=(w_{jk}(M_{jk})-B_{jk})\,u_{jk}(M_{jk}),\\
&(hv)_j(M_{jk})=(w_j(M_{jk})-B_{jk})\,v_j(M_{jk}),&&(hv)_{jk}(M_{jk})=(w_{jk}(M_{jk})-B_{jk})\,v_{jk}(M_{jk}).
\end{aligned}
$$
See Section 2 in \citep{LAEK} for more details on the reconstruction.	
	
Moreover, to design a well-balanced central-upwind scheme, a special
second-order  reconstruction of water surface is introduced in \citep{LAEK} which is positivity preserving for the steady-state solutions with partially flooded/dry cells. Hence, the linear approximation for the water surface is updated as follows:

\begin{itemize}
    \item In the dry cells in which $\xbar{w}_j=B_j$, the
      corresponding linear pieces for $w$ in \eref{pwlinear} are replaced by,
	\begin{equation}
	\widetilde{w}_j(x,y)=\widetilde{B}(x,y).
	\label{10eq12}
	\end{equation}
	\item If $T_j$ is a partially flooded which means $B_j<\xbar{w}_j<\max\{\widehat{B}_{j23},\widehat{B}_{j13},\widehat{B}_{j12}\}$, the water surface is reconstructed by using two linear pieces instead of one as,
	\begin{equation}
	\widetilde{w}_j(x,y)=\begin{cases}\mathring{w}_j(x,y),\quad &\text{if} \quad(x,y)\in T^{wet}_j,\\
	\widetilde{B}(x,y), \quad&\text{otherwise}, \end{cases}
	\label{10eq12}
	\end{equation}
	where $\mathring{w}_j(x,y)$ is a linear reconstruction of the
        water surface on the wet part $T^{wet}_j$ of the cell $T_j$.
	\item If $T_j$ is a fully flooded $\xbar{w}_j\geq\max\{\widehat{B}_{j23},\widehat{B}_{j13},\widehat{B}_{j12}\}$, no further modification for the linear approximation \eref{pwlinear} is needed.
\end{itemize}
See Section 3 in \citep{LAEK} for more details of the reconstruction of
the water surface $w$.
\par	In \eref{flux}, $a_{jk}^{\rm in}$ and $a_{jk}^{\rm out}$ are the one-sided local speeds of propagation in the directions $\pm\bm{n}_{jk}$.
	These speeds are related to the largest and smallest eigenvalues of the Jacobian matrix
	$J_{jk}=\cos(\theta_{jk})\,\frac{\partial {\bm F}}{\partial {\bm U}}+\sin(\theta_{jk})\,\frac{\partial {\bm G}}{\partial {\bm U}}$,
	denoted by $\lambda_+[J_{jk}]$ and $\lambda_-[J_{jk}]$, respectively, and are defined by
	\begin{equation}
	\begin{aligned}
	&a^{\rm in}_{jk}=-\min\{\lambda_-[J_{jk}(\mU_j(M_{jk}))],\,\lambda_-[J_{jk}(\mU_{jk}(M_{jk})],\,0\},\\
	&a^{\rm out}_{jk}=\max\{\lambda_+[J_{jk}(\mU_j(M_{jk}))],\,\lambda_+[J_{jk}(\mU_{jk}(M_{jk})],\,0\},
	\end{aligned}\label{speed}
	\end{equation}
	where
	$$
	\begin{aligned}
	&\lambda_\pm[J_{jk}(\mU_j(M_{jk}))]=\cos(\theta_{jk})u_j(M_{jk})+\sin(\theta_{jk})v_j(M_{jk})\pm\sqrt{gh_j(M_{jk})},\\
	&\lambda_\pm[J_{jk}(\mU_{jk}(M_{jk}))]=\cos(\theta_{jk})u_{jk}(M_{jk})+\sin(\theta_{jk})v_{jk}(M_{jk})\pm\sqrt{gh_{jk}(M_{jk})}.
	\end{aligned}
	$$
	\begin{remark}
		In order to avoid division by 0 (or by a very small positive number), the numerical flux \eref{flux} is replaced with
		$$
		\begin{aligned}
		\mH_{jk}=&\frac{\ell_{jk}\cos(\theta_{jk})}{2}\left[\mF(\mU_{jk}(M_{jk}),B_{jk})+\mF(\mU_j(M_{jk}),B_{jk})\right]\\
		+&\frac{\ell_{jk}\sin(\theta_{jk})}{2}\left[\mG(\mU_{jk}(M_{jk}),B_{jk})+\mG(\mU_j(M_{jk}),B_{jk})\right]
		\end{aligned}
		$$
		wherever $a_{jk}^{\rm in}+a_{jk}^{\rm out}<\sigma$. In
                all of the reported numerical examples in Section \ref{sect4}, we have taken $\sigma=10^{-6}$.
		\end{remark}	
\par	A fully discrete scheme can be obtained by numerically solving the ODE system \eref{10eq8}, \eref{flux} using a stable and sufficiently
		accurate ODE solver. The time-step size on each cell
                $T_j\in \mathcal{T}$ should satisfy the CFL-type condition (see \citep{MR2804645}), which can be expressed as,
		\begin{equation}
		\Delta t<\frac{1}{6}\min_{j,k}\left[\frac{r_{jk}}{\max(a^{\rm in}_{jk},a^{\rm out}_{jk})}\right],
		\label{1eqCFL}
		\end{equation}
		where $r_{j1}$, $r_{j2}$ and $r_{j3}$ are the three corresponding
		altitudes of the triangle $T_j$.  However, the condition
                \eref{1eqCFL} can become too restrictive on partially
                flooded cells. Thus, for partially flooded
                cells,  the ``draining'' time-step technique is used
                to ensure the positivity of the scheme without
                reducing the time step size \eref{1eqCFL},
                \citep{MR2799646,LAEK}. Namely, first  the ``draining'' time-step $\Delta t_j^{\rm drain}$ is defined by,
\begin{equation*}
\Delta t_j^{\rm drain}:=\frac{|T_j|\,\xbar h_j^{\,n}}{\sum\limits_{k=1}^3\max(0,H^{(1)}_{jk})}
\end{equation*}
where $H^{(1)}_{jk}$ is the first component of the numerical flux ${\bm H}_{jk}$ given by \eref{flux}. Notice that,
for fully flooded cells $\Delta t_j^{\rm drain}=\Delta t$, while for
dry cells $\Delta t_j^{\rm drain}=0$. Next, the local ``draining''
time-step $\Delta t_{jk}$ for each edge $k$ of the cell $T_j \in
\mathcal{T}$ is defined as, 
\begin{equation}
			\Delta t_{jk}=\begin{cases}
			\min(\Delta t, \Delta t^{drain}_j), \quad \text{if}\quad  H^{(1)}_{jk}> 0,\\			
			\min(\Delta t, \Delta t^{drain}_{jk}), \quad \text{if}\quad  H^{(1)}_{jk}\leq 0,\\
			\end{cases}
		\end{equation}
	where $\Delta t^{drain}_{jk}$ is the
        ``draining'' time-step in the neighboring triangle $T_{jk}\in
        \mathcal{T}$ of $T_j$ and $\Delta t$ is computed by
\eref{1eqCFL}, but with the minimum taken there over the flooded cells
only. This procedure of the draining time step is a part
of the adaptive SSPRK2 time evolution \eref{eq:rk1}-\eref{eq:rk2} in
Section \ref{sect3b}.
\par 	Finally, the cell average of the source term $\mS_j$ in \eref{10eq8},
	\begin{equation*}
	\xbar\mS_j(t)\approx\frac{1}{|T_j|}\iint\limits_{T_j}\mS\big(\mU(x,y,t),B(x,y)\big)\,dxdy,
	\end{equation*}
	has to be discretized in a well-balanced manner \citep{LAEK}:\\
Quadrature for $\xbar S^{\,(2)}_j$ is
\begin{equation}
\begin{aligned}
\xbar S^{\,(2)}_j&=\frac{g}{2|T_j|}\sum_{k=1}^3\ell_{jk}\cos(\theta_{jk})\cdot
\frac{\Delta t_{jk}}{\Delta t}\cdot\big[w(M_{jk})-B(M_{jk})\big]^2\\
&-\frac{g}{3}\Big[(w_{j12}-\widehat B_{j12})\,w_x(V_{j12})+(w_{j23}-\widehat B_{j23})\,w_x(V_{j23})+
(w_{j13}-\widehat B_{j13})\,w_x(V_{j13})\Big].
\label{qra1}
\end{aligned}
\end{equation}
A similar quadrature for $\xbar S^{\,(3)}_j$ is
\begin{equation}
\begin{aligned}
\xbar S^{\,(3)}_j&=\frac{g}{2|T_j|}\sum_{k=1}^3\ell_{jk}\sin(\theta_{jk})\cdot
\frac{\Delta t_{jk}}{\Delta t}\cdot\big[w(M_{jk})-B(M_{jk})\big]^2\\
&-\frac{g}{3}\Big[(w_{j12}-\widehat B_{j12})\,w_y(V_{j12})+(w_{j23}-\widehat B_{j23})\,w_y(V_{j23})+
(w_{j13}-\widehat B_{j13})\,w_y(V_{j13})\Big].
\label{qra2}
\end{aligned}
\end{equation}
\begin{remark}
Note, that in Section \ref{sect4}, we compare performance of the
developed adaptive central-upwind scheme Section \ref{sect3} with a performance of the
central-upwind scheme without adaptivity from \citep{LAEK} (see also
brief review above). We use
standard SSPRK2 time discretization \citep{GST} together with the draining time step
for the scheme without adaptivity from \citep{LAEK} in numerical experiments in Section \ref{sect4}.
\end{remark}


	
	\section{Adaptive Central-Upwind Scheme}\label{sect3}
	
	The traditional numerical schemes are
	based on the use of very fine fixed meshes to reconstruct delicate features of the solution. This
	can lead to high computational
	cost, as well as poor resolution of all small scale features of the
	problem. In many engineering and scientific applications,
	it is beneficial to use adaptive meshes for improving the accuracy of
	the approximation at a much lower cost. Therefore, in this
	section, we will introduce an efficient and accurate adaptive central-upwind algorithm.
	
	\subsection{Adaptive Central-Upwind Algorithm}\label{sect3a0}
	The adaptive central-upwind algorithm is described briefly by the following steps.
	
	\textbf{Step 0.} At time $t=t^0$, generate the initial uniform
        grid $\mathcal{T}^{0, 0}$.
	
	\textbf{Step 1.} On mesh $\mathcal{T}^{n, \mathcal{M}_n}$, evolve  the cell averages $\xbar{\mU}^n$ of
	the solution from time $t^n$ to $\xbar{\mU}^{n+1}$ at the
	next time level $t^{n+1} $ using adaptive central-upwind
        scheme \eref{eq:rk1}-\eref{eq:rk2}, see Section \ref{sect3}: 	
	\begin{itemize}
		\item  At time $t^n$,  determine the level $l=0, 1,
                  ..., L$ of each cell/triangle
                  $T^{n, \mathcal{M}_n}_j \in \mathcal{T}^{n, \mathcal{M}_n}$,
		\eref{eq:lev}, Section \ref{sect3b}.
		\item At each time level $t^{n,p}_l, p=0, 1, ..., \mathcal{P}_l-1$,  perform the piecewise polynomial reconstruction
		\eref{pwlinear} and compute the point values, Section
                \ref{sect2}, Section \ref{sect3b}. 
		\item At each time level $t^{n,p}_l, p=0, 1, ..., \mathcal{P}_l-1$,  calculate the one-sided local speeds of
		propagation using \eref{speed}, Section \ref{sect2},
                Section \ref{sect3b}.
		\item At time $t^n$, calculate the reference time step
		$\Delta t$ using \eref{eq:CFLdtref}, Section \ref{sect3b}.
		\item  At each time level $t^{n,p}_l, p=0, 1, ..., \mathcal{P}_l-1$, compute the local time step for
		each cell level, \eref{eq:dtl},  Section \ref{sect3b}.
		\item  At each time level $t^{n,p}_l, p=0, 1, ...,
                  \mathcal{P}_l-1$, compute numerical fluxes and source term in
		the adaptive central-upwind scheme
                \eref{eq:rk1}-\eref{eq:rk2}, \eref{flux},
                \eref{qra1}-\eref{qra2}, Section \ref{sect2},  Section
                \ref{sect3b}.			
	\end{itemize} 
	
	\textbf{Step 2.} On mesh $\mathcal{T}^{n, \mathcal{M}_n}$, compute WLR error using \eref{ej} in Section
        \ref{sect3c} and update
	the refinement/de-refinement status for each cell/triangle,
	Section \ref{sect3c}.
	
	\textbf{Step 3.} Generate the new adaptive mesh $\mathcal{T}^{n+1, \mathcal{M}_{n+1}}$ at  $t^{n+1}$,
	Section \ref{sect3a}. This step includes coarsening of some cells, refinement of
	some cells, and the appropriate projection of
	the cell averages from the mesh $\mathcal{T}^{n,
          \mathcal{M}_n}$ at $t^n$ onto a new adaptive mesh  $\mathcal{T}^{n+1, \mathcal{M}_{n+1}}$
	at time $t^{n+1}$,
	Section \ref{sect3a}.

       \textbf{Step 4.}  Repeat \textbf{Step 1} - \textbf{Step 3}
       until final time.

	\subsection{Adaptive Mesh Refinement/Coarsening}\label{sect3a}
	
	\par
\begin{center}
		\begin{figure}[h!]
			\centering
			\subfigure[Triangle $T^{n,m}_j$ (left) is split into four ``children'' cells $T^{n,m+1}_{js}, s=1,2,3,4$  (right).]{\includegraphics[width=0.8\textwidth]{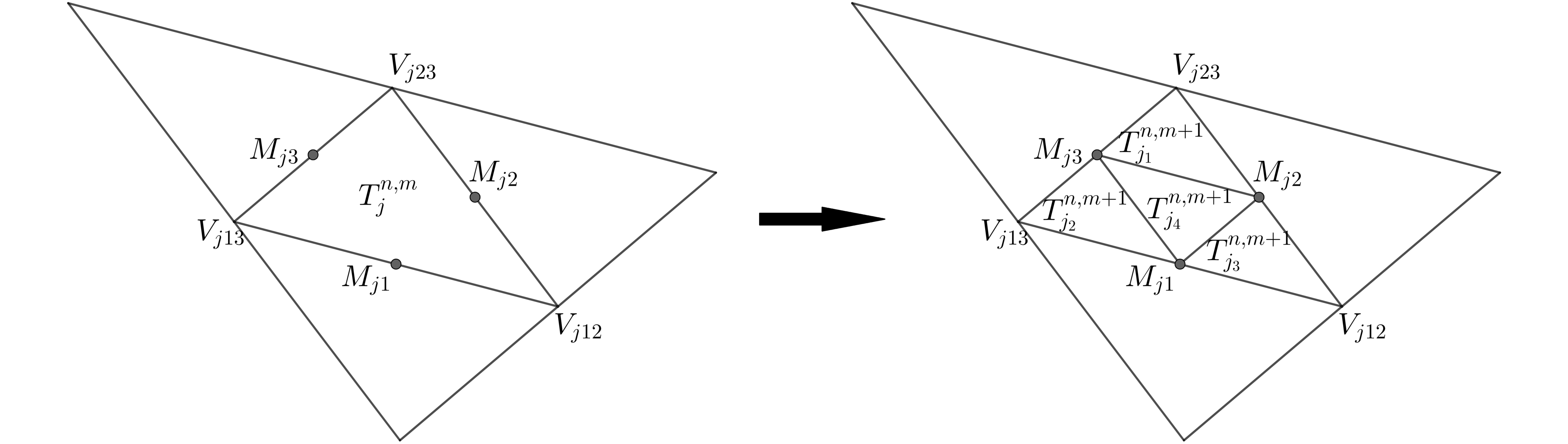}}\label{fig:ref1b}	
	\subfigure[Refinement in the neighboring cells  of $T^{n,m}_j$.]{\includegraphics[width=0.4\textwidth]{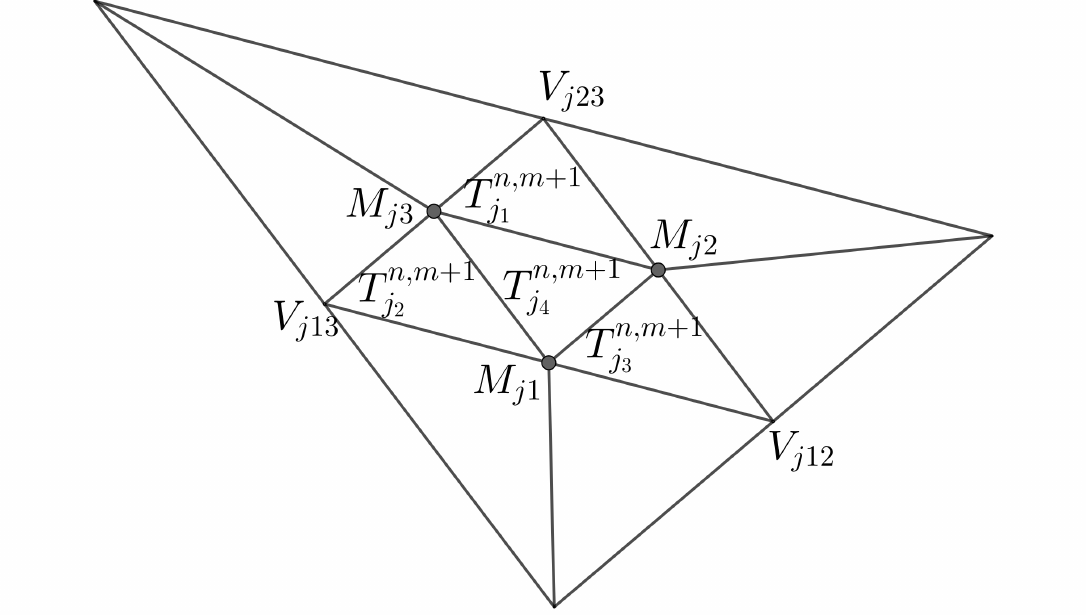}}\label{fig:ref1c}	
			\vspace*{2mm}
			\caption{An outline of the ``regular refinement''.}\label{fig:ref1}
		\end{figure}
	\end{center}

 The main idea of the proposed adaptive mesh refinement algorithm
	is as follows. At time $t^n$, we start with the given mesh, denoted  as
	$\mathcal{T}^{n,m}=\{ T^{n,m}_j\}$, where $T^{n,m}_j$ is a
	triangular cell of size $|T^{n,m}_j|$ with the barycenter $(x^{n, m}_j, y^{n,m}_j)$ within the
	initial	mesh $\mathcal{T}^{n,m}$,  and index $m=0,1,2 ...$ is the level of refinement ($m=0$ corresponds to the mesh with no refinement and $\mathcal{T}^{n,0}\equiv\mathcal{T}^{0,0}$ for all $n$). To flag triangular cells in the mesh
	$\mathcal{T}^{n,m}$ for the refinement/de-refinement (or
        coarsening), we use weak local residual (WLR)
	error estimate,  see Section \ref{sect3c}. We apply ``regular
	refinement'' on the triangles flagged for refinement to obtain a new mesh
	$\mathcal{T}^{n,m+1}$ with the refinement level $m+1$. The ``regular refinement'' on a fully
	flooded triangle is obtained by
	splitting each flagged triangle (``parent'' triangle) into four smaller
	triangles (``children'' triangles) by inserting a new node at the
	mid-point of each edge of the ``parent'' triangle.  We
	illustrate this idea using \fref{fig:ref1} (a), where we show
	an example of splitting a flagged triangle $T^{n,m}_j$ by
	using the mid-points of the sides to obtain the ``children'' cells $T^{n,m+1}_{j_s}, s=1,2,3,4$. In addition, the insertion of
	new nodes on the edges means that non-flagged triangles adjacent to refined
	triangles get hanging nodes and must also be refined. This is done by
	inserting a new edge between the hanging node and the opposite
	corner as illustrated in \fref{fig:ref1} (b). 	
        	\vspace*{-11mm}
        	
	In practice, we may want to reach a higher level of refinement for some cells. This happens when those cells have very large WLR error \eref{ej}, and we need to add more data points. We can obtain a finer cell by repeating the refinement for the flagged triangles in the refined mesh  $\mathcal{T}^{n, m+1}$ to get the mesh with higher level $\mathcal{T}^{n, m+2}, m=0,1,2,....$. \fref{3fig1} is the illustration of the ``regular
	refinement'' procedure with two levels of refinement.	
	\begin{center}\label{3fig1}
		\begin{figure}[h!]
			\centering
			\includegraphics[scale=0.65]{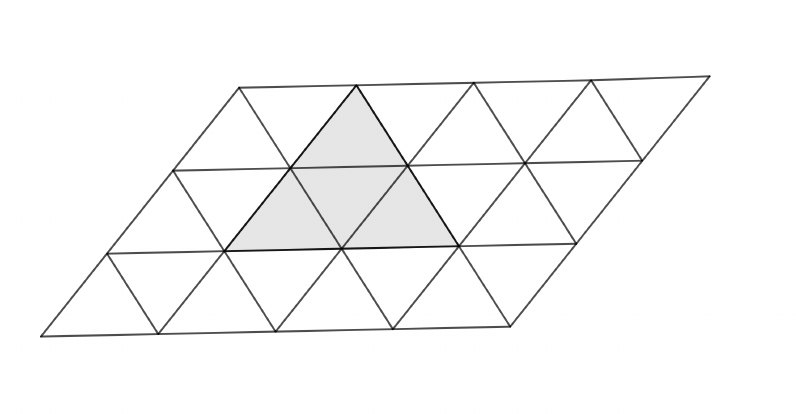}
			\includegraphics[scale=0.65]{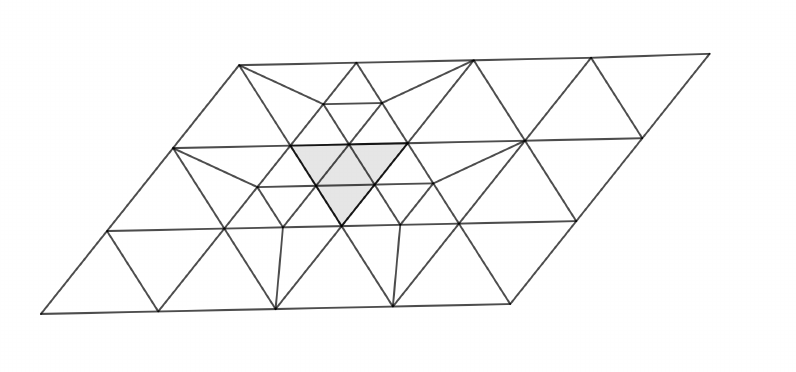}
			\includegraphics[scale=0.65]{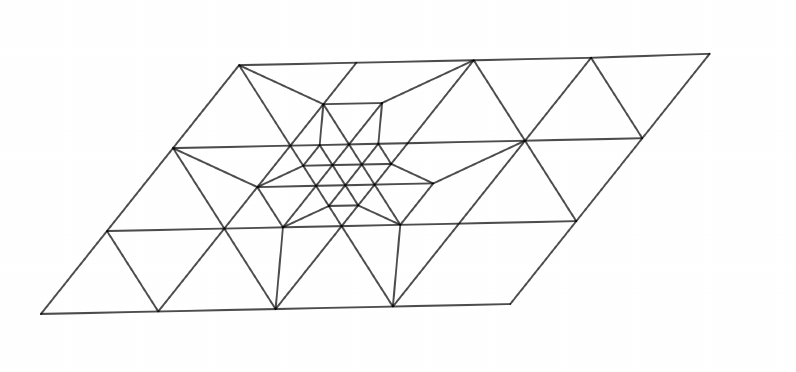}
			\caption{An example of the ``regular
				refinement'' procedure with two levels of refinement. The left figure is the initial coarse mesh
				$\mathcal{T}^{n,0}$ with the region
				flagged for the refinement (gray). The middle figure is the ``first''
				level mesh
				$\mathcal{T}^{n,1}$  with the region
				flagged for higher level of the refinement
				(gray). The right figure is the
				``second'' level mesh $\mathcal{T}^{n,2}$.}\label{fig:mull}
		\end{figure}
	\end{center}		

\begin{center}
		\begin{figure}[h!]
			\centering
			\includegraphics[width=0.4\textwidth]{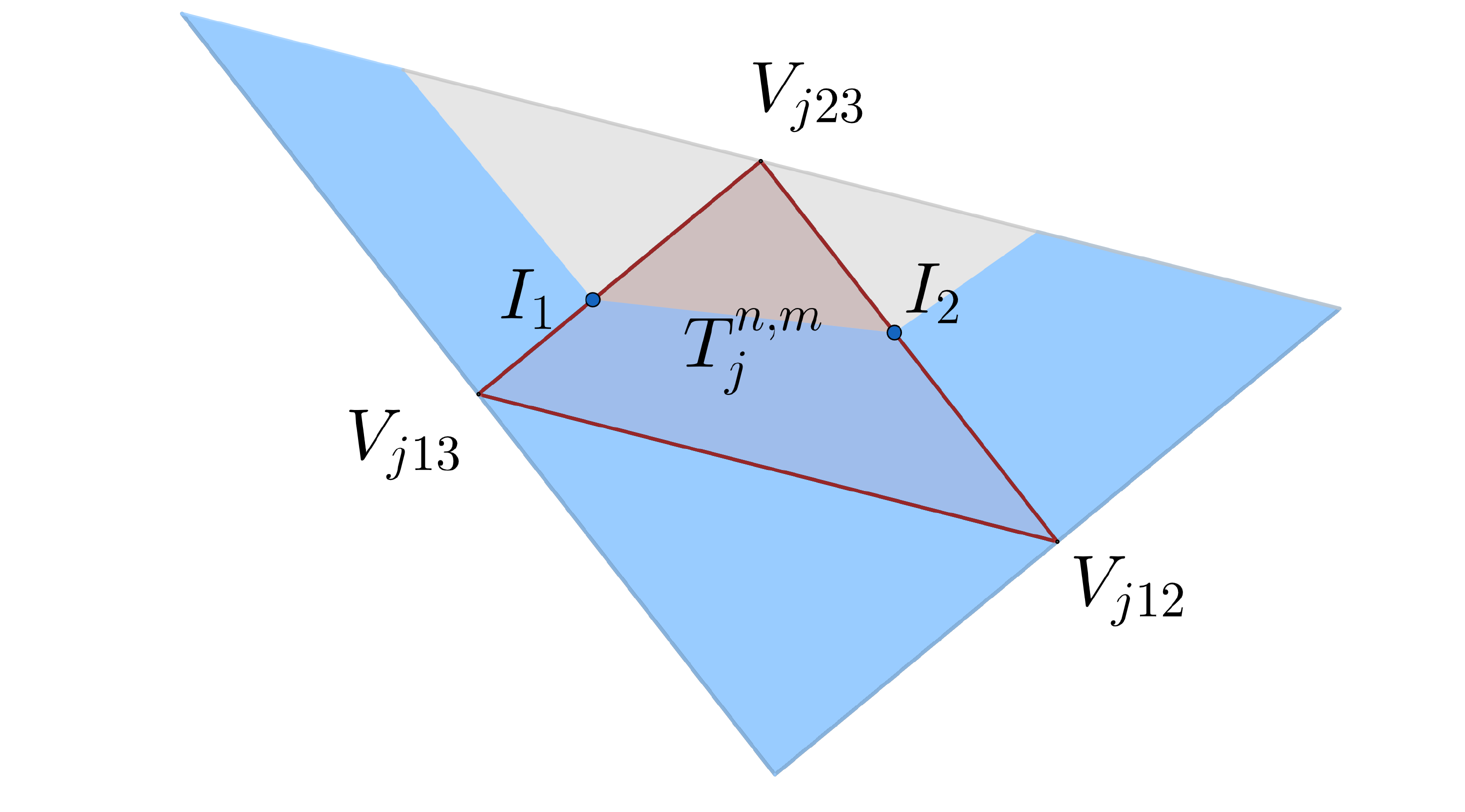}
			\caption{An example of a partially flooded triangle $T^{n,m}_j$ with  wet (blue) and dry (gray) regions, where $T^{n,m}_j$ is flagged for refinement (red) and its neighboring cells are not flagged.}\label{fig:parttri}
		\end{figure}
	\end{center}

\begin{center}
		\begin{figure}[h!]
			\centering
			\includegraphics[width=0.4\textwidth]{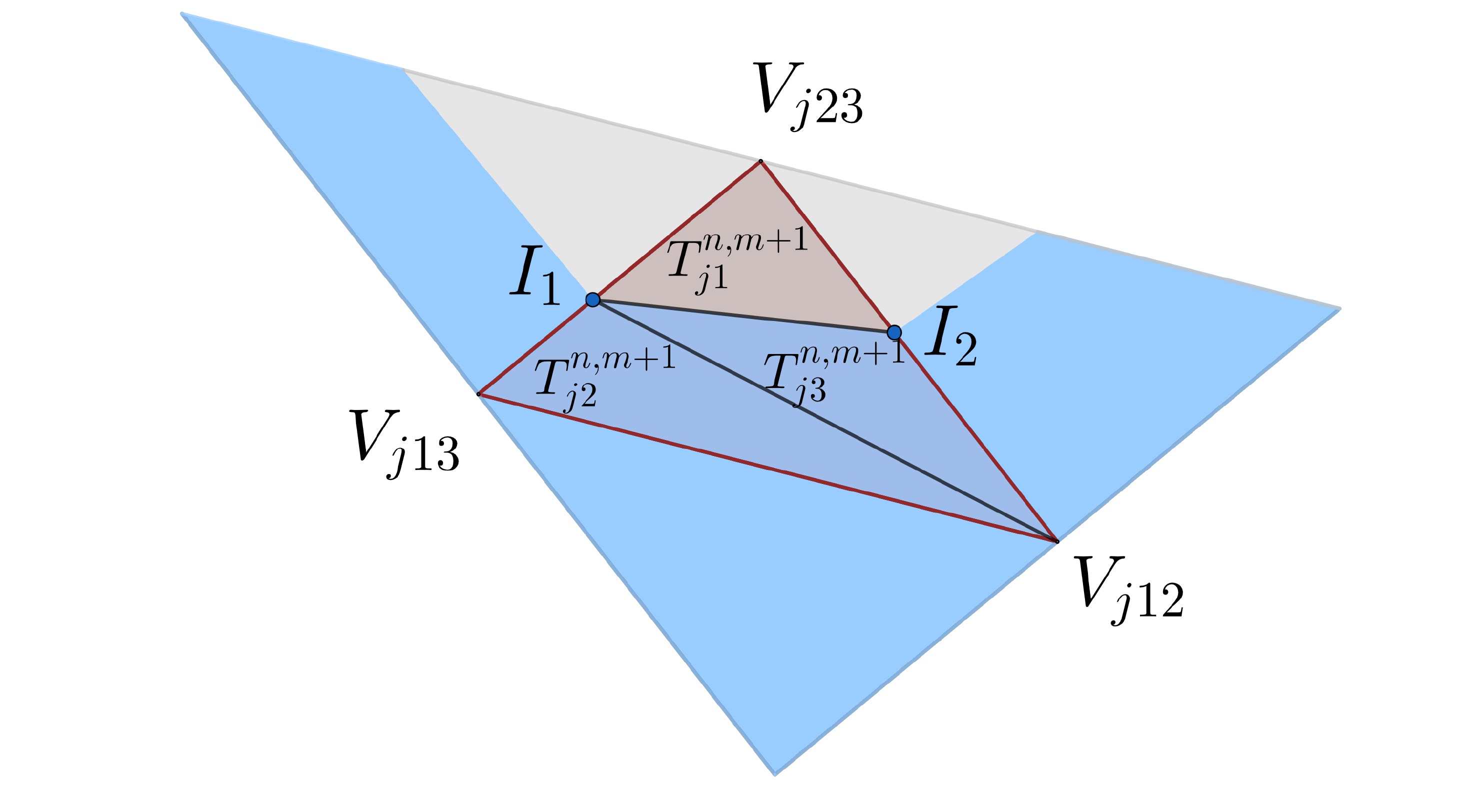}
			\includegraphics[width=0.4\textwidth]{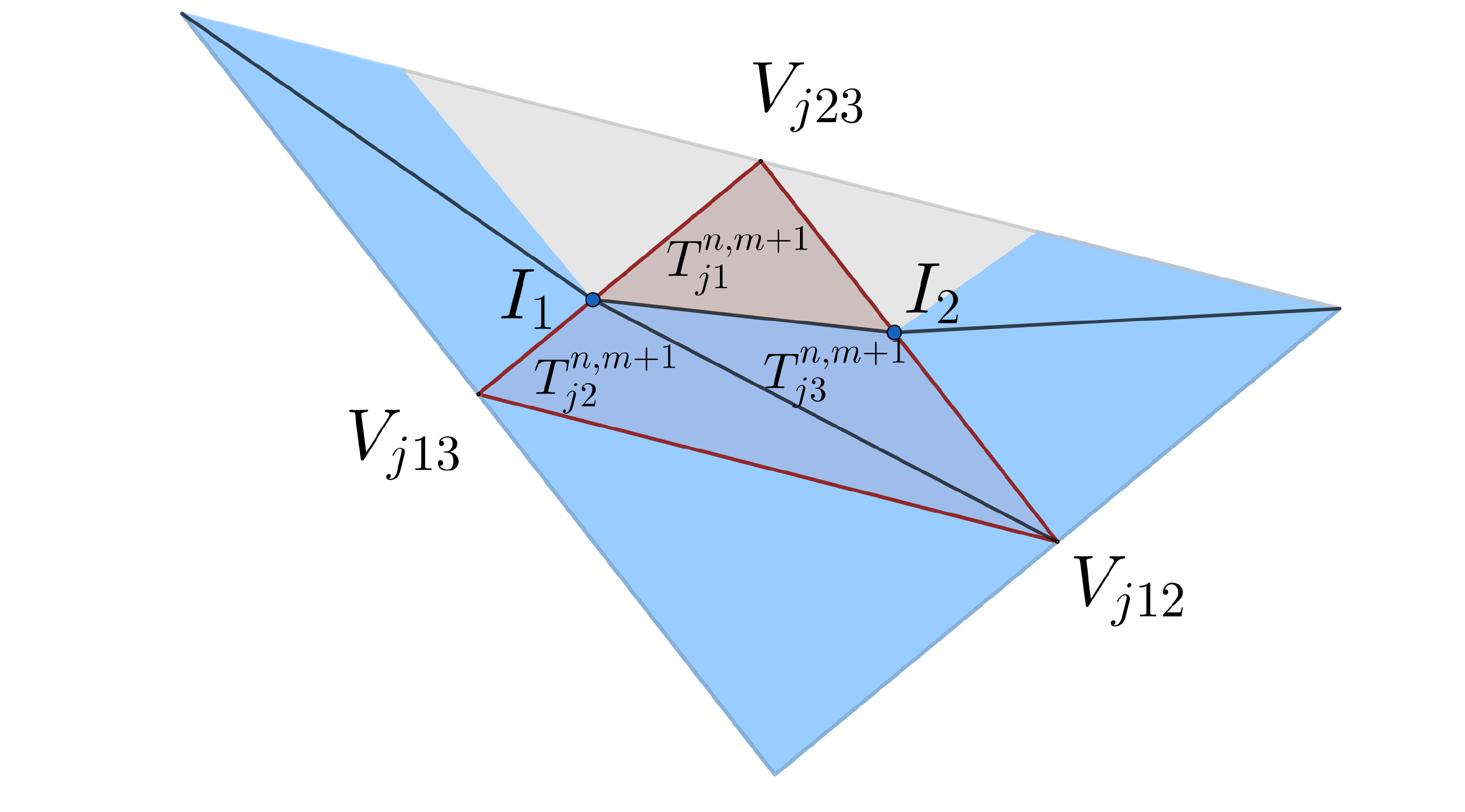}
			\caption{An example of the refinement of a
				partially flooded triangle by using the
				wet/dry interface $I_1I_2$ (left) and the
				refinement of the neighboring triangles (right).}\label{fig:refpart}
		\end{figure}
	\end{center}

\begin{center}
		\begin{figure}[h!]
			\centering
			\includegraphics[width=1\textwidth]{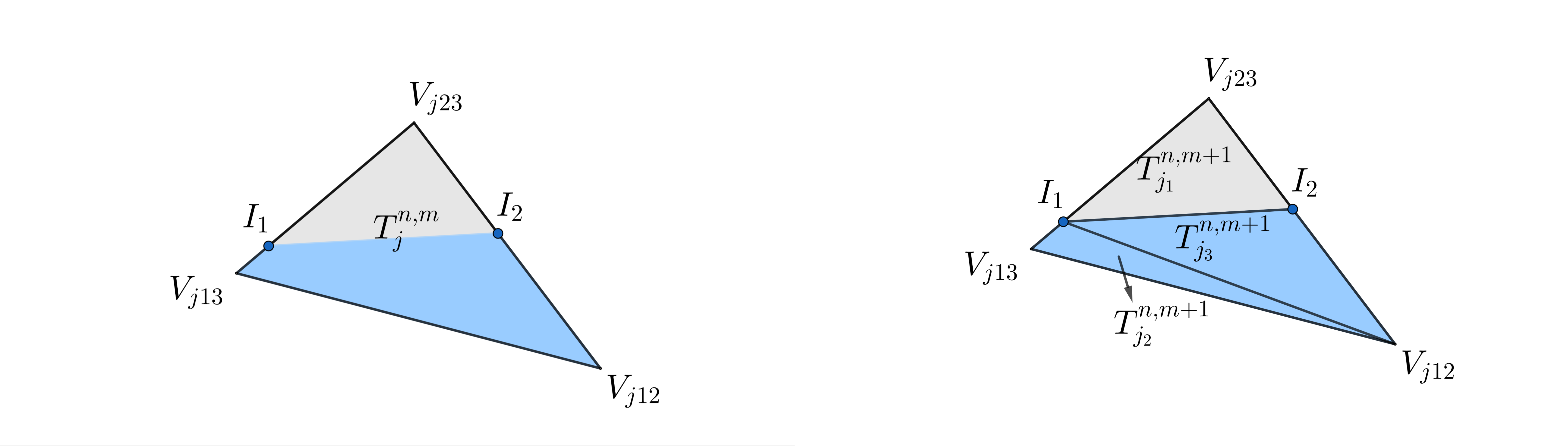}		
			\caption{An example of a refinement of a
				partially flooded triangle (left) using idea of
				``wet/dry refinement'' that produces
				``child'' cell with large obtuse angle
                                (right).}\label{fig:refpartdiff}
		\end{figure}
	\end{center}

      	\vspace*{-11mm}
	In the partially flooded cells, Section \ref{sect2}, the
	approximation of the water surface $\widetilde{w}(x,y)$ at each
	time level $t^n$ consists of two linear pieces, the piece for the wet
	and for the dry region,  see \fref{fig:parttri}. This
	motivates an idea of the ``wet/dry refinement'' which uses the
	boundary between the wet region and the dry region of the cell
	to refine the partially flooded triangles as shown in the
	example in \fref{fig:refpart} (left). Namely, consider a
	partially flooded triangle $T^{n,m}_j$ which is flagged for
	the refinement and has three non-flagged neighboring cells in the grid $\mathcal{T}^{n,m}$. The segment $I_1I_2$ is the
	boundary between the wet and dry interface in that
	triangle. Note that, the location of the nodes $I_1$ and $I_2$
	is determined by the second order water surface reconstruction
	developed in \citep{LAEK}, see also Section \ref{sect2}. During
	``regular refinement'' of the partially flooded cell,  we
	first split the flagged cell $T^{n,m}_j$ into a smaller
	triangle and a quadrilateral using the wet/dry interface
	$I_1I_2$. We then continue to refine the quadrilateral by its
	diagonal. As can be seen in \fref{fig:refpart} (left), the
	flagged triangle $T^{n,m}_j$ has three ``children'' which are
	either fully flooded or dry. Similarly, as in \fref{fig:ref1},
	the appearance of two hanging nodes  $I_1$ and $I_2$ on two
	sides leads to the need of the further splitting of the
	neighboring cells as presented in
	\fref{fig:refpart} (right). The ``regular wet/dry refinement''
	will capture the features of the wet/dry fronts and will
        minimize number of partially flooded ``children'' cells. However, this method may
	give us difficulties in controlling the shape of triangles in
	the adaptive mesh. In some cases, it may produce ``children'' cells with
	unexpected large obtuse angles or very small altitudes,  as shown in
	\fref{fig:refpartdiff}. For cells where such situation happens, we instead use the
	``regular refinement''  for the
	fully flooded cells as described above.
	
      	\vspace*{-30mm}
      	
	Very often in the numerical simulations of the wave phenomena,
	the regions of the domain that need to be refined move with
	time. Hence, the refinement in some cells may become no longer
	needed. The de-refinement or coarsening procedure is then
	introduced to deactivate unnecessarily fine cells in the
	mesh. The de-refinement is performed by coarsening
	(by deactivating ``children'' cells) in the triangles
	of the mesh flagged for coarsening (and possibly deactivating finer neighboring
	triangles due to removal of the hanging nodes). At time $t^n$,
	``children'' cells in the mesh $\mathcal{T}^{n,m+1}, m=0,1,...,M_n-1$, are
	deactivated based on the WLR and the corresponding ``parent'' cell from the mesh
	$\mathcal{T}^{n,m}$ is activated back. 
	In order to minimize the complexity of the adaptive grid generation,
	the de-refinement should be applied on all cells flagged for
	coarsening prior to the refinement, see for example \citep{NHA}.

	The refinement/de-refinement process at time $t^n$ produces a
	hierarchical system of grids
	$\mathcal{S}^n=\{\mathcal{T}^{n,0}, \mathcal{T}^{n,1},
	\mathcal{T}^{n,2},..., \mathcal{T}^{n,\mathcal{M}_n}\}$, where
	$\mathcal{T}^{n,m}$, $m=1,2,...,\mathcal{M}_n$ is the grid
	with the level of refinement $m$ obtained by refining the grid
	$\mathcal{T}^{n,m-1}$. The final mesh $ \mathcal{T}^{n, \mathcal{M}_n} \in
	\mathcal{S}^n$ is the mesh that is used in the
	adaptive central-upwind scheme \eref{eq:rk1}-\eref{eq:rk2} at time level $t^n$. After the evolution of the numerical solution
	from time $t^{n}$ to time $t^{n+1}$ using mesh  $ \mathcal{T}^{n, \mathcal{M}_n} \in
	\mathcal{S}^n$, we proceed with the generation of a new
	adaptive grid $\mathcal{T}^{n+1, \mathcal{M}_{n+1}} \in
	\mathcal{S}^{n+1}$ from the mesh $\mathcal{T}^{n,
		\mathcal{M}_n}$,  using WLR in Section \ref{sect3c}. After a new adaptive mesh $\mathcal{T}^{n+1, \mathcal{M}_{n+1}}$ is constructed, the obtained cell
	averages $\xbar{\mU}^{n+1}$ on the mesh $ \mathcal{T}^{n,
		\mathcal{M}_n}$ need to be projected accurately on the new
	mesh $\mathcal{T}^{n+1, \mathcal{M}_{n+1}}$, using the ideas
	as summarized briefly below.
	
	\textbf{Case 1.} If a triangle $T^{n+1,
		\mathcal{M}_{n+1}}_j\in \mathcal{T}^{n+1,
		\mathcal{M}_{n+1}}$ at $t^{n+1}$ is  the same cell as in the grid
	$\mathcal{T}^{n,
		\mathcal{M}_{n}}$, we will keep without any change the cell
	averages for that triangle at $t^{n+1}$.
	
	\textbf{Case 2.} A cell $T^{n+1, \mathcal{M}_{n+1}}_j\in
	\mathcal{T}^{n+1, \mathcal{M}_{n+1}}$ is obtained by
	coarsening some finer cells $T^{n, \mathcal{M}_n}_{j_s} \in
	\mathcal{T}^{n, \mathcal{M}_n}, s=1,2,.., S$. In order to
	enforce the conservation, the solution, $\xbar{U}^{n+1}_j$ in
	the cell $T^{n+1, \mathcal{M}_{n+1}}_j$,  is computed as
	$$\begin{aligned}
	\xbar{U}^{n+1}_j&=\dfrac{1}{|T^{n+1,\mathcal{M}_{n+1}}_j|}\sum \limits_{s=1}^{S}\xbar{U}^{n}_{j_s}|T^{n,\mathcal{M}_n}_{j_s}|, 
	\end{aligned}$$
	where $\xbar{U}^{n}_{j_s}$ is the solution in $T^{n,\mathcal{M}_n}_{j_s}$.
	
	\textbf{Case 3.} A triangle $T^{n+1,
		\mathcal{M}_{n+1}}_j\in \mathcal{T}^{n+1,
		\mathcal{M}_{n+1}}$ is obtained from the refinement of the cell  $T^{n,\mathcal{M}_n}_i\in \mathcal{T}^{n,\mathcal{M}_n}$. The
	approximation of the
	cell averages of the solution at $t^{n+1}$ in $\mathcal{T}^{n+1,
		\mathcal{M}_{n+1}}$ is obtained using the
	evaluation of  the piecewise linear reconstruction
        \eref{pwlinear} of the
	solution at $t^{n+1}$ in the triangle
	$T^{n,\mathcal{M}_n}_i$. 	
		\subsection{Second-order Adaptive Time Evolution}\label{sect3b}
	The CFL-type condition \eref{1eqCFL} is needed due to numerical
        stability. Hence, use of a global time step in the adaptive
        algorithm may lead to 
	the time step  that can become very small due presence of much
        finer cells in the mesh. To improve the computational cost of
        the algorithm, we consider the approach based on the adaptive
        time step from \citep{domingues2008adaptive, Donat, MORE}. The main idea of the adaptive time evolution is that we group cells
	into different levels based on the cell sizes. After that, we evolve the
        solution on each cell level individually with its local time
        step. This approach does not violate the stability of the
        explicit time discretization scheme as was shown in
        \citep{domingues2008adaptive, MORE}. Below, we present the
        brief summary of the adaptive time evolution algorithm based
        on  the second-order strong stability preserving Runge-Kutta
        methods (SSPRK2) in \citep{GST},
        \citep{domingues2008adaptive, Donat, MORE}.
\begin{center}
		\begin{figure}[ht!]
			\centering
			\includegraphics[width=0.7\textwidth]{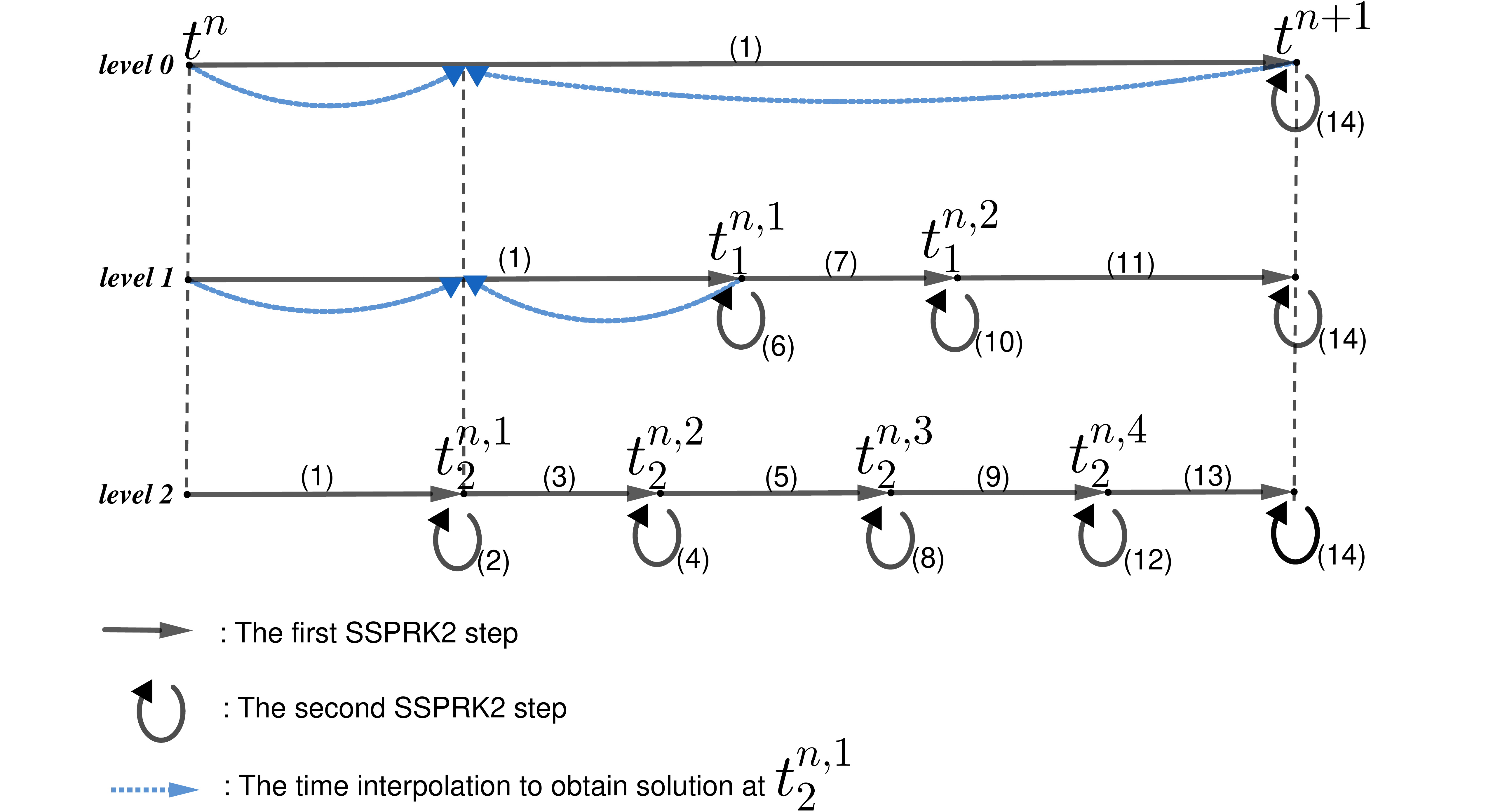}	
\caption{The example of SSPRK2 on mesh with
                          three cell levels,  $l=0,1,2$.}\label{fig:order}
		\end{figure}
	\end{center}
        	\vspace*{-8mm}
        	
  The idea and the order of steps of adaptive SSPRK2 are illustrated on
 the example in \fref{fig:order}.    First, we group all cells
 $T^{n,\mathcal{M}_n}_j\in\mathcal{T}^{n,\mathcal{M}_n}$ at time $t^n$
 in cell levels $l=0, 1,
 .., L$ based on their sizes. We define cell levels at $t^n$ in mesh
$\mathcal{T}^{n,\mathcal{M}_n}$ as follows.  A cell
$T^{n,\mathcal{M}_n}_j\in\mathcal{T}^{n,\mathcal{M}_n}$ belongs to the
level $l$, if $l$ is the smallest positive integer satisfying,	
	\begin{equation}
	2^l\geq \dfrac{\displaystyle  \max_j\left(\displaystyle
            \min_{k}\left(r_{jk}\right)\right)}{\displaystyle\min_{k}(r_{j,k})}, \label{eq:lev}
	\end{equation}	
	where $r_{jk}, k=1,2,3$ are three altitudes of triangle
        $T^{n,\mathcal{M}_n}_j$. Thus, the cell levels with larger index $l$
        will contain finer cells which will require a
        smaller time step \eref{1eqCFL} for the evolution from  $t^n
        \to t^{n+1}$. Next, at $t^n$, we define the
        reference time step $\Delta t$ as the local time step on the
        coarsest level $l=0$  of cells in the mesh
        $\mathcal{T}^{n,\mathcal{M}_n}$ by considering the CFL-type
        condition \eref{1eqCFL} locally on level $l=0$. Namely, define
        first $a_{max}$,
	\begin{equation}
	a_{max}:= \displaystyle \max_{j,k}(a^{\rm in}_{jk},a^{\rm out}_{jk}),\label{eq:amax}
	\end{equation}	
	where $(a^{\rm in}_{jk},a^{\rm out}_{jk})$ are the local
        one-sided speeds of propagation \eref{speed} at $t^n$ for
        sides $k=1,2,3$ in the triangle
        $T^{n,\mathcal{M}_n}_j\in \mathcal{T}^{n,\mathcal{M}_n}
        $. Then, the
        reference time step $\Delta t$ is computed as,	
	\begin{equation}
	\begin{aligned}
	\Delta t\equiv \Delta t^{n,0}_0 =\frac{0.9\displaystyle \max_j\left(\displaystyle \min_{k}\left(r_{jk}\right)\right)}{6a_{max}}.
	\end{aligned}
	\label{eq:CFLdtref}
	\end{equation}
 We set, $t^{n+1}=t^n+\Delta t$.

Assume next, that $\mathcal{P}_l$ is the
 number of steps taken on higher levels $l=1, .. L$ to evolve from $t^n$ to
 $t^{n+1},$ namely $[t^n, t^{n+1}]=\cup[t_l^{n,
   p}, t_l^{n, p+1}], p=0, ..., \mathcal{P}_l-1$ with $t^{n,0}_l\equiv
        t^n$, $t_l^{n, \mathcal{P}_l}\equiv t^{n+1} \quad \forall l$.  We define the local time step for cells on these levels $l=1, ..., L$ as,
	\begin{equation}\Delta t^{n, p}_l=\dfrac{2^{-l} \Delta
            t}{\max(\mu^{n, p}_l,1)},\label{eq:dtl}
	\end{equation}
where parameter $\mu^{n,p}_l$ takes into account change in the local one-sided
speeds of the propagation, 
 $$\mu^{n, p}_l=\dfrac{	\displaystyle\max_{j,k}(a^{\rm in}_{jk},a^{\rm
       out}_{jk})^{n, p}_l}{a_{max}},$$ where $(a^{\rm
   in}_{jk},a^{\rm out}_{jk})^{n, p}_l$ are the local one-sided
 speeds of propagation at $t^{n, p}_l$ of the cell
 $T^{n,\mathcal{M}_n}_j$ in the level $l$. Therefore, on
 each cell $T^{n,\mathcal{M}_n}_j$ of level $l$,  for each substep $[t_l^{n,p},t_l^{n, p+1}]\equiv [t_l^{n,p},t_l^{n, p}+\Delta
        t_l^{n, p}],
 p=0,1,2,3,...,\mathcal{P}_l-1$ of the evolution from $t^n$ to  $t^{n+1}$, we apply the following two adaptive steps of SSPRK2 method,
	\begin{subequations}
		\begin{align}
		\xbar{\bm{U}}_j^{(1)}&=\xbar{\bm{U}}_j^{n,
                                       p}-\frac{1}{|T^{n,\mathcal{M}_n}_j|}\sum_{k=1}^3\Delta
                                       t^{n, p}_{jk}{\bm
                                       H}^{n,p}_{jk}+\Delta
                                       t_l^{n,p}\,\xbar{\bm{S}}^{n,
                                       p}_j:={\bm{R}}(\xbar{\bm{U}}_j^{n,
                                       p},\Delta t_l^{n, p}),\label{eq:rk1}\\
		\xbar{\bm{U}}_j^{n,
                  p+1}&=\dfrac{1}{2}\xbar{\bm{U}}_j^{n,
                        p}+\dfrac{1}{2}{\bm{R}}(\xbar{\bm{U}}_j^{(1)},\Delta
                        t_l^{n, p}),\label{eq:rk2}
		\end{align}
	\end{subequations}
	Here, $\Delta
        t_l^{n, p}$ is the local time step of the cells of
        level $l$ at time $t_l^{n, p}$ \eref{eq:dtl},  and $\Delta
        t^{n, p}_{jk}$ is the ``draining'' time step in cell
        $T^{n,\mathcal{M}_n}_j$ for the local time step
        $\Delta t_l^{n, p}$ in level $l$.  The flux term ${\bm
                                       H}^{n,p}_{jk}$ in \eref{eq:rk1}
                                     -\eref{eq:rk2} is the flux
                                     \eref{flux} computed at $t=t^{n,
                                       p}_l$. The source term $\xbar{{\bm
                                       S}}^{n,p}_{j}$ in \eref{eq:rk1}
                                     -\eref{eq:rk2} is the source
                                     \eref{qra1} - \eref{qra2}
                                     computed at $t=t^{n, p}_l$ with
                                     the time step $\Delta t_l^{n, p}$ and
                                     with the corresponding local
                                     ``draining'' time step  $\Delta
        t^{n, p}_{jk}$.
Note that,
        $\xbar{\bm{U}}_j^{n, 0}\equiv
        \xbar{\bm{U}}_j^{n}$ and 	$\xbar{\bm{U}}_j^{n,
          \mathcal{P}_l}\equiv \xbar{\bm{U}}_j^{n+1}$.
\begin{remark}  If cells from different cell levels are neighbors, we
  use linear interpolation in time to match the time levels of such
  cells, see also Fig. 3.6, for the illustration of the interpolation. 
\end{remark}	
	\subsection{A Posteriori Error Estimator}\label{sect3c}
	\par Here, using the idea of Weak Local Residual (WLR) from
        \citep{MR1111445, MR2126235},
	we will derive local error estimator that is used as the
        robust indicator for the adaptive mesh refinement in our work.	
		\par Let us recall that the weak form of the mass conservation
	equation for the system \eref{eqs2} in $\Omega\times [0, T]$ takes the integral form,
	\begin{equation}
	\int_0^T \int_{\Omega} (w \phi_t(x, y,
	t)+hu\phi_x(x, y, t)+hv\phi_y(x, y, t) d\Omega dt+\int_{\Omega} w(x,
	y, 0)\phi(x, y, 0) d\Omega=0,
	\label{Errapp}
	\end{equation}
	for all sufficiently smooth test functions $\phi(x,y,t)$ with
	compact support on $\Omega\times[0,T)$.
	
	\par	Consider example of localized test function in time and space,
	
	\begin{equation}
	\begin{aligned}
	\phi^{n+\frac1{2}}_{i}(x,y,t)=\frac 1{\Delta}f_{i}(x,y)f^{n+\frac1{2}}(t),
	\end{aligned}
	\label{Errtest}
	\end{equation}

	where $\Delta:=\max(\displaystyle\max_{j, k}(r_{jk}), \Delta t)$,
	$\Delta t=t^{n+1}-t^{n}=t^{n+\frac 1{2}}-t^{n-\frac 1{2}},
	n=1,2,3...$, and $r_{jk}$ are the heights of the triangle/cell $T^{n,\mathcal{M}_n}_j\in\mathcal{T}^{n,\mathcal{M}_n}$ . Function $f^{n+\frac1{2}}(t)$ is a linear function in time
	with a local support defined as,	
	\begin{equation}
	\begin{aligned}
	f^{n+\frac1{2}}(t)&=\begin{cases}\dfrac{t-t^{n-\frac1{2}}}{\Delta t}, &\mbox{if} \quad t^{n-\frac1{2}}\leq t\leq  t^{n+\frac1{2}},\\\\\dfrac{t^{n+\frac3{2}}-t}{\Delta t}, &\mbox{if} \quad t^{n+\frac1{2}}\leq t\leq  t^{n+\frac3{2}},\\\\0, &\mbox{otherwise}.\end{cases}
	\end{aligned}
	\label{Errtest1}
	\end{equation}

	Function $f_{i}(x, y), i=1,2,3,...$ is a ``hat function'',
        namely, a
	piecewise linear function  with compact support over all
	triangles with common  vertex $N_{i}=(\widetilde{x}_i,\widetilde{y}_i)$. The function
	$f_i(x, y)$ takes value $1$ at the vertex $N_i$ and $0$ at all
	other nodes. More precisely, assume that there are $C_i$
	triangles $T^{n,\mathcal{M}_n}_{j_1}$, $T^{n,\mathcal{M}_n}_{j_2}$, $T^{n,\mathcal{M}_n}_{j_3}$,..., $T^{n,\mathcal{M}_n}_{j_{C_i}}\in \mathcal{T}^{n, \mathcal{M}_{n}}$
	which share common vertex $N_{i}$. Thus, the function $f_{i}(x,
	y)$ is defined as,
	\begin{equation}
	\begin{aligned}
	f_{i}(x,y)&=\begin{cases}a^{(i)}_c(x-\widetilde{x}_i)+b^{(i)}_c(y-\widetilde{y}_i)+1,&\mbox{if}\quad (x,y)\in T^{n,\mathcal{M}_n}_{j_c}, \quad c=1,2,...C_i\\0, &\mbox{otherwise}\end{cases},\\
	\end{aligned}
	\label{Errtest2}
	\end{equation}
	The quantity $(a^{(i)}_c, b^{(i)}_c)$ is the gradient of the
	linear piece of $f_i(x, y)$ restricted to $T^{n,\mathcal{M}_n}_{j_c}$,	
	\begin{equation}
	\begin{aligned}
	a^{(i)}_c&=\dfrac{\widetilde{y}_2-\widetilde{y}_3}{(\widetilde{y}_3-\widetilde{y}_i)(\widetilde{x}_2-\widetilde{x}_i)-(\widetilde{y}_2-\widetilde{y}_i)(\widetilde{x}_3-\widetilde{x}_i)},\\
	b^{(i)}_c&=\dfrac{\widetilde{x}_3-\widetilde{x}_2}{(\widetilde{y}_3-\widetilde{y}_i)(\widetilde{x}_2-\widetilde{x}_i)-(\widetilde{y}_2-\widetilde{y}_i)(\widetilde{x}_3-\widetilde{x}_i)},\\
	\end{aligned}\label{gradfor}
	\end{equation}
	where $N_i=(\widetilde{x}_i, \widetilde{y}_i), (\widetilde{x}_2, \widetilde{y}_2), $ and $(\widetilde{x}_3, \widetilde{y}_3) $ are the three vertices of triangle $T^{n,\mathcal{M}_n}_{j_c}$.
	Next, define the following piecewise constant approximation for
	the solution $\mU=(w,hu,hv)$,
	
	\begin{equation}
	\mU^\Delta:=\xbar{\mU}_{j_c}^{n} , \mbox{ if }  (x,y,t)\in
	T^{n,\mathcal{M}_n}_{j_c}\times[
        t^{n-\frac1{2}},t^{n+\frac1{2}}]. 
	\label{Errapps1}
	\end{equation}
	
	\par Now, using the localized test function
	$\phi^{n+\frac1{2}}_{i}(x,y,t)$, \eref{Errtest} together
	with the piecewise constant approximation $\mU^\Delta$, \eref{Errapps1} in
	\eref{Errapp}, we can define {\it the weak form of the truncation error, (WLR)},
	which will be used as the error indicator for refinement/de-refinement in the adaptive grid,

	\begin{equation}
	\begin{aligned}
	E^{n+\frac1{2}}_{i}:=E(\mU^\Delta, \phi^{n+\frac{1}{2}}_i)&=\sum_{c=1}^{C_i}
	\int_{t^{n-\frac{1}{2}}}^{t^{n+\frac{1}{2}}}
        \int_{T^{n,  \mathcal{M}_n}_{j_c}}
	(\xbar{w}^{n}_{j_c}(\phi^{n+\frac1{2}}_{i})_t+(\xbar{hu})^{n}_{j_c}(\phi^{n+\frac1{2}}_{i})_x+(\xbar{hv})^{n}_{j_c}(\phi^{n+\frac1{2}}_{i})_y)d\Omega
	dt \\
	&+ \sum_{c=1}^{C_i}
	\int_{t^{n+\frac{1}{2}}}^{t^{n+\frac{3}{2}}} \int_{T^{n,\mathcal{M}_n}_{j_c}}
	(\xbar{w}^{n+1}_{j_c}(\phi^{n+\frac1{2}}_{i})_t+(\xbar{hu})^{n+1}_{j_c}(\phi^{n+\frac1{2}}_{i})_x+(\xbar{hv})^{n+1}_{j_c}(\phi^{n+\frac1{2}}_{i})_y)d\Omega
	dt.
	\end{aligned}
	\label{Errapp1}
	\end{equation}

	\par After straightforward calculations, the WLR error
	$E^{n+\frac1{2}}_{i}$ on mesh $\mathcal{T}^{n, \mathcal{M}_{n}}$ is given by the formula,
	
	\begin{equation}\label{wlrerr}
	\begin{aligned}
	E^{n+\frac 1{2}}_{i}&=\frac1{\Delta}(\mathcal{U}^{n+\frac 1{2}}_{i}+\mathcal{F}^{n+\frac 1{2}}_{i}+\mathcal{G}^{n+\frac 1{2}}_{i}),\\\mathcal{ \mU}^{n+\frac 1{2}}_{i}&=\sum \limits_{c=1}^{C_i} \frac1{3}|T^{n, \mathcal{M}_n}_{j_c}| (\xbar{w}^{n}_{j_c}-\xbar{w}^{n+1}_{j_c}),\\
	\mathcal{ F}^{n+\frac 1{2}}_{i}&=\sum \limits_{c=1}^{C_i}
        a^{(i)}_c\frac{\Delta t}{2}|T^{n, \mathcal{M}_n}_{j_c}| ((\xbar{hu})^{n}_{j_c}+(\xbar{hu})^{n+1}_{j_c}),\\
	\mathcal{ G}^{n+\frac 1{2}}_{i}&=\sum \limits_{c=1}^{C_i}
        b^{(i)}_c\frac{\Delta t}{2}|T^{n, \mathcal{M}_n}_{j_c}| ((\xbar{hv})^{n}_{j_c}+(\xbar{hv})^{n+1}_{j_c}).
	\end{aligned}
	\end{equation}

	The error in a cell $T^{n,\mathcal{M}_n}_j\in \mathcal{T}^{n.\mathcal{M}_n}$ is given by, 
	\begin{equation}
	e_j=\max_{k}\left|E^{n+\frac 1{2}}_{jk}\right|, \quad
        k=1, 2, 3, \label{ej}
	\end{equation}
	where $E^{n+\frac 1{2}}_{jk}$ is the  WLR error
        computed in (\ref{wlrerr}) at a node $k$ of triangle $T_j$.

In our numerical experiments, we define an error tolerance,  $\omega$
as,   
\begin{equation} \label{errtol}
\omega=\sigma \max_j(e_j),
\end{equation}
 where $\sigma<1 $ is a given problem-dependent constant (see Section
 \ref{sect4}),  and $e_j$ is
 the WLR error in the triangle $T^{n.\mathcal{M}_n}_j$, \eref{ej}. The
 error $e_j$ in each cell $T^{n.\mathcal{M}_n}_j \in \mathcal{T}^{n,
   \mathcal{M}_{n}}$ is compared to the error tolerance
 \eref{errtol}, and the cell
 is either ``flagged'' for refinement/de-refinement or ``no-change''. 

 Note that, in this work, we consider only the equation
   \eref{eq:swe1} to obtain WLR error.  The full system of shallow
   water equations can be used to derive WLR too, however it  will make the computation
   of the error indicator more complex and more expensive. 
	\section{Numerical Examples}\label{sect4}

\begin{figure}[h!]	
		\centering
		\includegraphics[width=0.3\linewidth]{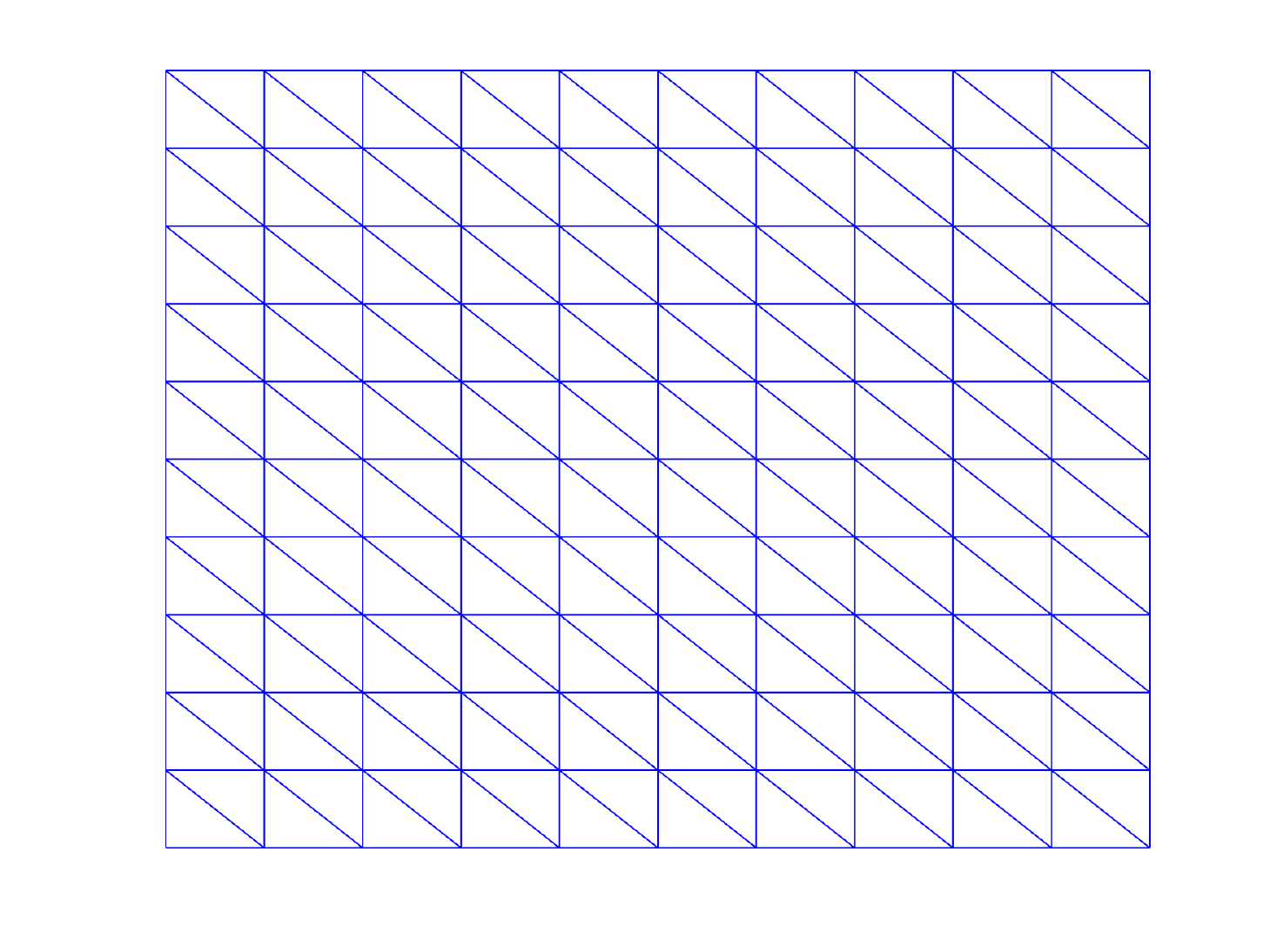}\\
		\vspace*{2mm}
		\caption{An outline of uniform triangular mesh.}
		\label{fig:ex1}
	\end{figure}

	In this section, we illustrate the performance of the
	designed adaptive central-upwind scheme. We compare the
	results of the adaptive central-upwind scheme developed in this
	work with the results
	of the central-upwind scheme from \citep{LAEK} on uniform
	triangular meshes (example of such uniform triangular mesh is
	outlined  in \fref{fig:ex1}). In addition, in all experiments,
	we compute ratio,
	$\mathcal{R}_{CPU}=\frac{CPU_{uniform}}{CPU_{adaptive}}$,
	which is the ratio of the CPU times of the central-upwind
	algorithm without adaptivity to the CPU time of the adaptive
	central-upwind algorithm. To  compare $L^1$-errors, as well as
        to compare
        the CPU times and to
	compute $\mathcal{R}_{CPU}$, we consider uniform mesh and the
	adaptive mesh with the same size of the smallest cells. More
        precisely, in Table \ref{tab:ex1a},  $L^1$-errors,  and in Tables
        \ref{tab:ex1b}-\ref{tab:ex3ccpu}, $\mathcal{R}_{CPU}$ are computed using the uniform meshes $2\times N\times N, N=100, 200, 400$ and using the corresponding adaptive meshes which are obtained from the coarser uniform mesh  $2\times N/2^\mathcal{M}\times N/2^\mathcal{M}$ ($\mathcal{M}=1,2$ is the highest level of refinement in the adaptive mesh). In Example \ref{sec4_1} and the first
	two cases in Example \ref{sec4_2} with \eref{eq:ex2B}-\eref{eq:ex2I}, the gravitational acceleration
	is set, $g = 1.0$, while in the last case of Example
	\ref{sec4_2} with \eref{eq:ex2cB}-\eref{eq:ex2cI} and
	Example \ref{sec4_3}, we take $g=9.8$. We set the
	desingularization parameters $\tau$ and $\varepsilon$ for
	calculations of the velocity components $u$ and $v$ to be
	$\tau=\max_j\{|T^{n,\mathcal{M}_n}_j|^2\}$ and $\varepsilon=10^{-4}$, except for
	the Example \ref{sec4_2}, \eref{eq:ex2cB}-\eref{eq:ex2cI},
	where $\varepsilon=10^{-2}$ (see Section 2.1 formula (2.6) in \citep{LAEK}).
	
	\subsection{Example 1---Accuracy  Test}\label{sec4_1}
	
	Here, we consider the example from \citep{MR2804645}, and we verify experimentally the order of accuracy of  the designed
	adaptive central-upwind scheme. We also compare the
	computational efficiency of the adaptive central-upwind scheme
	and the central-upwind scheme without adaptivity on uniform and non-uniform triangular meshes. 
	
	We consider the following initial data
	and the bottom topography,
	
	\begin{equation*}
	w(x,y,0)=1, \quad u(x,y,0)=0.3, \quad v(x,y,0)=0,
	\end{equation*}
	
	\begin{equation*}
	B(x,y)=0.5\exp(-25(x-1)^2-50(y-0.5)^2).
	\end{equation*}
	
	The computational domain is $[0,2]\times[0,1]$. A zero-order
	extrapolation is used at all boundaries. The error tolerance
	\eref{errtol} for the
	mesh refinement in this example is set to $\omega = 0.01\max_j(e_j)$.

\begin{figure}[h!]
		\centering
		\subfigure[Uniform mesh $2\times 25\times 25$.]{\includegraphics[width=1\textwidth]{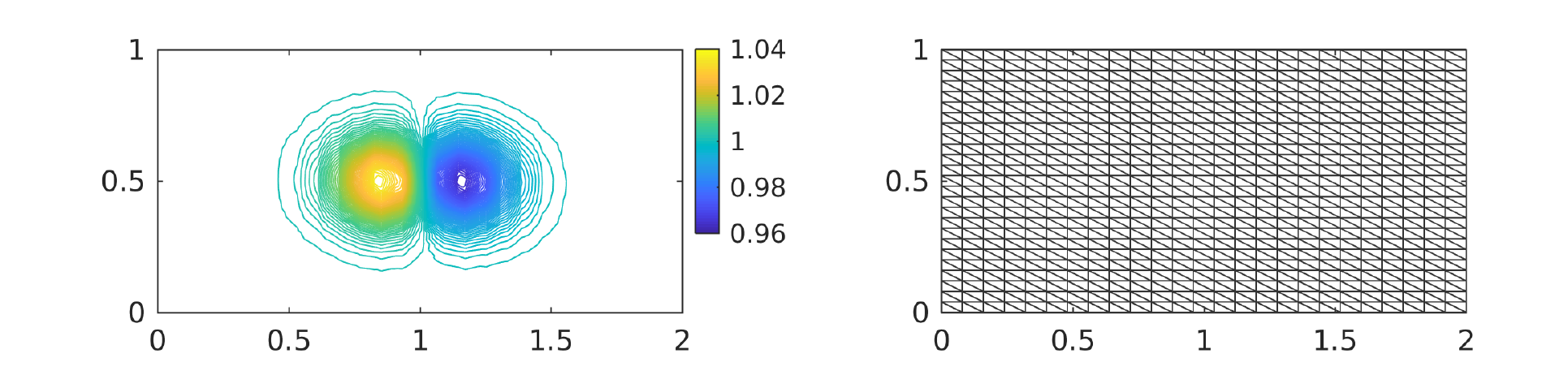}}\label{fig:1a}\\	
		\subfigure[Uniform mesh $2\times 50\times 50$.]{\includegraphics[width=1\textwidth]{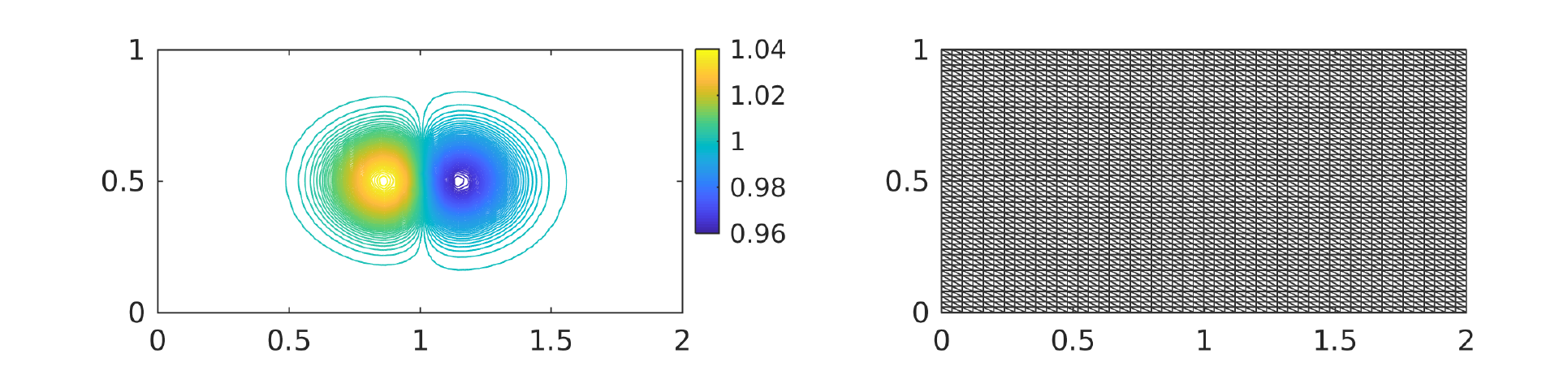}}\label{fig:1b}\\
		\subfigure[Adaptive mesh with $\mathcal{M}=1$.] {\includegraphics[width=1\textwidth]{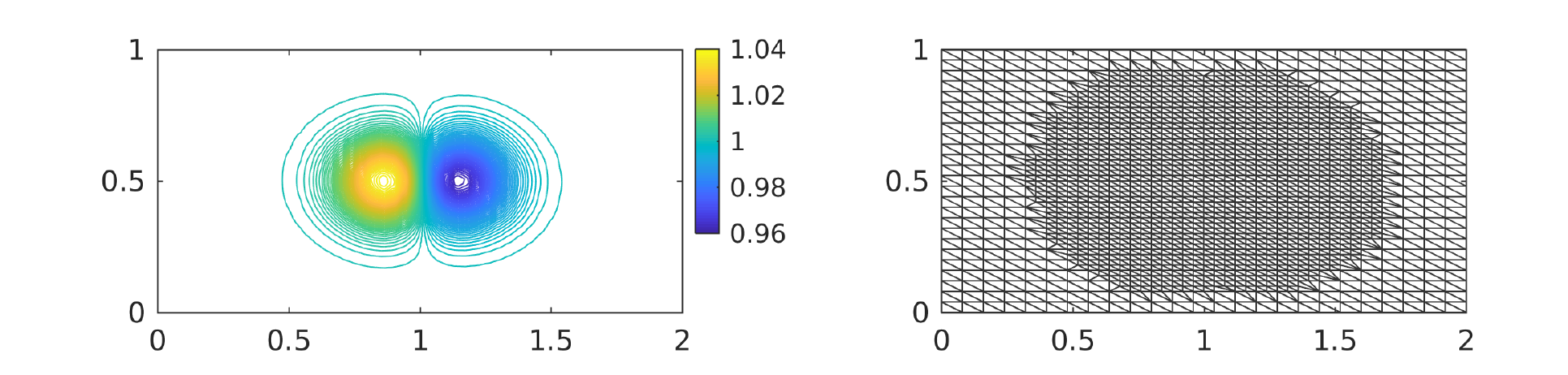}}\label{fig:1c}\\
		\subfigure[Adaptive mesh with $\mathcal{M}=2$.]{\includegraphics[width=1\textwidth]{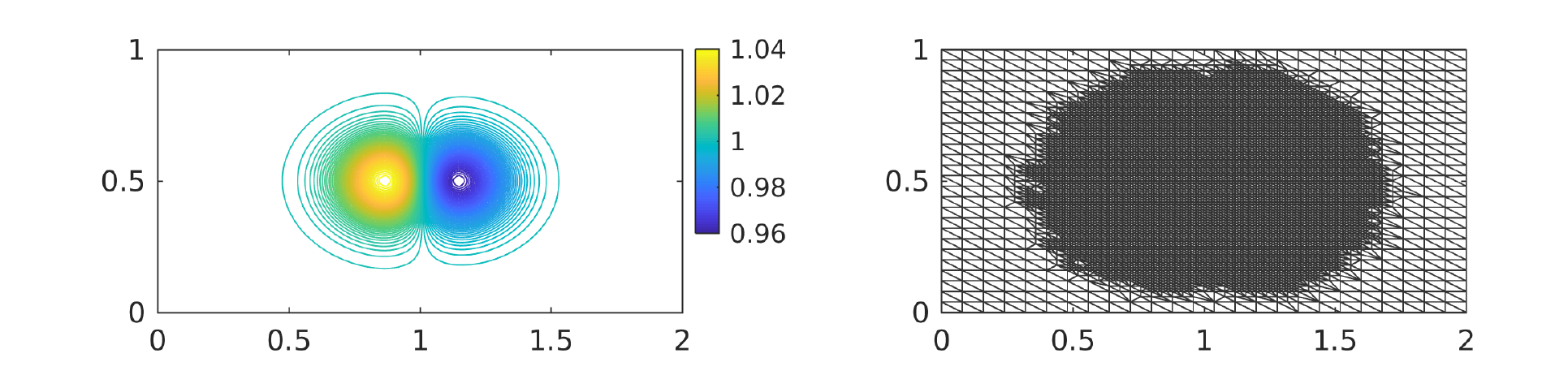}}\label{fig:1d}\\
		\vspace*{2mm}
		\caption{Example 1: Computational water surface $w(x,y,0.07)$ (left column) with the corresponding meshes (right column).}\label{fig:ex1a}
	\end{figure}

	From the result reported in \citep{MR2804645}, by $t=0.07$, the
	numerical solution reaches the steady state.  In
	\fref{fig:ex1a} (left column), we show the numerical solution of water
	surface at $t=0.07$. The solution is computed using the central-upwind scheme
	on uniform meshes on \fref{fig:ex1a} (a, b)  and using the
	adaptive central-upwind scheme on \fref{fig:ex1a} (c, d). The
	adaptive meshes in \fref{fig:ex1a} are obtained from the
	uniform mesh $2\times 25\times 25$, \fref{fig:ex1a} (a). The
	adaptive mesh  with one level of
	refinement $\mathcal{M}=1$ (as the highest level of
	refinement) is on \fref{fig:ex1a} (c),  and with two levels
	of refinement $\mathcal{M}=2$ (as the highest level of
	refinement) is on \fref{fig:ex1a} (d).
	
	Next, in \tref{tab:ex1a} we compute the $L^1$-errors of the central-upwind
	scheme on uniform meshes and of the adaptive central-upwind
	scheme.  To obtain the errors, the reference solution is calculated on the
	uniform mesh with 
	$2\times 800\times 800$ triangles. In  \tref{tab:ex1b}, we present the $\mathcal{R}_{CPU}$ ratio to compare the
	computational efficiency of the two methods. From \tref{tab:ex1a} and
	\tref{tab:ex1b}, we observe that the adaptive algorithm produces
	similar accuracy as the scheme on fixed uniform triangular
	meshes,  but at a  less
	computational cost. Also, as expected, the adaptive central-upwind
	scheme achieves second-order accuracy in space, similar to the
	central-upwind scheme without adaptivity.
	
	\begin{table}[h!]
		\vspace*{5mm}
		\centering
		\begin{tabular}{ |c| c c|c|c c|c|c c|}
			\hline
			\multicolumn{3}{|c|}{algorithm without adaptivity}&\multicolumn{6}{c|}{adaptive algorithm}\\
			\cline{4-9}
			\multicolumn{3}{|c|}{ }
			&\multicolumn{3}{c|}{one level $\mathcal{M}=1$}&\multicolumn{3}{c|}{two levels $\mathcal{M}=2$}\\
			\hline
			cells &$L^1$-error  &rate&cells	&$L^1$-error  &rate&cells	&$L^1$-error  &rate\\
			\hline
			$2\times100\times100$& 2.22e-05&&11,808&2.33e-05&&10,404&3.41e-05&	\\  
			$2\times200\times200$&4.74e-06	&2.22&46,892& 5.62e-06& 2.05&40,728& 5.91e-06&2.53\\
			$2\times400\times400$&9.84e-07&2.26&187,726&1.55e-06&1.85&160,846&1.54e-06&1.94\\
			\hline	
		\end{tabular}
		\caption{Example 1: $L^1$-errors of the water surface
			$w$ at $t=0.07$, and the convergence rates of the
			central-upwind scheme without adaptivity (uniform mesh $2\times N\times N, N=100,200,400$) and the
			adaptive scheme (the corresponding adaptive
                        mesh is reconstructed from the uniform mesh $2\times N/2^\mathcal{M}\times N/2^\mathcal{M}$).}\label{tab:ex1a}
		\vspace*{2mm}
	\end{table}

	\begin{table}[h!]
		\vspace*{5mm}
		\begin{tabular}{ |c|c| c c|c |c c|}
			\hline
			uniform mesh &\makecell{adaptive mesh\\ $\mathcal{M}=1$}&\multicolumn{2}{c|}{$\mathcal{R}_{CPU}$ with $\mathcal{M}=1$} &\makecell{adaptive mesh\\ $\mathcal{M}=2$}&\multicolumn{2}{c|}{$\mathcal{R}_{CPU}$ with $\mathcal{M}=2$}\\
			\cline{3-4} \cline{6-7}
			(cells)&(cells)&total &\makecell{without\\grid generation} &(cells)&total &\makecell{without\\grid generation}  \\
			\hline
			$2\times100\times100$&11,808&2.18&2.36&10,404&2.58& 2.70 \\
			$2\times200\times200$&46,892&1.85&2.05&40,728& 2.48 &2.62 \\
			$2\times400\times400$&187,726&1.77& 1.98 &160,846&2.25& 2.40 \\
			\hline
			\multicolumn{2}{|c}{$\mathcal{R}_{CPU}$ average:}&1.93&\multicolumn{1}{c}{2.13 }&\multicolumn{1}{c}{} &2.44& 2.57 \\
			\hline	
		\end{tabular}
		\caption{Example 1: $\mathcal{R}_{CPU}$ ratio at
			$t=0.07$, where for adaptive central-upwind scheme,
			we consider the total CPU times and CPU times
			without the grid generation.}\label{tab:ex1b}
	\end{table}
	
\begin{figure}[h!]
		\centering
		\subfigure[Non-uniform mesh with 1824 cells.]{\includegraphics[width=1\textwidth]{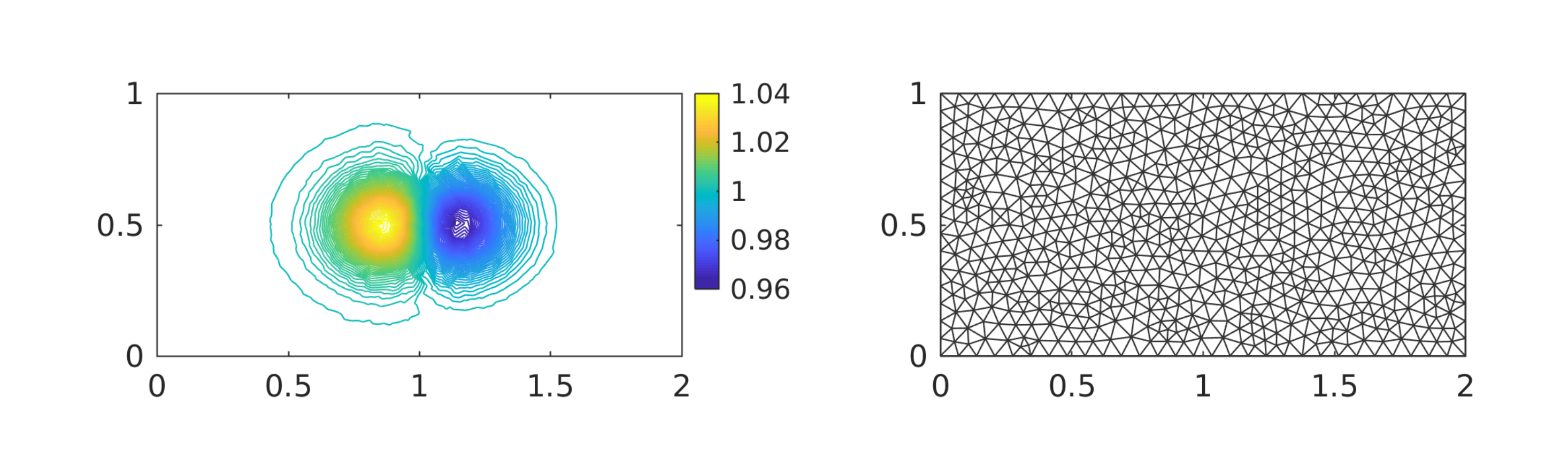}}\label{fig:ex1nona}\\	
		\subfigure[Adaptive mesh with $\mathcal{M}=1$.]{\includegraphics[width=1\textwidth]{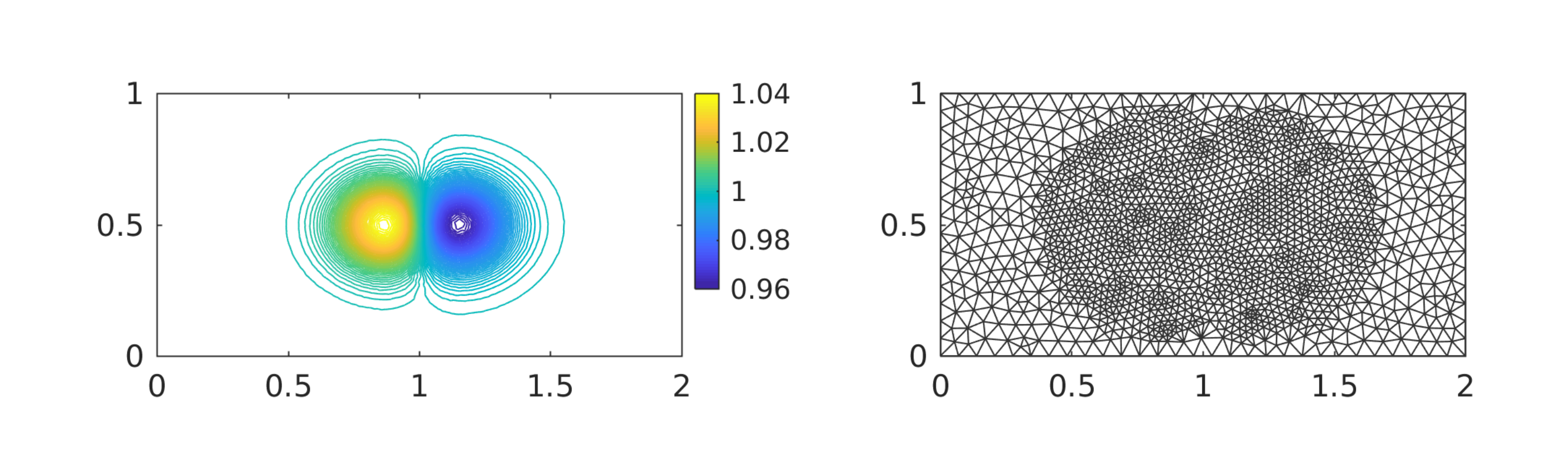}}\label{fig:ex1nonb}\\
        \subfigure[Adaptive mesh with $\mathcal{M}=2$.]{\includegraphics[width=1\textwidth]{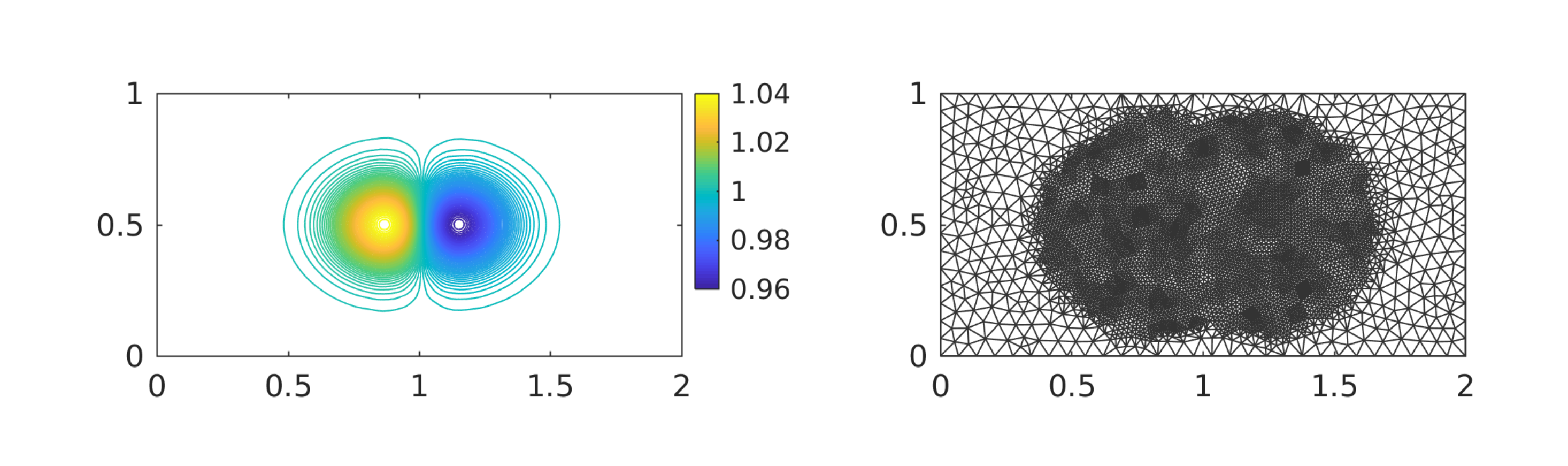}}\label{fig:ex1nonc}\\
        \vspace*{2mm}
			\caption{Example 1: Computational water surface $w(x,y,0.07)$ (left column) with the corresponding meshes (right column) obtained by applying the adaptive central-upwind scheme on a non-uniform mesh.}\label{fig:ex1non}
	\end{figure}

In addition, we will use this  example to show
          that the adaptive central-upwind scheme is also effective on
          the unstructured triangular meshes. On \fref{fig:ex1non}
          (left), we plot the numerical solution at $t=0.07$ computed
          using the central-upwind method on a non-uniform mesh,
          \fref{fig:ex1non} (a), the adaptive scheme with one level of
          refinement, \fref{fig:ex1non} (b) and the adaptive scheme with two 
	levels of refinement, \fref{fig:ex1non} (c). The adaptive meshes on \fref{fig:ex1non} (b, c) are reconstructed from the non-uniform mesh shown on the right \fref{fig:ex1non} (a).
	
	We then recompute the accuracy of the solution, see
        \tref{tab:ex1nona}, and the CPU time ratio, see
        \tref{tab:ex1nonb}, obtained by using the non-uniform
        meshes. From \tref{tab:ex1nona} and \tref{tab:ex1nonb}, we
        observe that the advantage of the adaptive scheme is still maintained on non-uniform meshes as it reduces the
	computational cost of the method. Next, we compare the results
        presented in \tref{tab:ex1nona} and in \tref{tab:ex1a}.  As
        expected, we
        see that the errors obtained by using scheme on non-uniform
        meshes is slightly larger than the ones using the corresponding
        uniform meshes (with approximately the same number of cells).

	\begin{table}[h!]
		\vspace*{5mm}
		\centering
		\begin{tabular}{ |c| c c|c|c c|c|c c|}
			\hline
			\multicolumn{3}{|c|}{algorithm without adaptivity}&\multicolumn{6}{c|}{adaptive algorithm}\\
			\cline{4-9}
			\multicolumn{3}{|c|}{ }
			&\multicolumn{3}{c|}{one level $\mathcal{M}=1$}&\multicolumn{3}{c|}{two levels $\mathcal{M}=2$}\\
			\hline
			cells &$L^1$-error  &rate&cells	&$L^1$-error  &rate&cells	&$L^1$-error  &rate\\
			\hline
			22,438& 2.64e-05&&13,036&2.53e-05&&11,630&3.59e-05&	\\  
			90,434&6.53e-06	&2.02&51,656&5.98e-06&2.08&44,400&7.06e-06&2.35\\
			314,708&1.15e-06&2.51&204,044&1.68e-06&1.83&176,278&1.60e-06&2.14\\
			\hline	
		\end{tabular}
		\caption{Example 1: $L^1$-errors of the water surface
			$w$ at $t=0.07$, and the convergence rates of the
			central-upwind scheme without adaptivity in non-uniform mesh and the
			adaptive scheme (the corresponding adaptive
			meshes has the same size of the smallest cells with the compared non-uniform meshes).}\label{tab:ex1nona}
		\vspace*{2mm}
	\end{table}

	\begin{table}[h!]
		\vspace*{5mm}
		\begin{tabular}{ |c|c| c c|c |c c|}
			\hline
			\makecell{Non-uniform and\\non-adaptive mesh} &\makecell{adaptive mesh\\ $\mathcal{M}=1$}&\multicolumn{2}{c|}{$\mathcal{R}_{CPU}$ with $\mathcal{M}=1$} &\makecell{adaptive mesh\\ $\mathcal{M}=2$}&\multicolumn{2}{c|}{$\mathcal{R}_{CPU}$ with $\mathcal{M}=2$}\\
			\cline{3-4} \cline{6-7}
			(cells)&(cells)&total &\makecell{without\\grid generation} &(cells)&total &\makecell{without\\grid generation}  \\
			\hline
			22,438&13,036&1.86&1.94&11,630&3.77& 3.90 \\
			90,434&51,656&2.20&2.20&44,400&3.34 &3.43\\
			314,708&204,044&2.45&2.57 &176,278&3.36& 3.46 \\
			\hline
			\multicolumn{2}{|c}{$\mathcal{R}_{CPU}$ average:}&2.17&\multicolumn{1}{c}{2.24 }&\multicolumn{1}{c}{} &3.49&3.60 \\
			\hline	
		\end{tabular}
		\caption{Example 1: $\mathcal{R}_{CPU}$ ratio at
			$t=0.07$, where for adaptive central-upwind scheme,
			we consider the total CPU times and CPU times
			without the grid generation in non-uniform meshes.}\label{tab:ex1nonb}
	\end{table}
	\vspace*{5mm}

	\vspace*{5mm}

	\subsection{Example 2---Well-Balanced  Tests and Test with Wet/Dry
		Interfaces}\label{sec4_2}
	The first numerical example here was proposed in \citep{LeV} to test
	capability  of the numerical scheme to accurately resolve small perturbations of a steady state solution. We take a computational domain $[0, 2]\times[0, 1]$ and the bottom topography, 
	
	\begin{equation}\label{eq:ex2B}
	B(x,y)=0.8\exp(-5(x-0.9)^2-50(y-0.5)^2).
	\end{equation}
	
	The initial conditions describe a flat surface of water with a small perturbation in $0.05<x<0.15$:
	
	\begin{equation}
	w(x,y,0)=\begin{cases}
	1 + \epsilon, \quad& 0.05<x<0.15,\\ 1, & \mbox{otherwise,}
	\end{cases}\quad u(x,y,0)\equiv v(x,y,0)\equiv 0,\label{eq:ex2I}
	\end{equation}
	
	where  $ \epsilon$ is the perturbation height. We have used the
	zero-order extrapolation at the right and the left boundaries
	of the domain and the periodic boundary conditions for the top
	and the bottom ones. 
	
\begin{figure}[h!]
		\centering
		\centering
		\subfigure[Central-upwind scheme on a uniform mesh $2\times 25\times 25$.]{\includegraphics[width=1\textwidth]{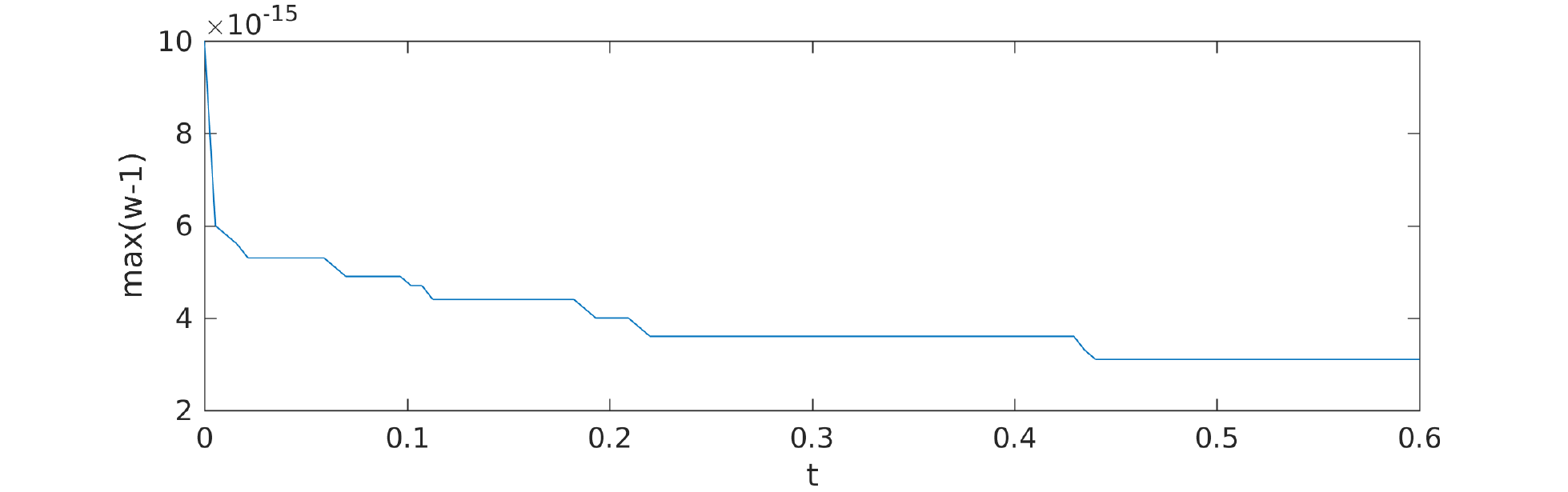}}
		\subfigure[Adaptive central-upwind scheme on a mesh
		obtained from the uniform mesh $2\times 25\times 25$ with $\mathcal{M}=1$.] {\includegraphics[width=1\textwidth]{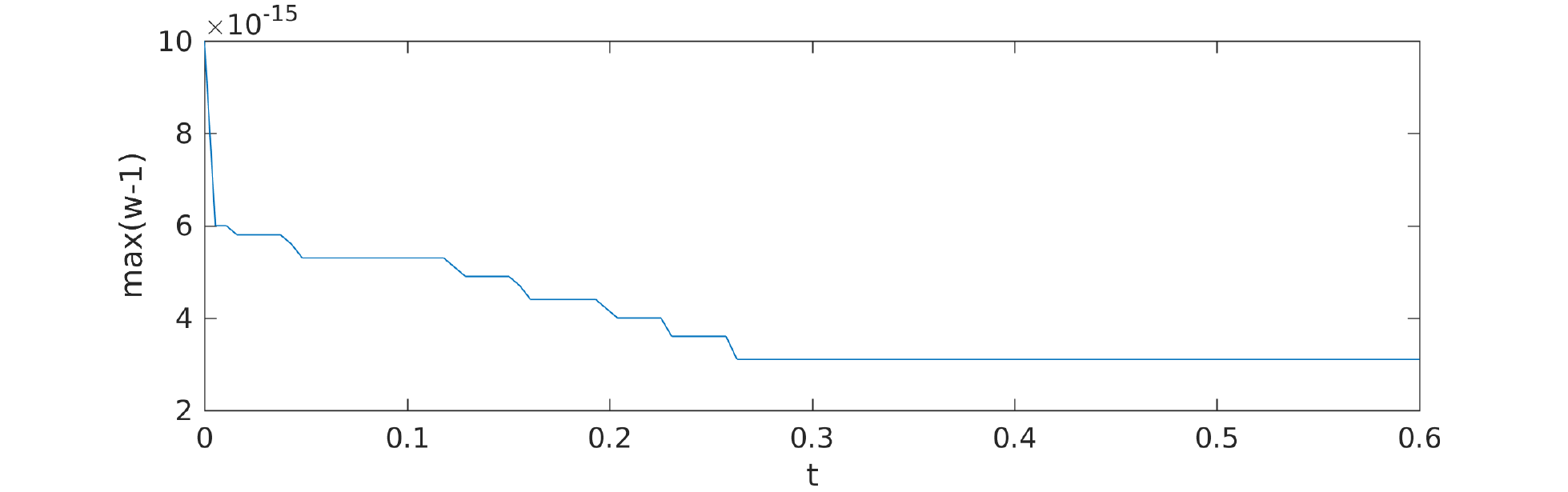}}
		\subfigure[Adaptive central-upwind scheme on a mesh
		obtained from the uniform mesh $2\times 25\times 25$ with $\mathcal{M}=2$.] {\includegraphics[width=1\textwidth]{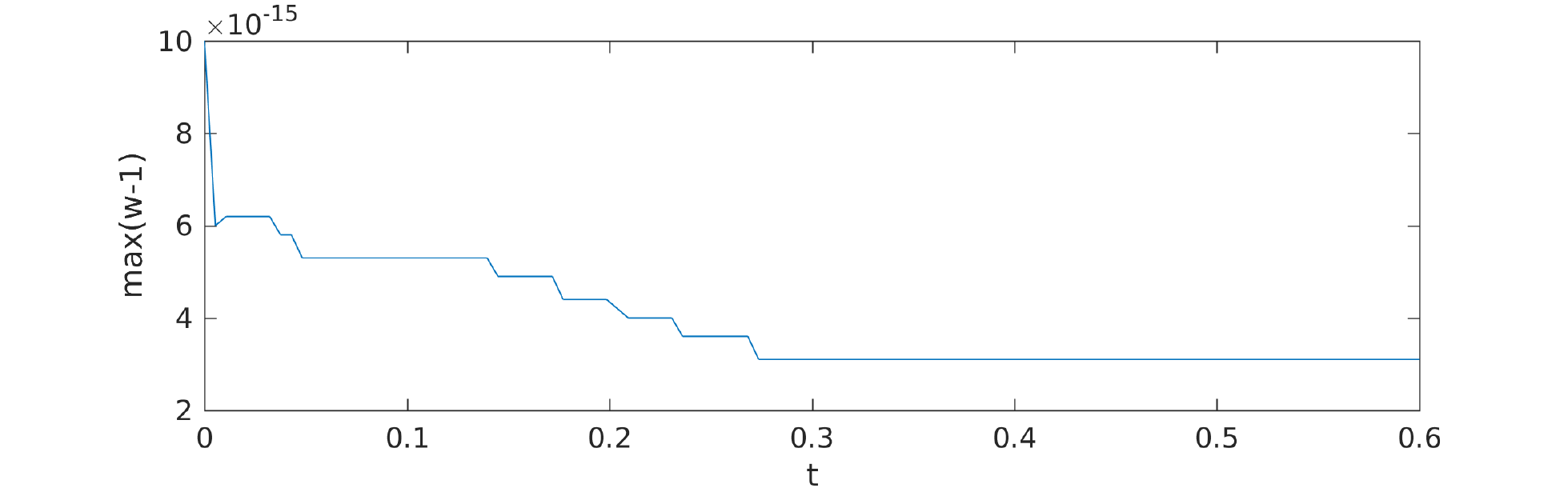}}
		\vspace*{5mm}
		\caption{Example 2: $\displaystyle
			\max_{x,y}(w(x,y,t)-1)$ is plotted as a function of $t$ on the uniform grid and the adaptive grids.}\label{fig:ex2a}
		\vspace*{2mm}
	\end{figure}

	To verify well-balanced property of the adaptive scheme, we
	first consider a very small perturbation $ \epsilon=10^{-14}$. The
	adaptive meshes  with levels $\mathcal{M}=1, 2$ are obtained from
	a coarse uniform mesh $2\times 25\times 25$. The threshold for
	mesh refinement in this example is $\displaystyle\omega =
	0.1\max_j(e_j)$.  We plot $\displaystyle
	\max_{x,y}(w-1)$ as a function of time for the central-upwind
	scheme without adaptivity on uniform mesh in \fref{fig:ex2a}
	(a), for the adaptive central-upwind scheme with $\mathcal{M}=1$
	in \fref{fig:ex2a} (b),  and with $\mathcal{M}=2$ in \fref{fig:ex2a}
	(c). The results of the test show that the adaptive scheme is
	stable and preserves numerically the balance between the
	fluxes and the source term, similar to the scheme without adaptivity.
	
\begin{figure}
		\centering
		\includegraphics[width=1\textwidth]{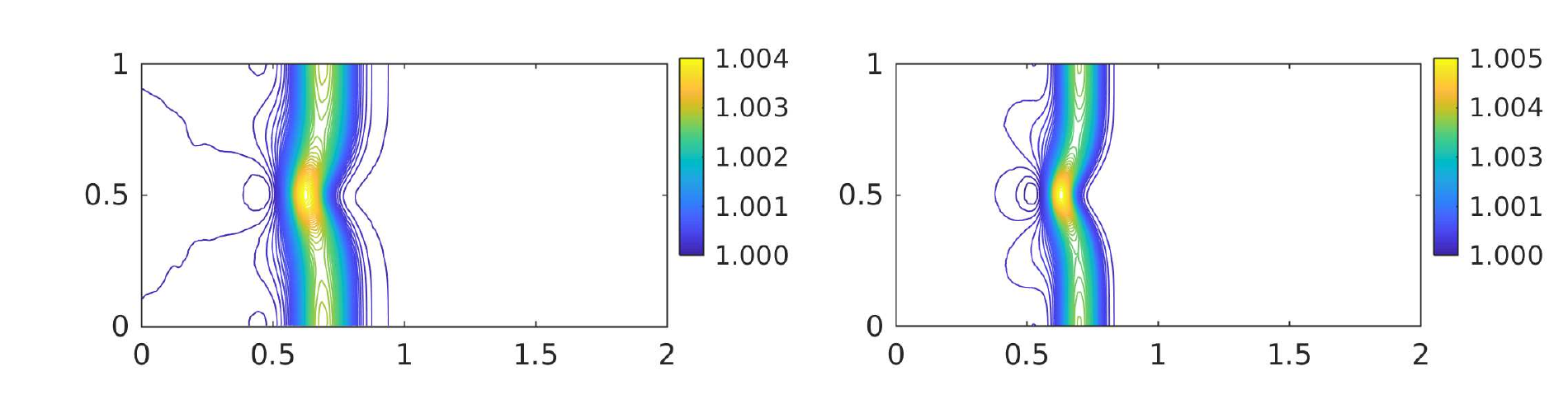}\\
		\includegraphics[width=1\textwidth]{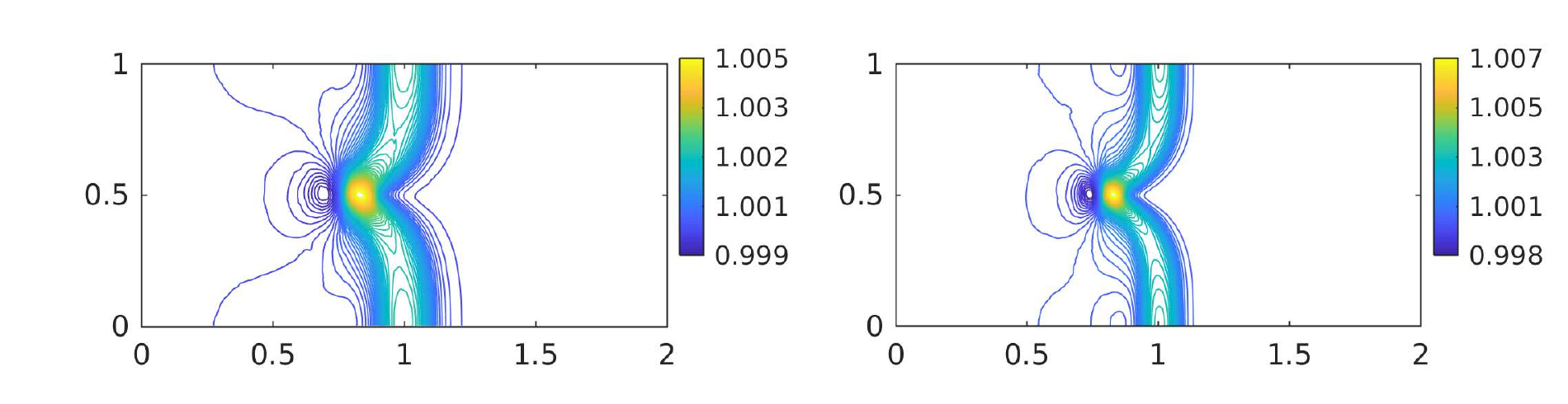}\\
		\includegraphics[width=1\textwidth]{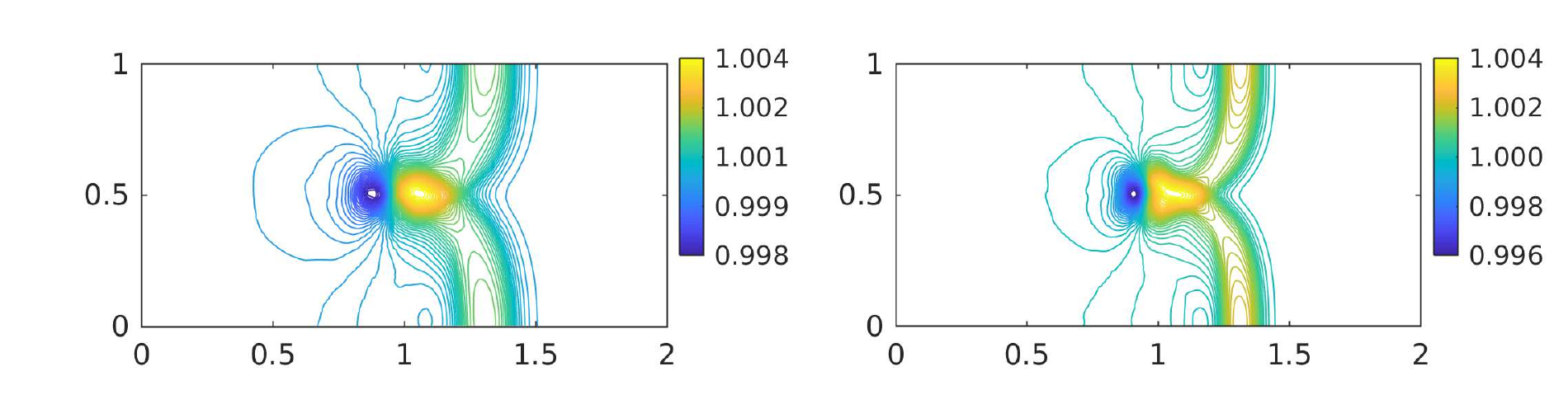}\\
		\includegraphics[width=1\textwidth]{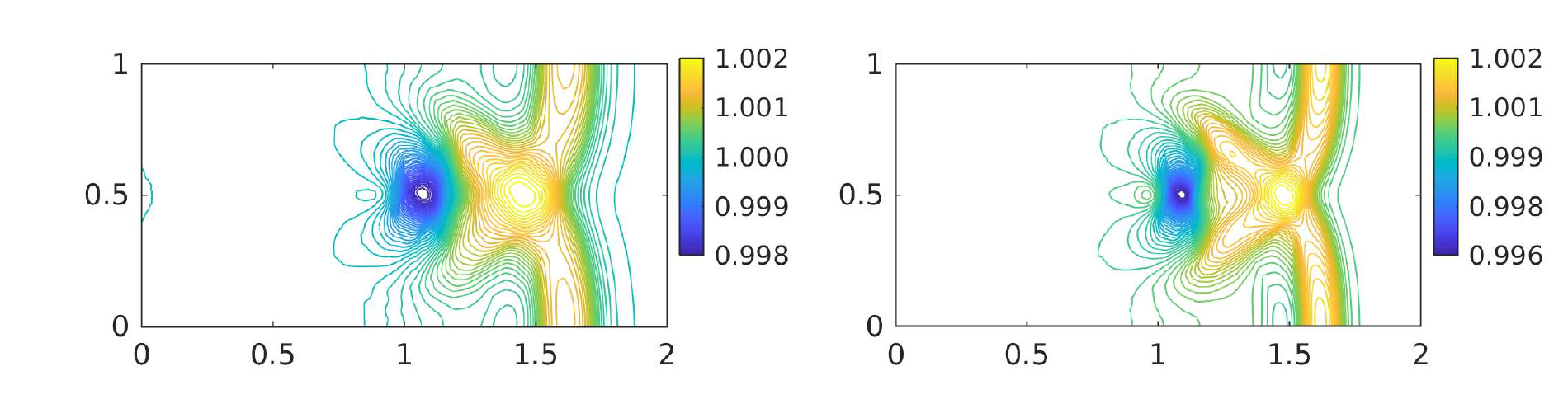}\\
		\includegraphics[width=1\textwidth]{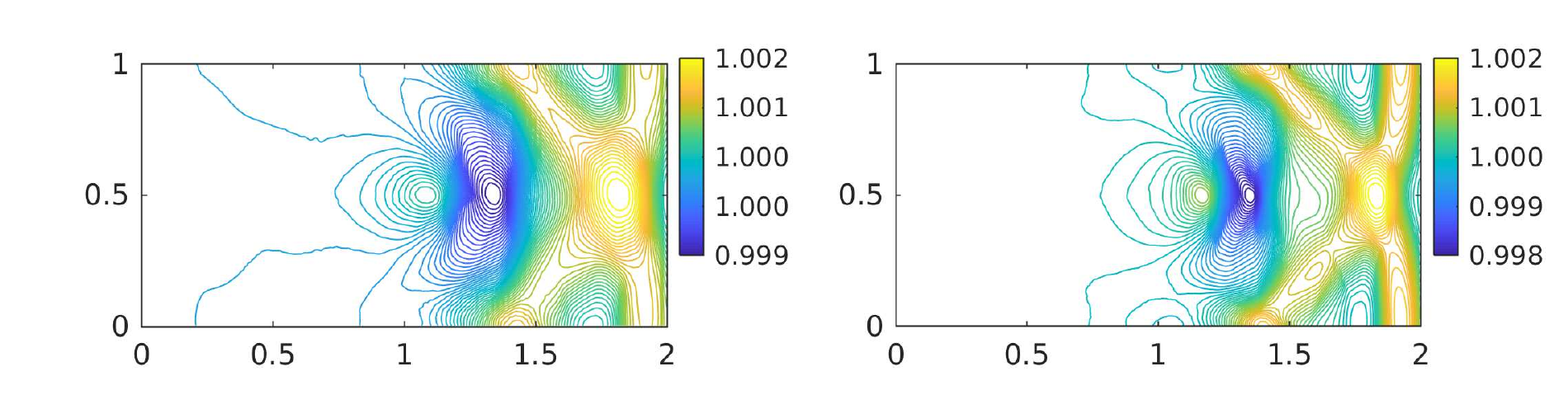}
		\caption{Example 2: $w$ component of the solution of
			the IVP
			\eref{eq:swe1}-\eref{eq:swe3},\eref{eq:ex2B}-\eref{eq:ex2I}
			with $\epsilon=10^{-2}$ at $t=0.6, 0.9, 1.2, 1.5$ and
			$1.8$ (from top to bottom) obtained by the
			central-upwind scheme without adaptivity on uniform meshes $2\times 100\times 100$ (left column) and $2\times 200\times 200$ (right column).} \label{fig:ex2u}
	\end{figure}

\begin{figure}[htbp]
		\centering
		\includegraphics[width=1\textwidth]{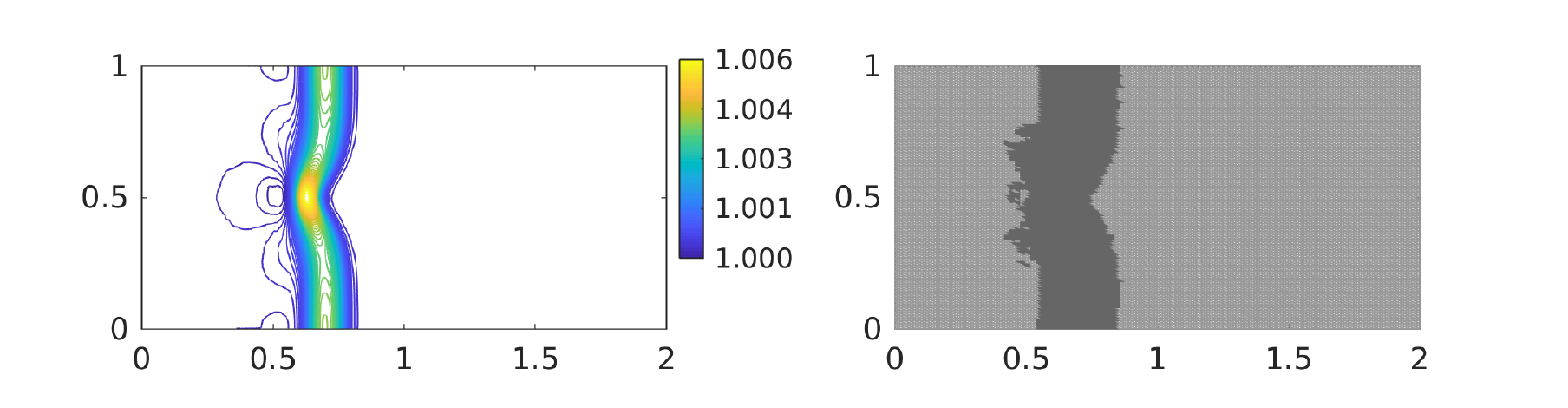}\\
		\includegraphics[width=1\textwidth]{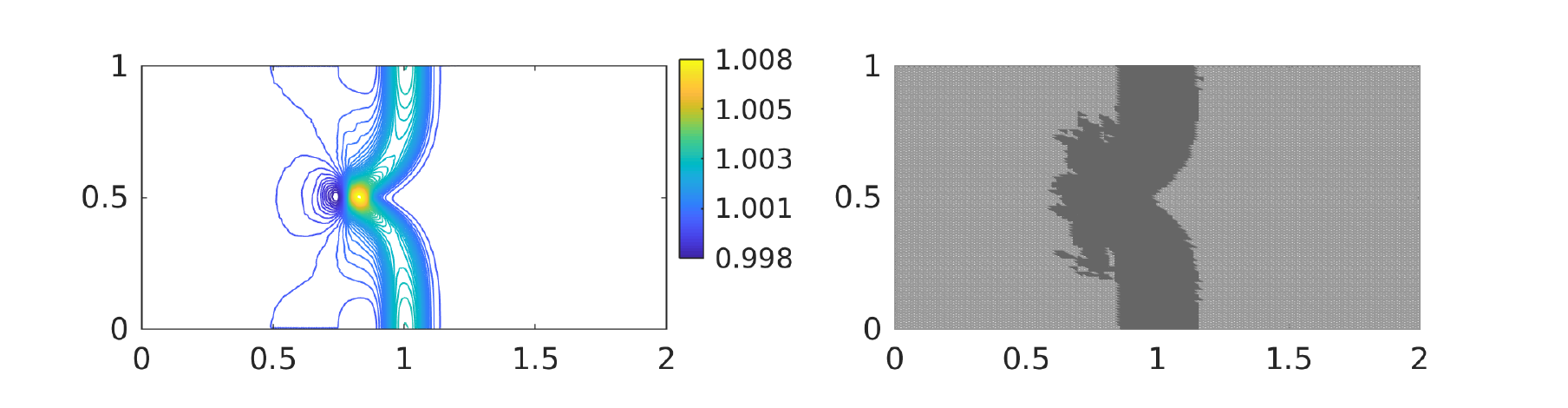}\\
		\includegraphics[width=1\textwidth]{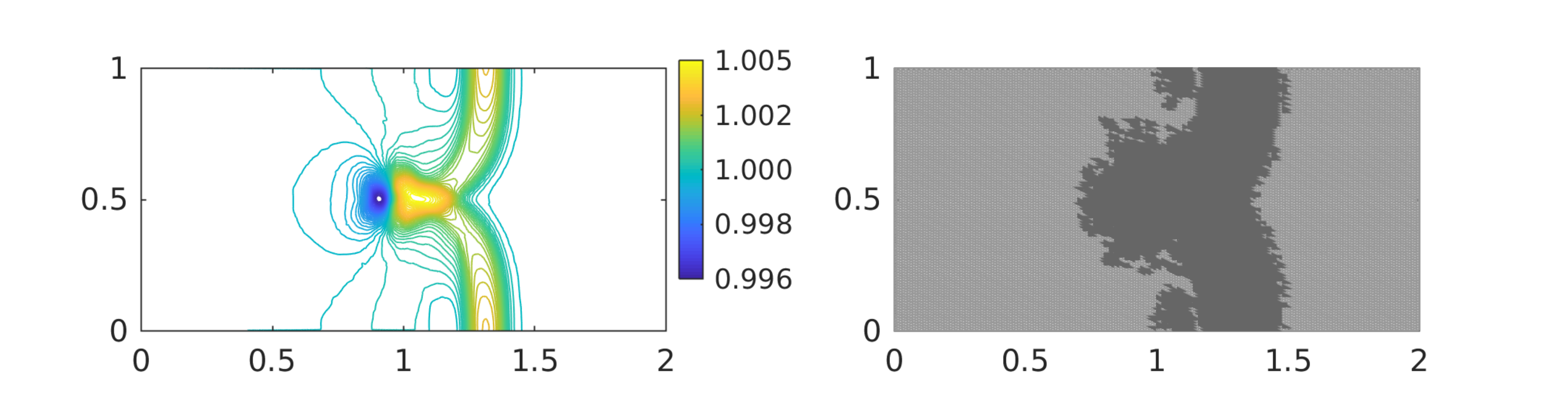}\\
		\includegraphics[width=1\textwidth]{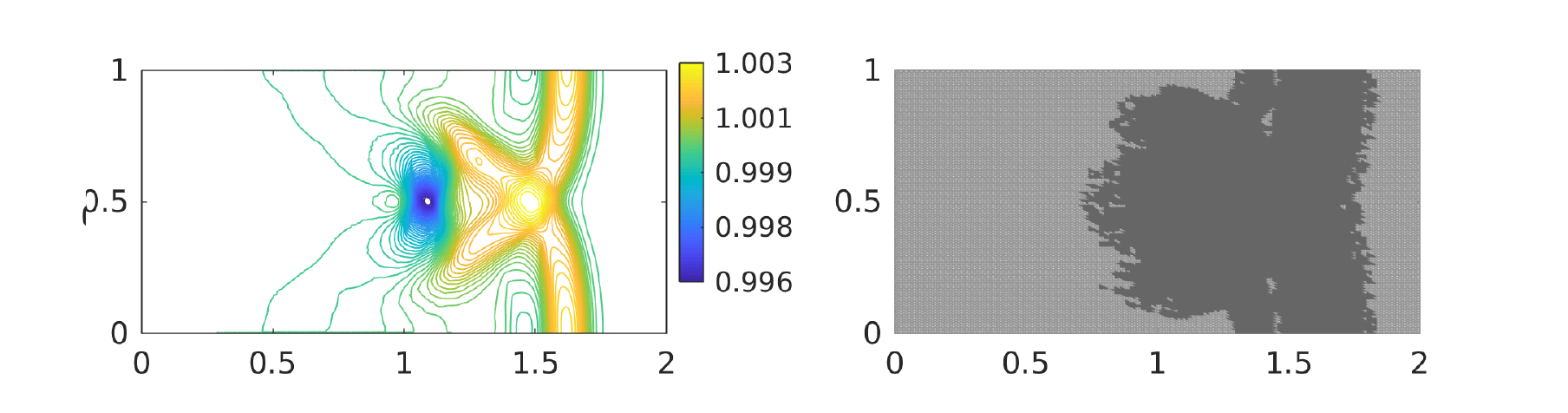}\\
		\includegraphics[width=1\textwidth]{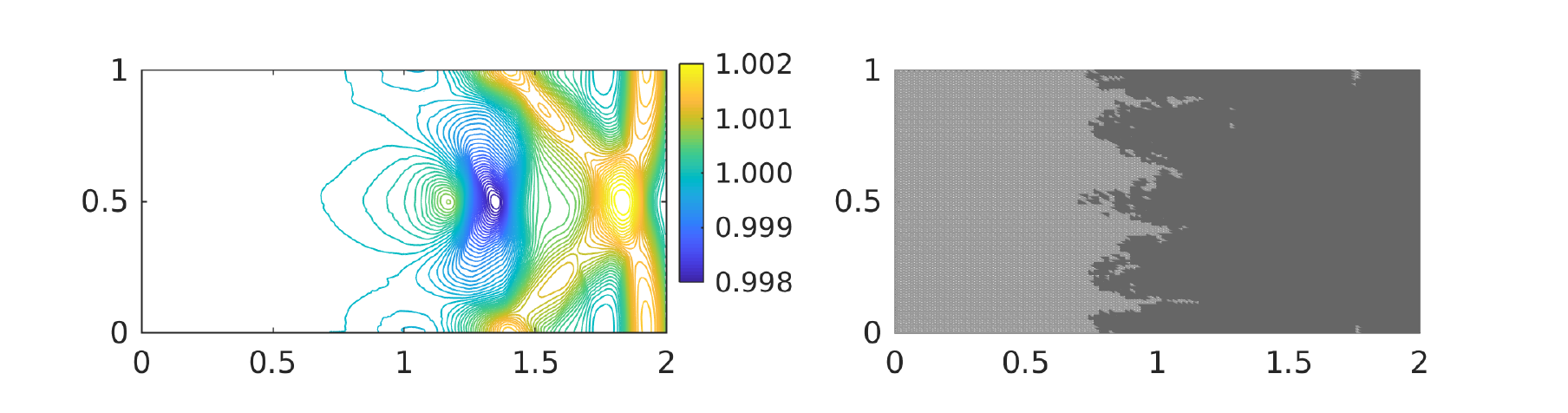}
		\caption{Example 2: $w$ component of the solution of
			the IVP
			\eref{eq:swe1}-\eref{eq:swe3},\eref{eq:ex2B}-\eref{eq:ex2I}
			with $\epsilon=10^{-2}$ (left column) at $t=0.6, 0.9,
			1.2, 1.5$ and $1.8$ (from top to bottom) obtained by the adaptive central-upwind scheme. The
			corresponding adaptive meshes with one level of
			refinement $\mathcal{M}=1$ (right column).}  \label{fig:ex2r1l}
	\end{figure}

\begin{figure}[htbp]
		\centering
		\includegraphics[width=1\textwidth]{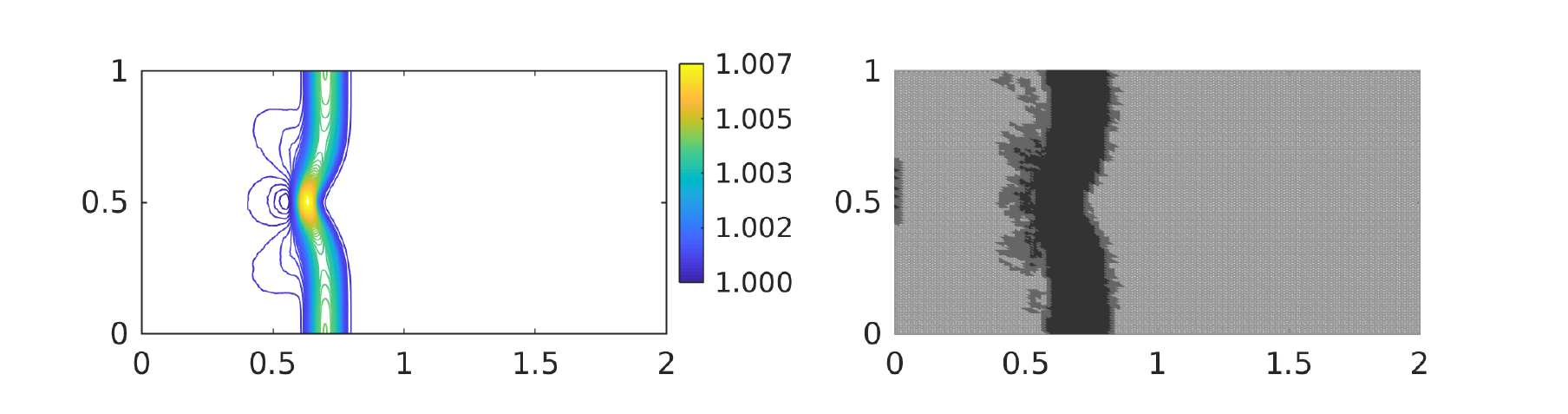}\\
		\includegraphics[width=1\textwidth]{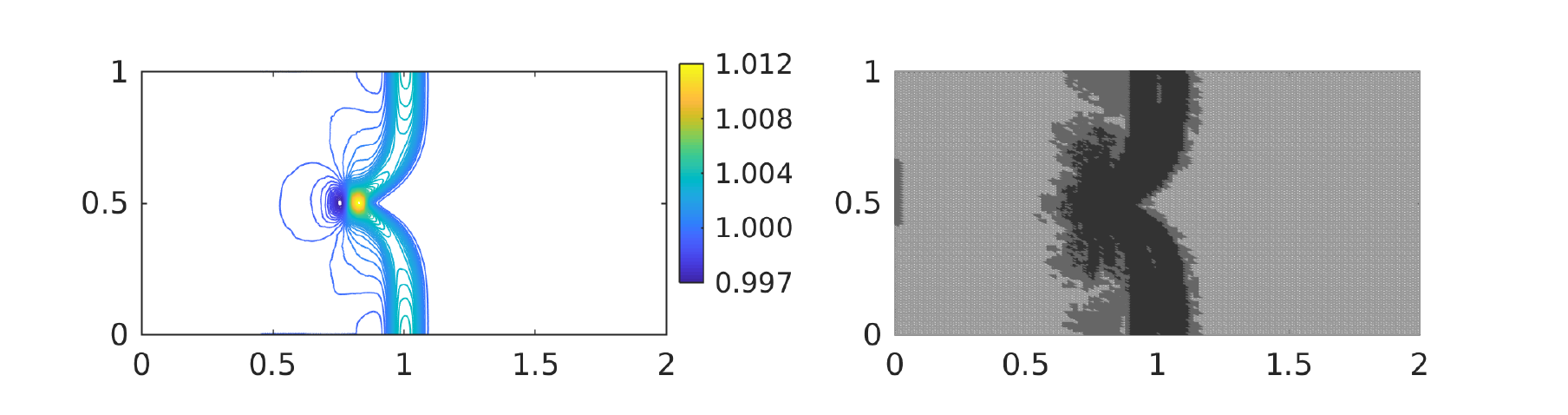}\\
		\includegraphics[width=1\textwidth]{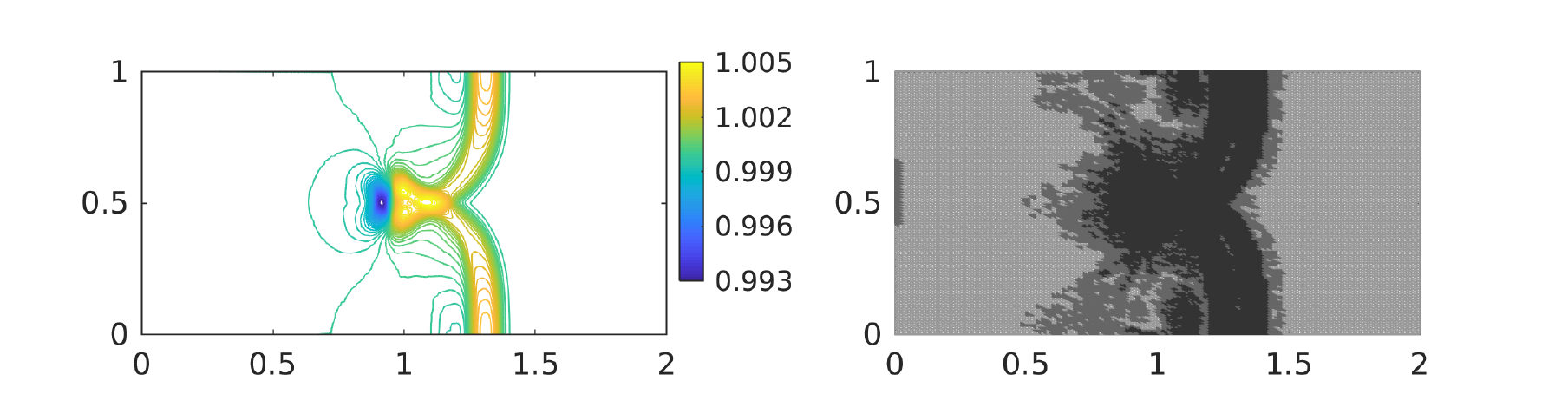}\\
		\includegraphics[width=1\textwidth]{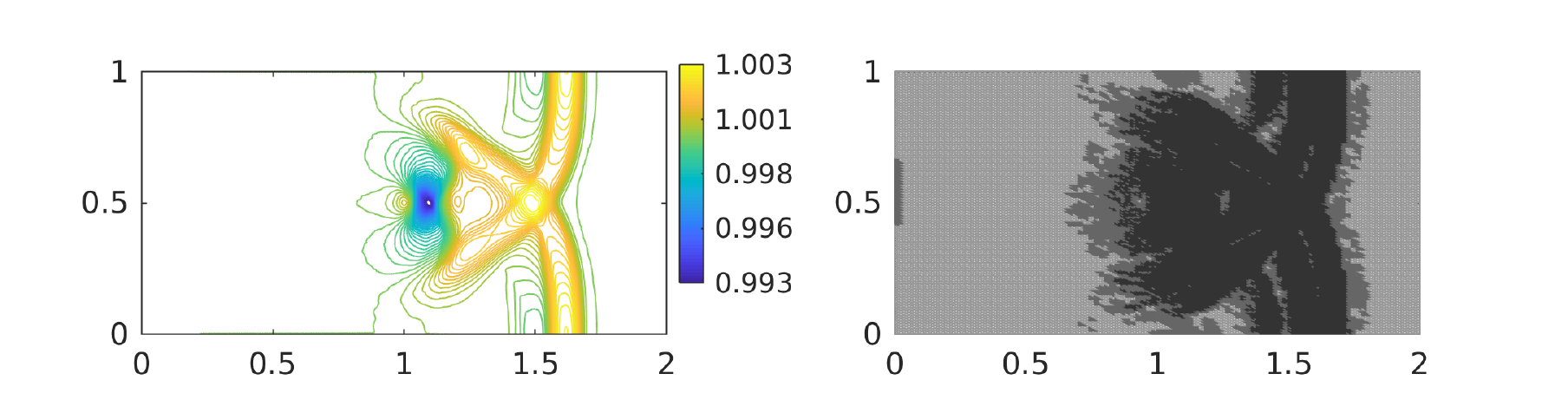}\\
		\includegraphics[width=1\textwidth]{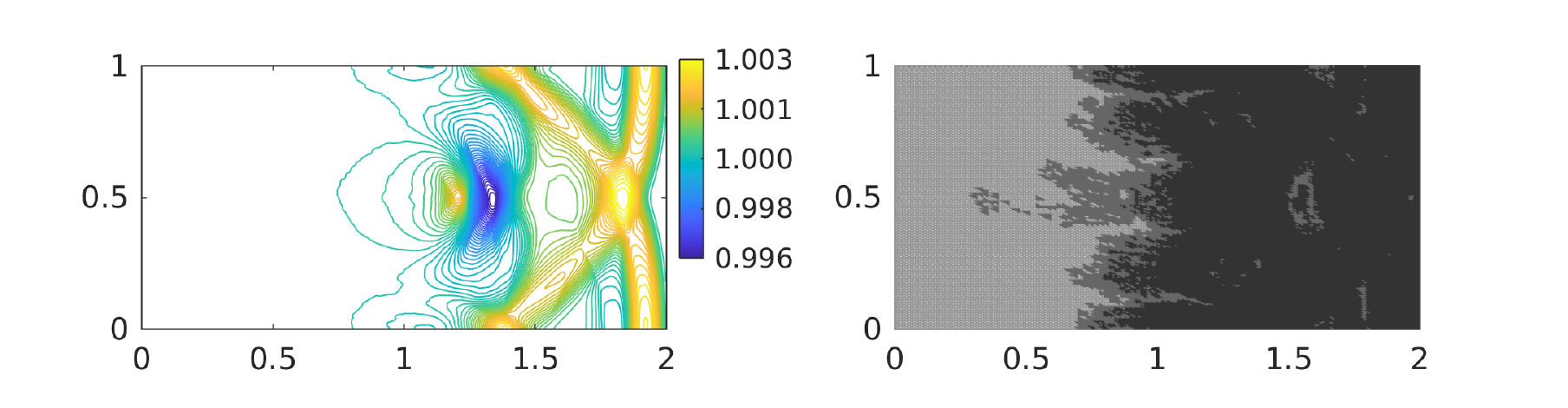}
		\caption{Example 2: $w$ component of the solution of
			the IVP
			\eref{eq:swe1}-\eref{eq:swe3},\eref{eq:ex2B}-\eref{eq:ex2I}
			with $\epsilon=10^{-2}$ (left column) at $t=0.6, 0.9,
			1.2, 1.5$ and $1.8$ (from top to bottom) obtained by
			the adaptive central-upwind scheme. The
			corresponding adaptive meshes with two levels of
			refinement $\mathcal{M}=2$ (right column).} \label{fig:ex2r2l}
	\end{figure}

	In the next numerical test, we take a larger perturbation
	value $\epsilon= 10^{-2}$. The purpose of this test is to
	demonstrate the capability of the adaptive algorithm to
	resolve the small scale features of the solution. In
	\fref{fig:ex2u}, we plot the water surface $w$ obtained by the
	method without adaptivity on two uniform meshes $2\times 100\times 100$ (left) and $2\times 200\times 200$ (right)  at $t=0.6, 0.9, 1.2, 1.5$ and
	$1.8$. The computed solutions of the water surface exhibit a
	right-going disturbance propagating past the hump. We then
	apply the adaptive algorithm and plot the numerical solution
	of the water surface $w$ on selected meshes at different times in
	\fref{fig:ex2r1l} (left) with $\mathcal{M}=1$ and in
	\fref{fig:ex2r2l}  (left) with $\mathcal{M}=2$ (the starting grid
	was a uniform mesh with $2\times 100\times 100$ and the threshold is $\displaystyle\omega =
	0.1\max_j(e_j)$).  We observe
	from  \fref{fig:ex2u}, \fref{fig:ex2r1l} and
	\fref{fig:ex2r2l} that the adaptive central-upwind scheme
	delivers high resolution of the complex features of a small
	perturbation of the "lake-at-rest" steady state. We note also,
	that by increasing the level of refinement from $\mathcal{M}=1$ to
	$\mathcal{M}=2$, the number of cells in the mesh increases from
	$30,912$ cells with $\mathcal{M}=1$ to $66,097$ cells with
	$\mathcal{M}=2$, but the accuracy is clearly improved with
	higher resolution as seen in \fref{fig:ex2r1l} and
	\fref{fig:ex2r2l}.  We also present the corresponding adaptive
	meshes in \fref{fig:ex2r1l} (right) and in \fref{fig:ex2r2l}
	(right). Clearly, the meshes are adapted to the behavior of
	the flow during time evolution from $t=0.6$ to $t=1.8$. This
	confirms the ability of the WLR error estimator to detect
	location of the steep local gradients in the solution.

	To  demonstrate further the advantage of the adaptive scheme,
	we have compared the CPU times of the central-upwind scheme
	without adaptivity and with adaptivity for the solution at
	$t=0.9$. The computed  $\mathcal{R}_{CPU}$ ratios are
	presented  in \tref{tab:ex2b}. The computed ratios
	$\mathcal{R}_{CPU}$ show that the adaptive central-upwind
	scheme can reduce the numerical cost by about three times with
	$\mathcal{M}=1$ and by about five times with $\mathcal{M}=2$ in
	comparison with the cost of the algorithm without adaptivity.

	\begin{table}[h!]
		\vspace*{2mm}
		\centering
		\begin{tabular}{ |c|c| c|c|c |}
			\hline
			\makecell{uniform mesh\\(cells)} &\makecell{adaptive mesh\\ $\mathcal{M}=1$\\(cells)}&\makecell{$\mathcal{R}_{CPU}$ \\with $\mathcal{M}=1$} &\makecell{adaptive mesh\\$\mathcal{M}=2$\\(cells)}&\makecell{$\mathcal{R}_{CPU}$ \\with $\mathcal{M}=2$}\\
			\hline
			$2\times100\times100$&9,127& 3.50 &8,274&4.06  \\
			$2\times200\times200$&30,912 &2.06&21,131& 5.29  \\
			$2\times400\times400$&108,297& 3.67 &66,097&6.38 \\
			\hline
			\multicolumn{2}{|c}{$\mathcal{R}_{CPU}$  average: }&\multicolumn{1}{c}{3.07 } &\multicolumn{1}{c}{ } &\multicolumn{1}{c|}{5.24}\\
			\hline
		\end{tabular}
		\caption{Example 2: The $R_{CPU}$ ratio obtained in solving the IVP \eref{eq:swe1}-\eref{eq:swe3},\eref{eq:ex2B}-\eref{eq:ex2I} with $\epsilon=10^{-2}$ at $t=0.9$.}\label{tab:ex2b}
		\vspace*{2mm}
	\end{table}

\begin{figure}
		\centering
		\includegraphics[width=0.5\textwidth]{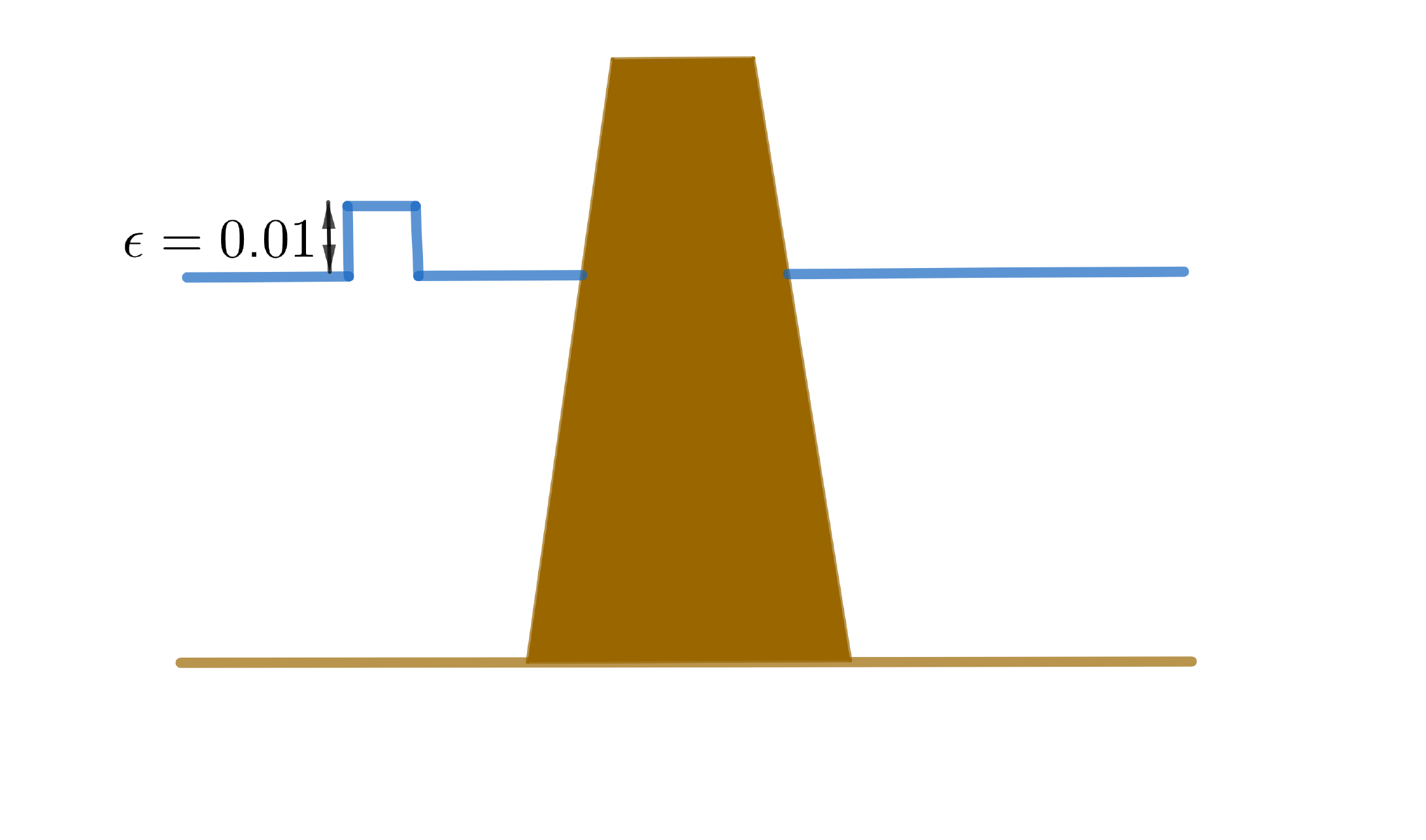}
		\caption{Example 2: 1-D slice of the bottom topography (brown) and water surface (blue) at $t=0$. The plot is not to
			scale.} \label{fig:ex2int}
	\end{figure}

	In the final test, we consider an example with a small
	perturbation that propagates around an island. Similar
	examples were considered in \citep{MR2804645,LAEK}.  We consider a
	hump partially submerged in water so that there is a
	disk-shaped island at the origin, see \fref{fig:ex2int}. Hence, the bottom topography is given by
	\begin{equation}\label{eq:ex2cB}
	B(x,y)=\begin{cases}
	1.1, \quad &r\leq 0.1,\\
	11\times (0.2-r)\quad  &0.1< r\leq 0.2,\\
	0, & \mbox{otherwise,}\end{cases}\quad r:=\sqrt{(x-0.5)^2+(y-0.5)^2}.\end{equation}
	
	We consider the following initial condition, 
	
	\begin{equation}\label{eq:ex2cI}
	w(x,y,0)=\begin{cases}
	1+\epsilon, \quad &0.1<x<0.2,\\ \max(1,B(x,y)), & \mbox{otherwise,}
	\end{cases}\quad u(x,y,0)\equiv v(x,y,0)\equiv 0,
	\end{equation}
	
	where $\epsilon=10^{-2}$ is the perturbation height. The
	homogeneous Neumann boundary conditions are used at all
	boundaries.
	
	We first obtain the water surface using central-upwind scheme
	without adaptivity. Different to the results in the submerged
	plateau case \eref{eq:ex2B}-\eref{eq:ex2I}, the right-going
	disturbance bends around the island and is being reflected by the island.

\begin{figure}
		\centering
		\includegraphics[width=1\textwidth]{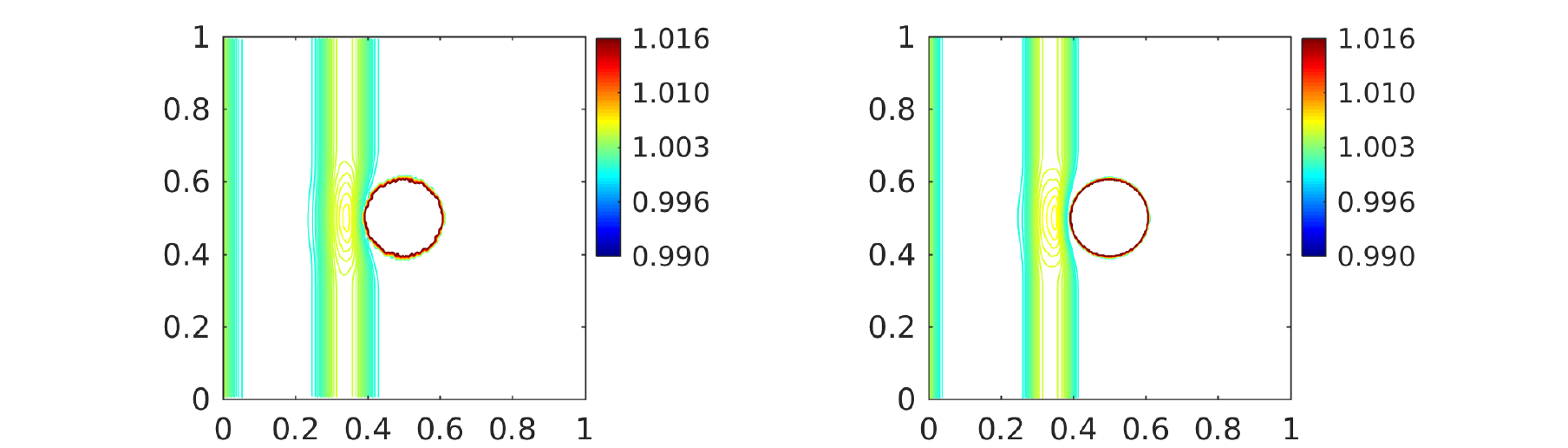}\\
		\includegraphics[width=1\textwidth]{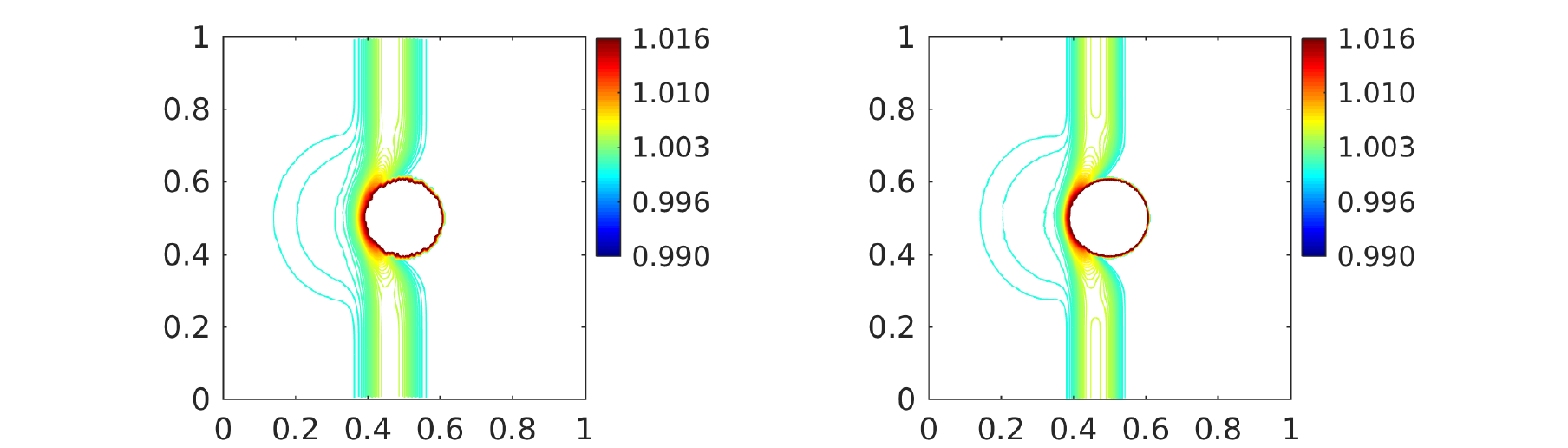}\\
		\includegraphics[width=1\textwidth]{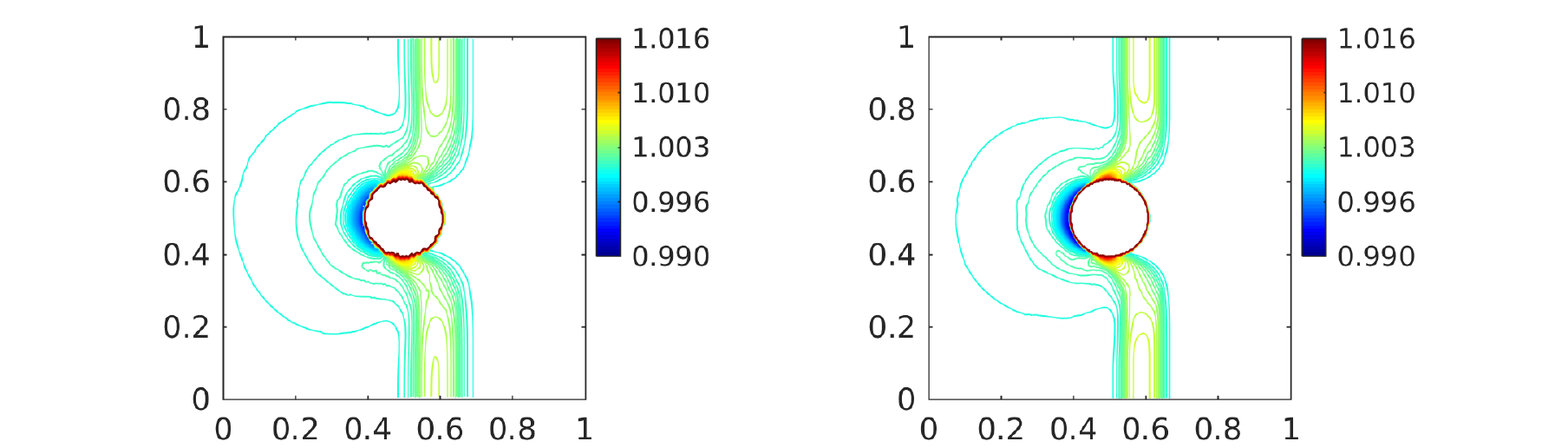}\\
		\includegraphics[width=1\textwidth]{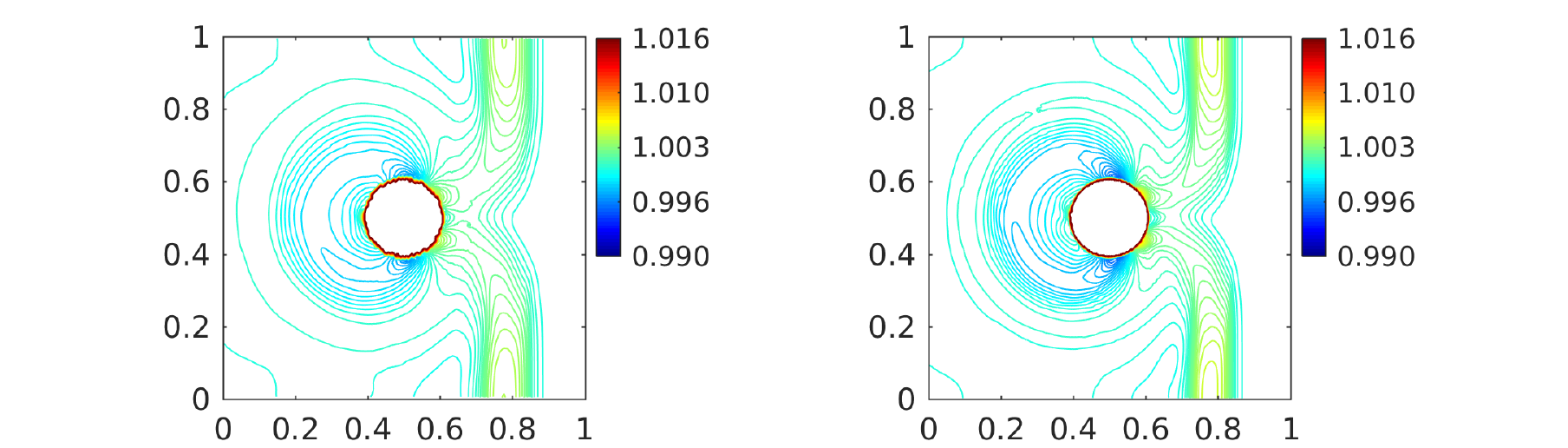}\\
		\vspace*{5mm}
		\caption{Example 2: $w$ component of the solution of
			the IVP \eref{eq:swe1}-\eref{eq:swe3},
			\eref{eq:ex2cB}-\eref{eq:ex2cI} at $t = 0.06, 0.1,
			0.14,$ and $0.2$ (from top to bottom) obtained by
			the central-upwind scheme without adaptivity on uniform meshes $2\times 100\times 100$ (left column) and $2\times 200\times 200$ (right column).} \label{fig:ex2cu}
	\end{figure}

\begin{figure}[htbp]
		\centering
		\includegraphics[width=1\textwidth]{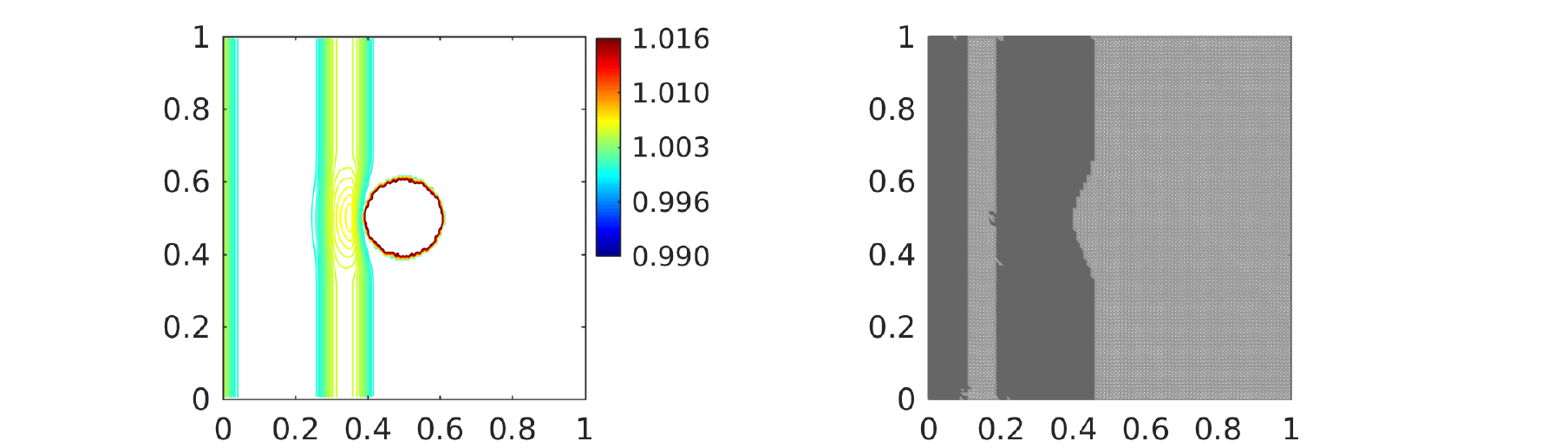}\\
		\includegraphics[width=1\textwidth]{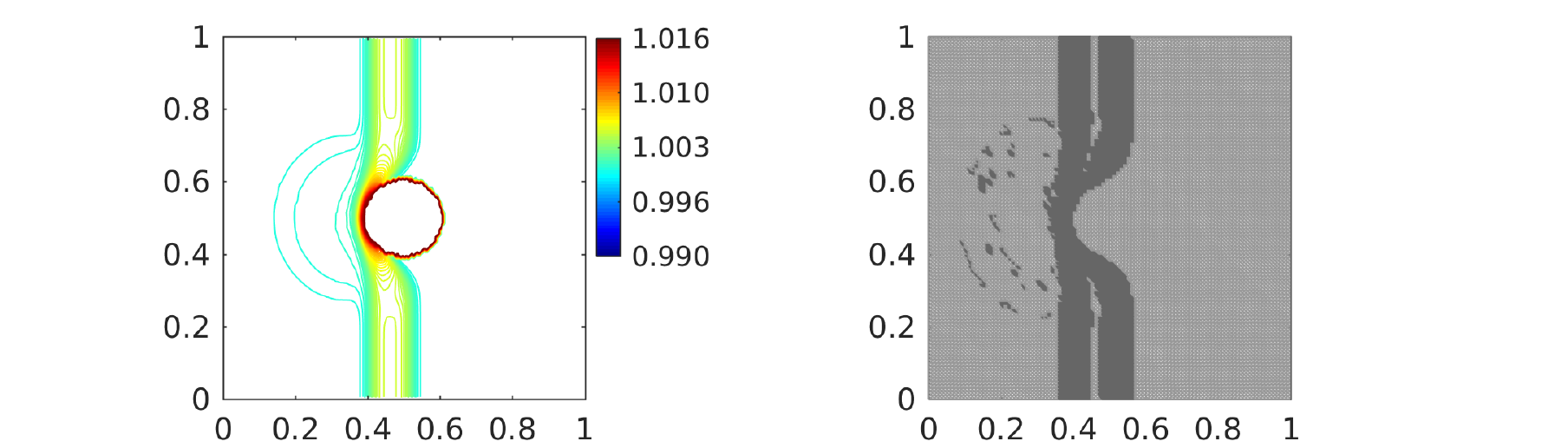}\\
		\includegraphics[width=1\textwidth]{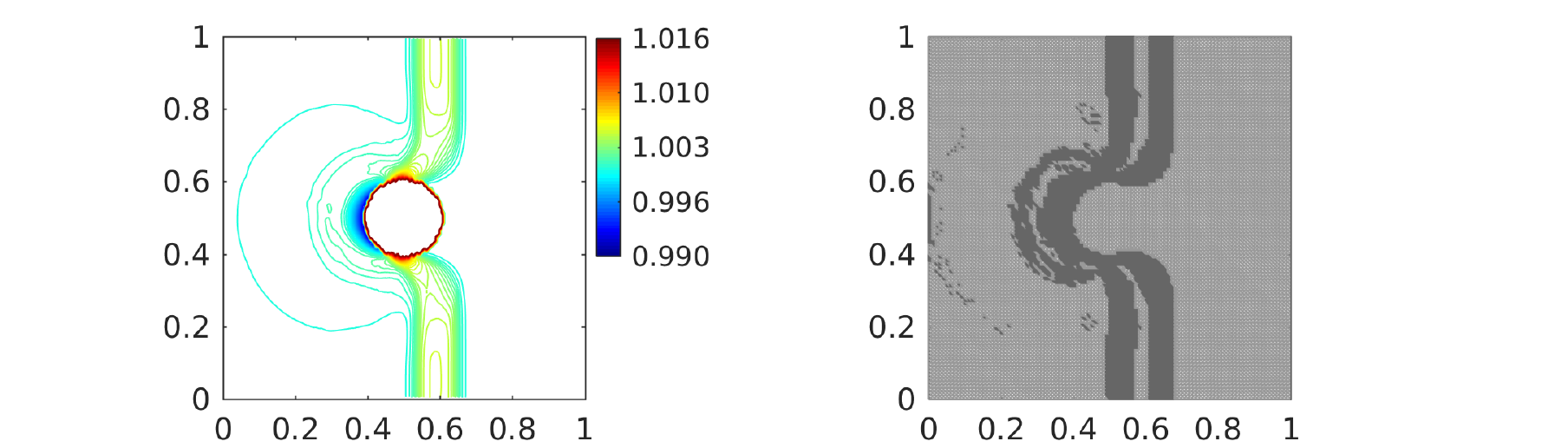}\\
		\includegraphics[width=1\textwidth]{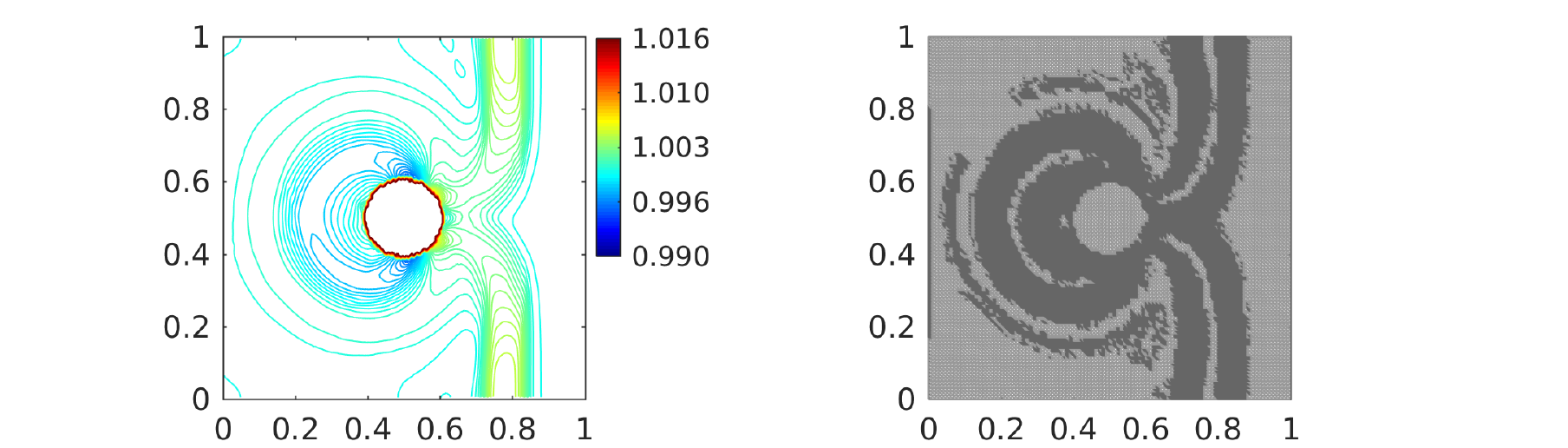}\\
		\vspace*{5mm}
		\caption{Example 2: $w$ component of the solution of
			the IVP \eref{eq:swe1}-\eref{eq:swe3},
			\eref{eq:ex2cB}-\eref{eq:ex2cI} at $t
			= 0.06, 0.1, 0.14,$ and $0.2$ (from top to bottom)
			obtained by the adaptive central-upwind scheme (left
			column)  and  the corresponding adaptive meshes with
			one level of refinement $\mathcal{M}=1$ (right column).}  \label{fig:ex2cr1l}
	\end{figure}

\begin{figure}[htbp]
		\centering
		\includegraphics[width=1\textwidth]{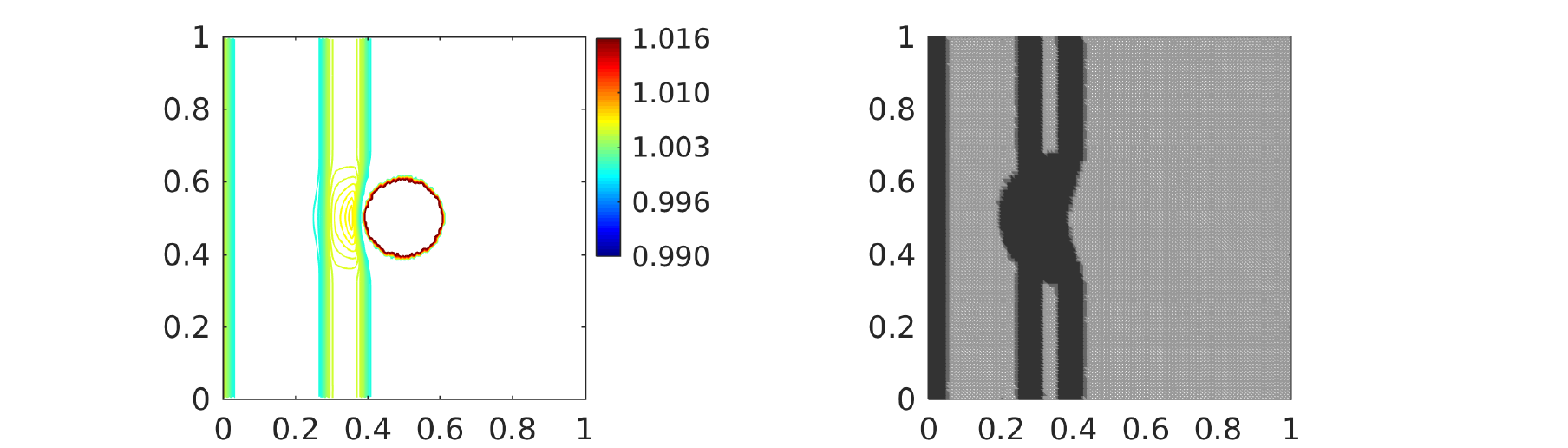}\\
		\includegraphics[width=1\textwidth]{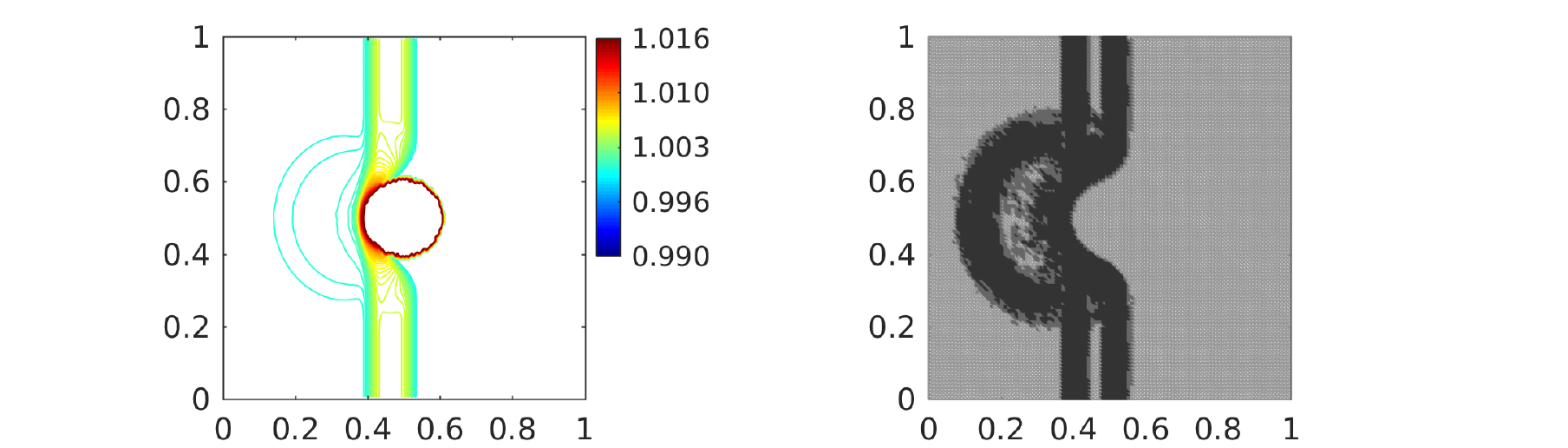}\\
		\includegraphics[width=1\textwidth]{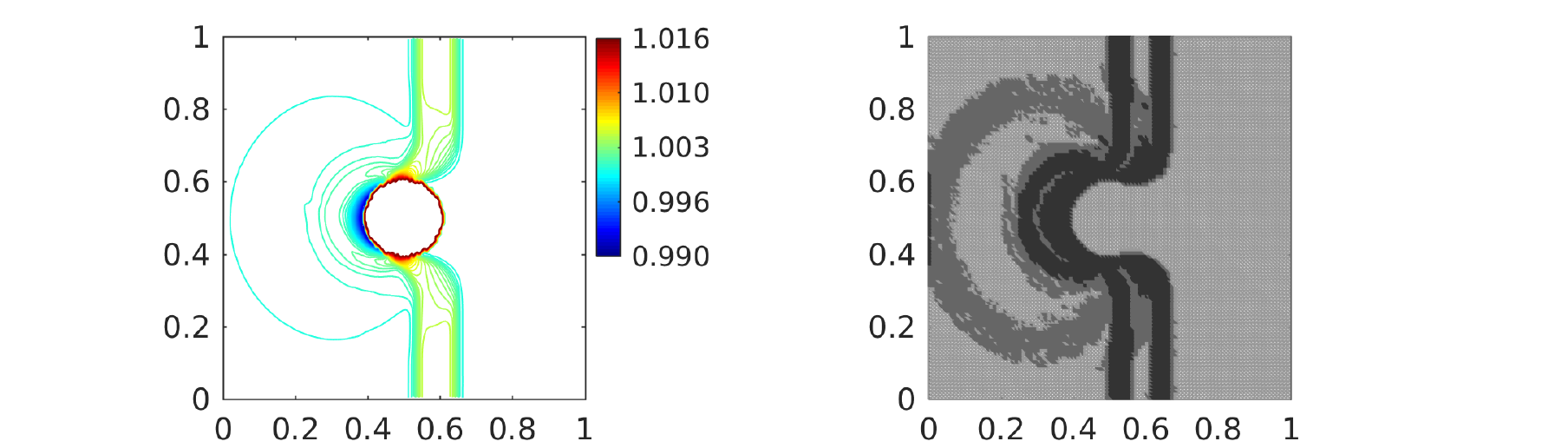}\\
		\includegraphics[width=1\textwidth]{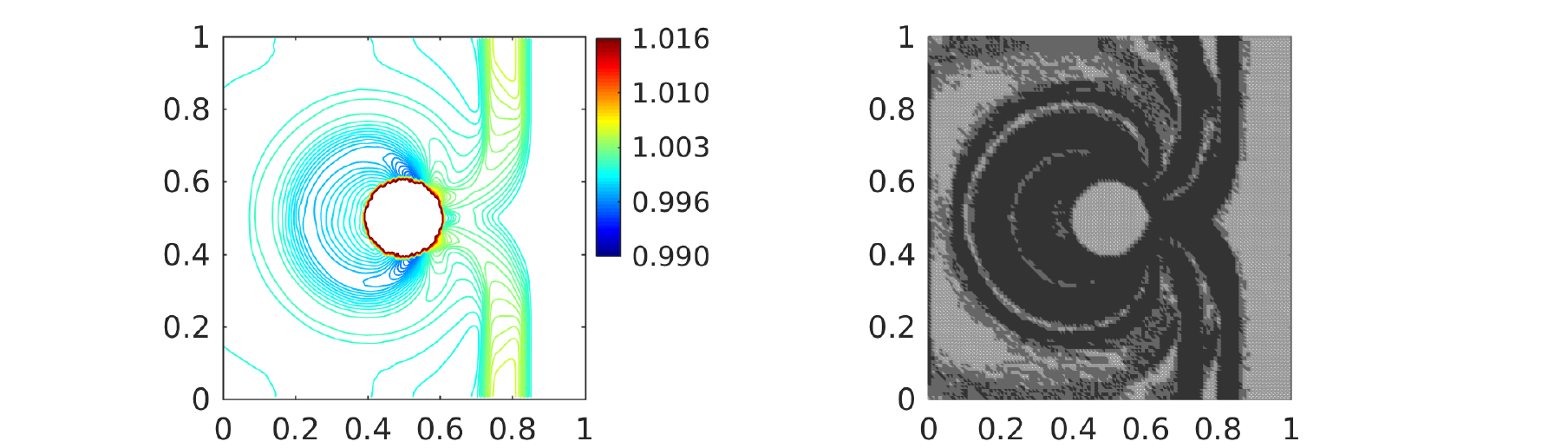}\\
		\vspace*{5mm}
		\caption{Example 2: $w$ component of the solution of
			the IVP \eref{eq:swe1}-\eref{eq:swe3},
			\eref{eq:ex2cB}-\eref{eq:ex2cI} at $t = 0.06, 0.1,
			0.14,$ and $0.2$ (from top to bottom)  obtained by the adaptive central-upwind scheme (left
			column)  and  the corresponding adaptive meshes with
			two levels of refinement $\mathcal{M}=2$ (right column).} \label{fig:ex2cr2l}
	\end{figure}

	We then compare with the performance of the adaptive scheme in
	this case. The adaptive grids are generated from the uniform
	mesh $2\times 100\times 100$ using the threshold
	$\omega=0.001\max_j(e_j)$ for one level of refinement
	$\mathcal{M}=1$ and $\omega=0.01\max_j(e_j)$ for $\mathcal{M}=2$. In \fref{fig:ex2cr1l} (left) and
	\fref{fig:ex2cr2l} (left), we plot the results for $w$ (left)
	obtained by the adaptive scheme and the corresponding adaptive
	meshes (right). The results are similar to the results of the
        central-upwind scheme without adaptivity in
	\fref{fig:ex2cu}. There are no non-physical spurious waves
	generated at the wet/dry front. The well-balanced property,
	positivity and stability of the proposed adaptive
	central-upwind method are maintained with the adaptive mesh
	refinement. As in previous examples, WLR error estimator
	accurately detects regions in the domain which are marked for
	adaptive refinement/coarsening.

	Next, in \tref{tab:ex2ccpu}, we present ratio of CPU times
	$\mathcal{R}_{CPU}$ to compute the solution at $t=0.1$ by the
	central-upwind method without adaptivity  and by the adaptive
	algorithm. 
	In order to compare the computational costs and calculate
	$\mathcal{R}_{CPU}$, we consider uniform and adaptive meshes with
	the same size of the smallest cells. From \tref{tab:ex2ccpu}, one can
	see that the adaptive algorithm reduces the CPU times up to four
	times. In our experiments, we considered only $\mathcal{M}=1$ and
	$\mathcal{M}=2$, but one can consider higher levels of refinement to
	further enhance the accuracy of the numerical solution at the reduced
	computational cost.

	\begin{table}[h!]
		\vspace*{2mm}
		\centering
		\begin{tabular}{ |c|c| c|c|c |}
			\hline
			\makecell{uniform mesh\\(cells)} &\makecell{adaptive mesh\\ $\mathcal{M}=1$\\(cells)}&\makecell{$\mathcal{R}_{CPU}$ \\with $\mathcal{M}=1$} &\makecell{adaptive mesh\\$\mathcal{M}=2$\\(cells)}&\makecell{$\mathcal{R}_{CPU}$ \\with $\mathcal{M}=2$}\\
			\hline
			$2\times100\times100$&11,831&1.91 &6,155& 3.04\\
			$2\times200\times200$&31,050 &2.08&25,753&  3.14 \\
			$2\times400\times400$&154,616&3.16&94,357 &5.82 \\
			\hline
			\multicolumn{2}{|c}{$\mathcal{R}_{CPU}$ average:}&
				\multicolumn{1}{c}{2.38 }&\multicolumn{1}{c}{} &\multicolumn{1}{c|}{4.00}\\
			\hline
		\end{tabular}
		\caption{Example 2: The $R_{CPU}$ ratio for solving the IVP  \eref{eq:swe1}-\eref{eq:swe3}, \eref{eq:ex2cB}-\eref{eq:ex2cI} at $t=0.1$.}\label{tab:ex2ccpu}
	\end{table}
	
\begin{figure}[htbp]
		\centering
		\includegraphics[width=0.3\textwidth]{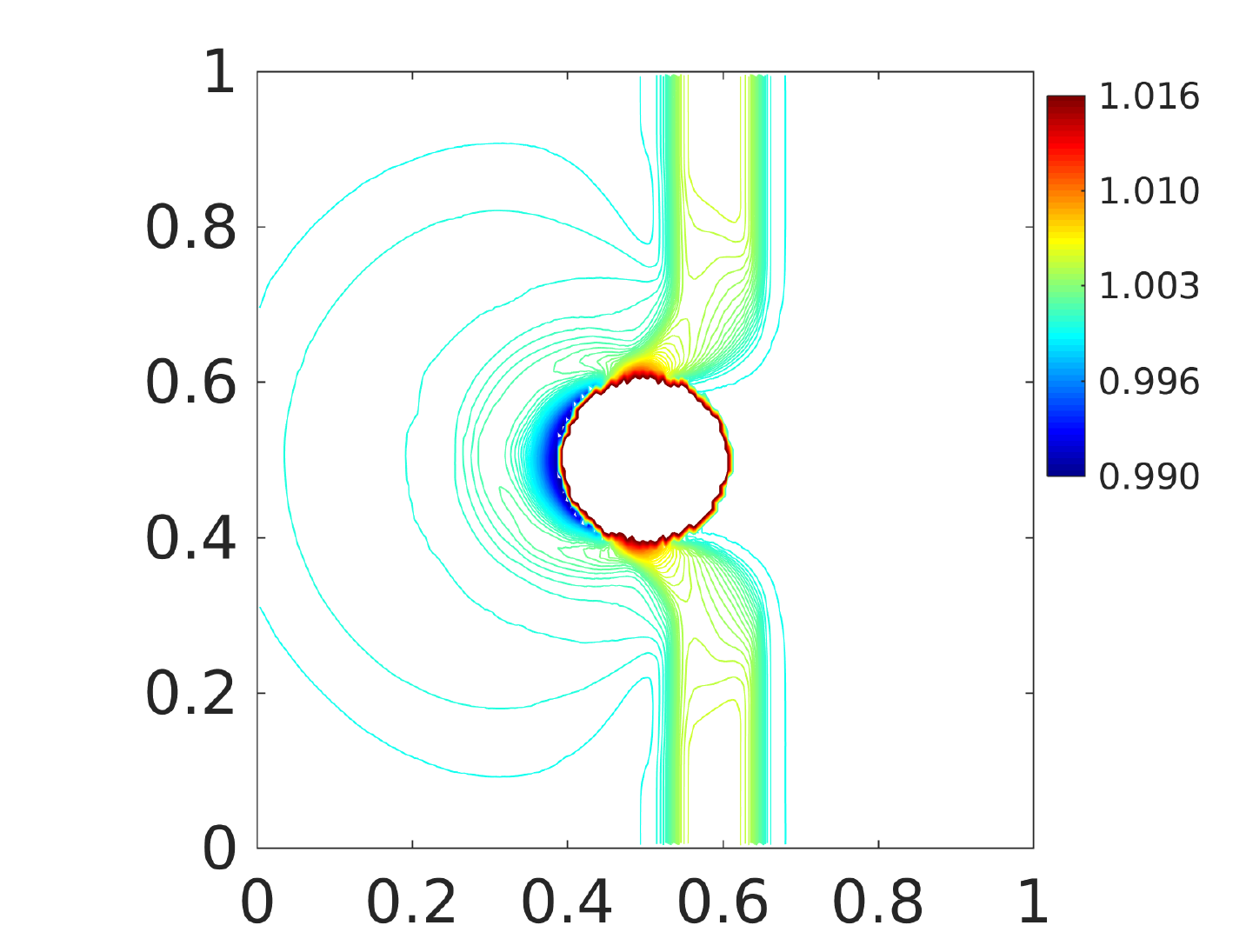}
		\includegraphics[width=0.3\textwidth]{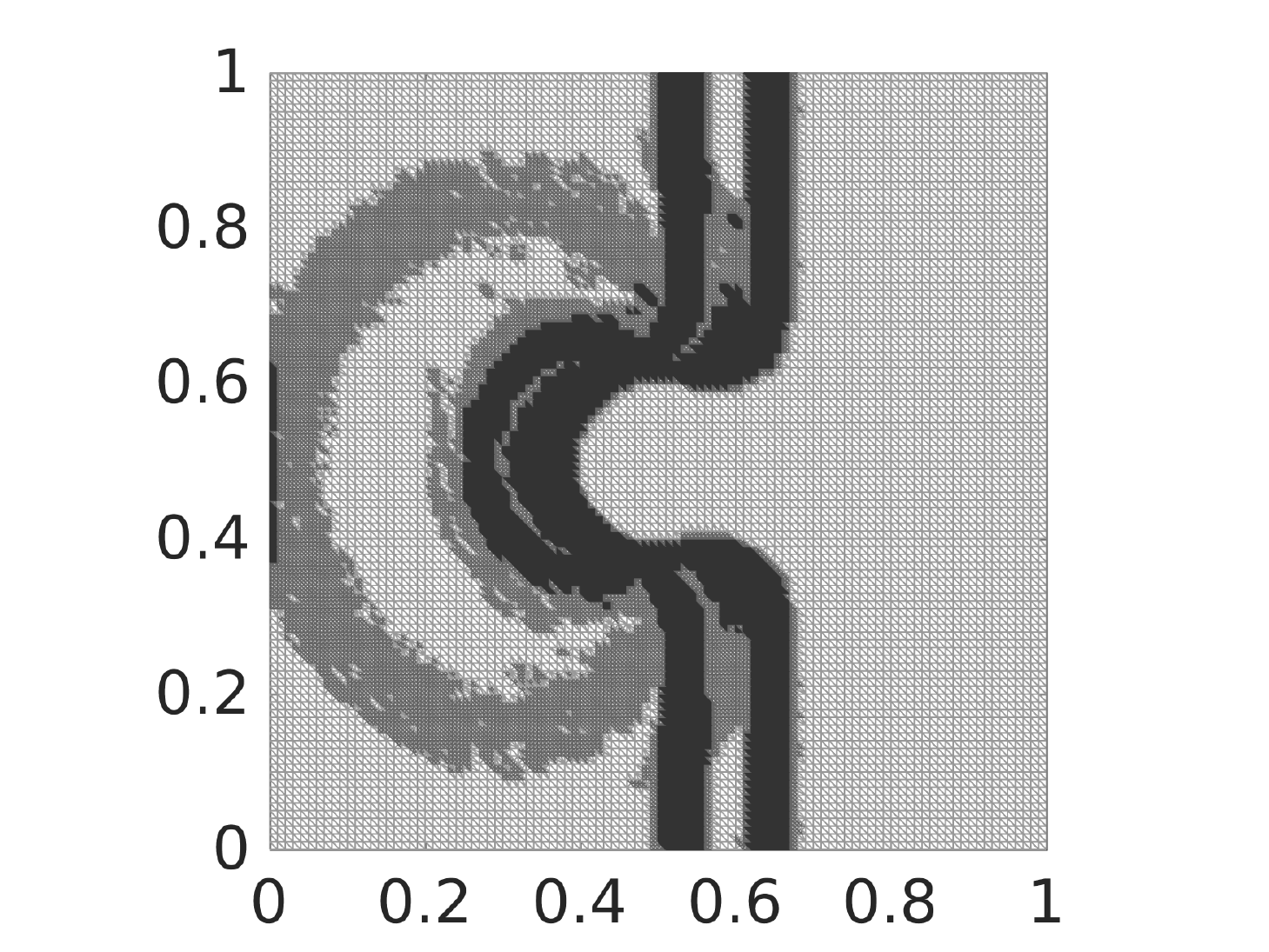}
		\includegraphics[width=0.3\textwidth]{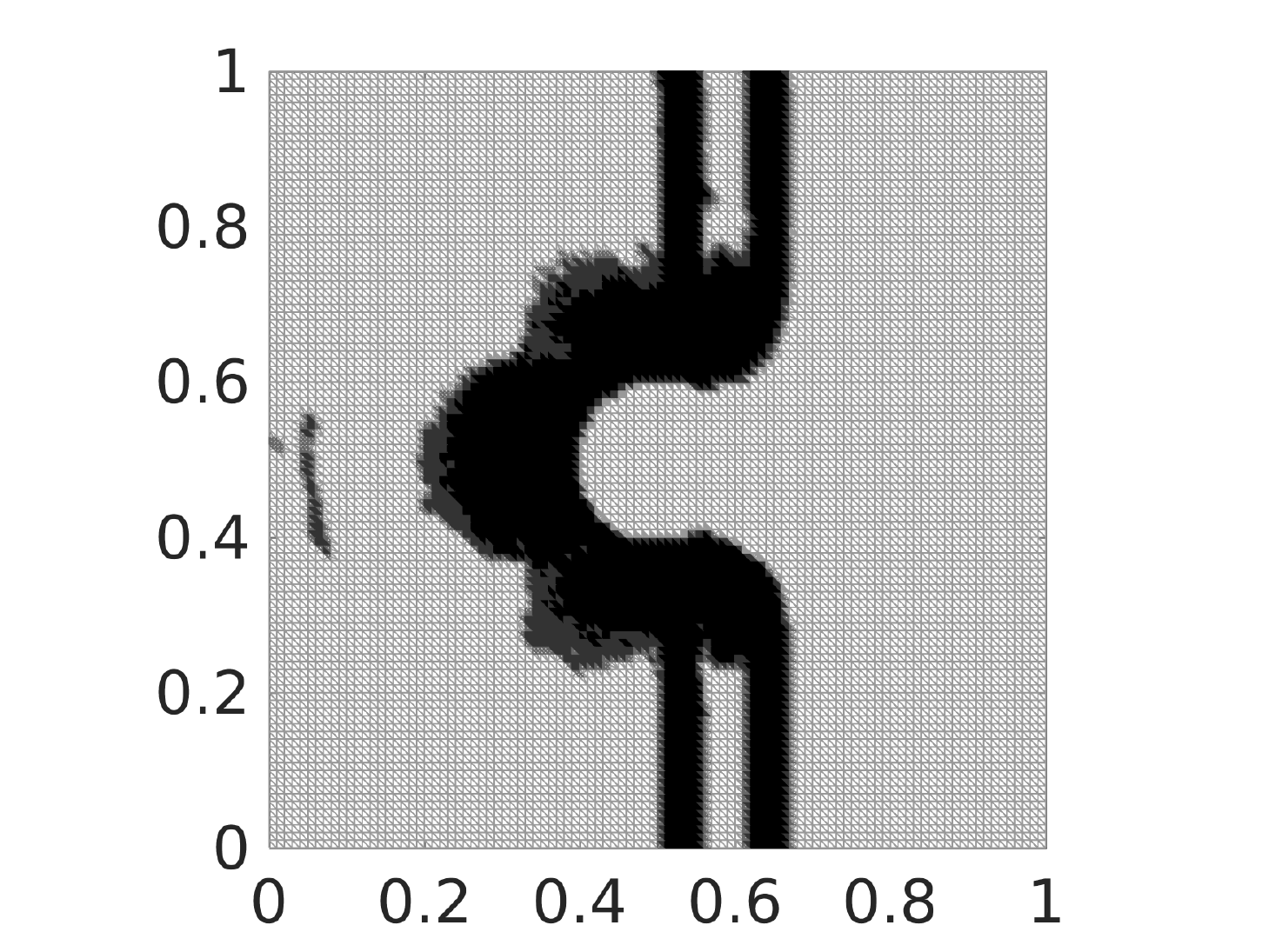}\\		\includegraphics[width=0.3\textwidth]{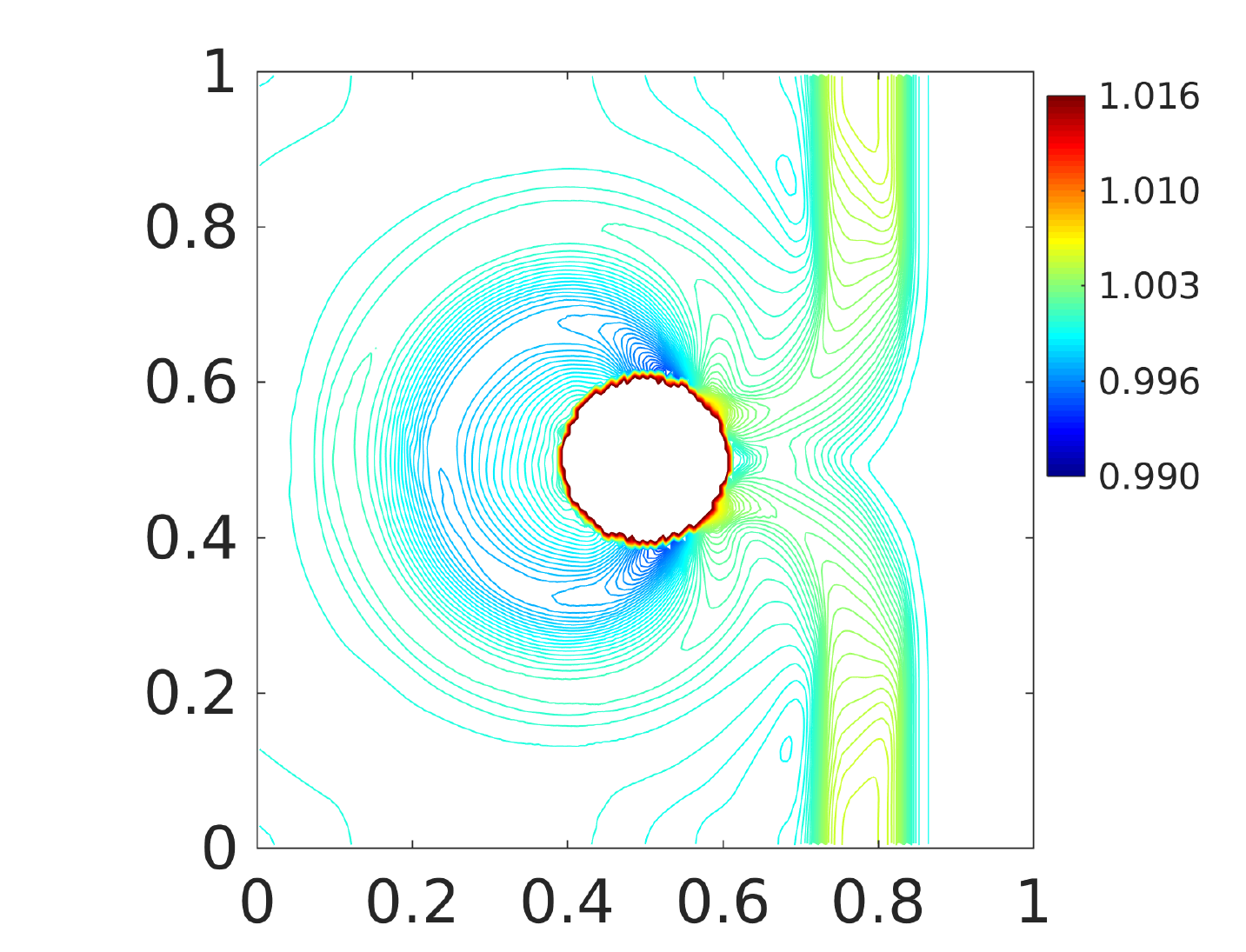}
		\includegraphics[width=0.3\textwidth]{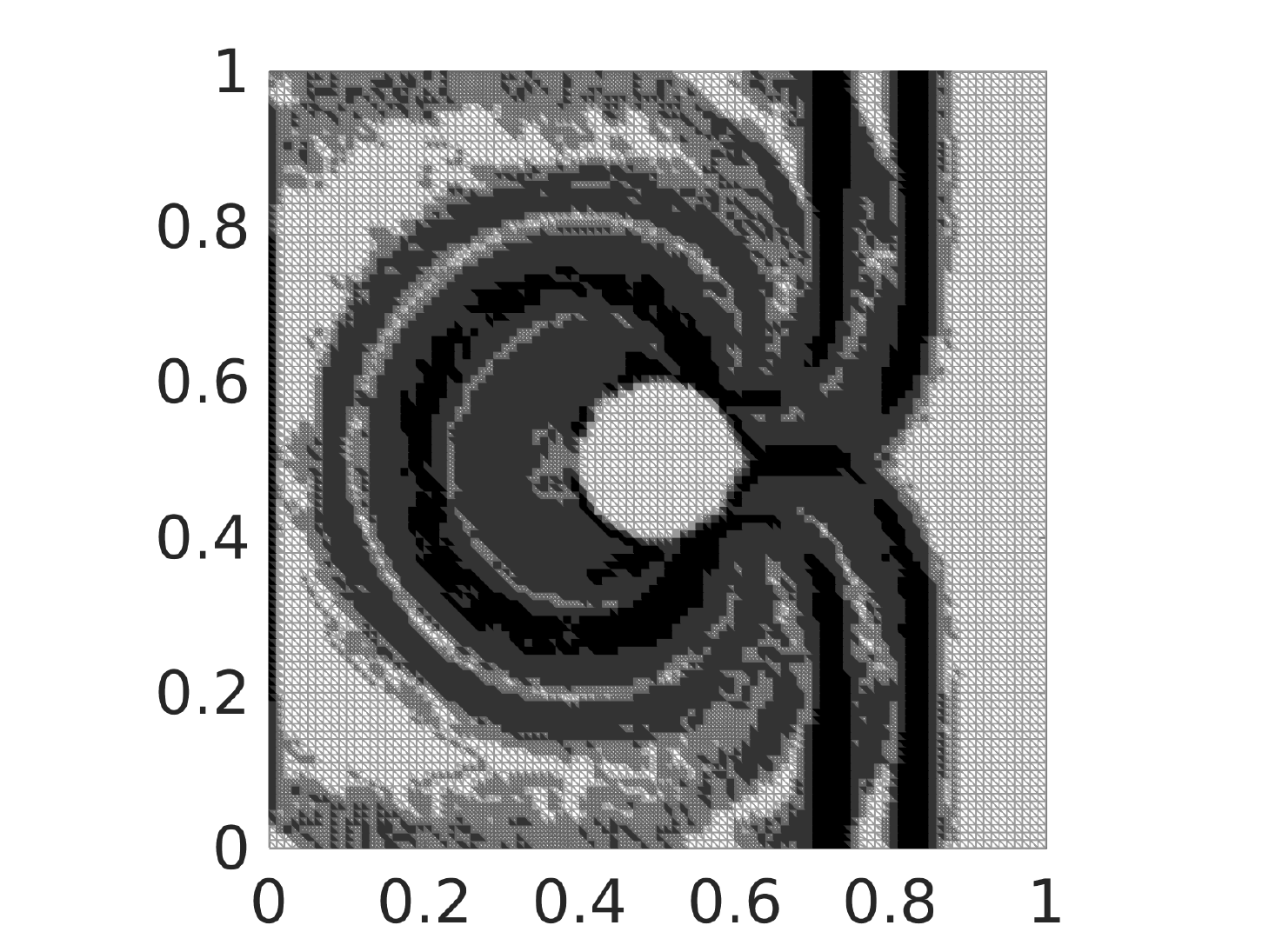}
		\includegraphics[width=0.3\textwidth]{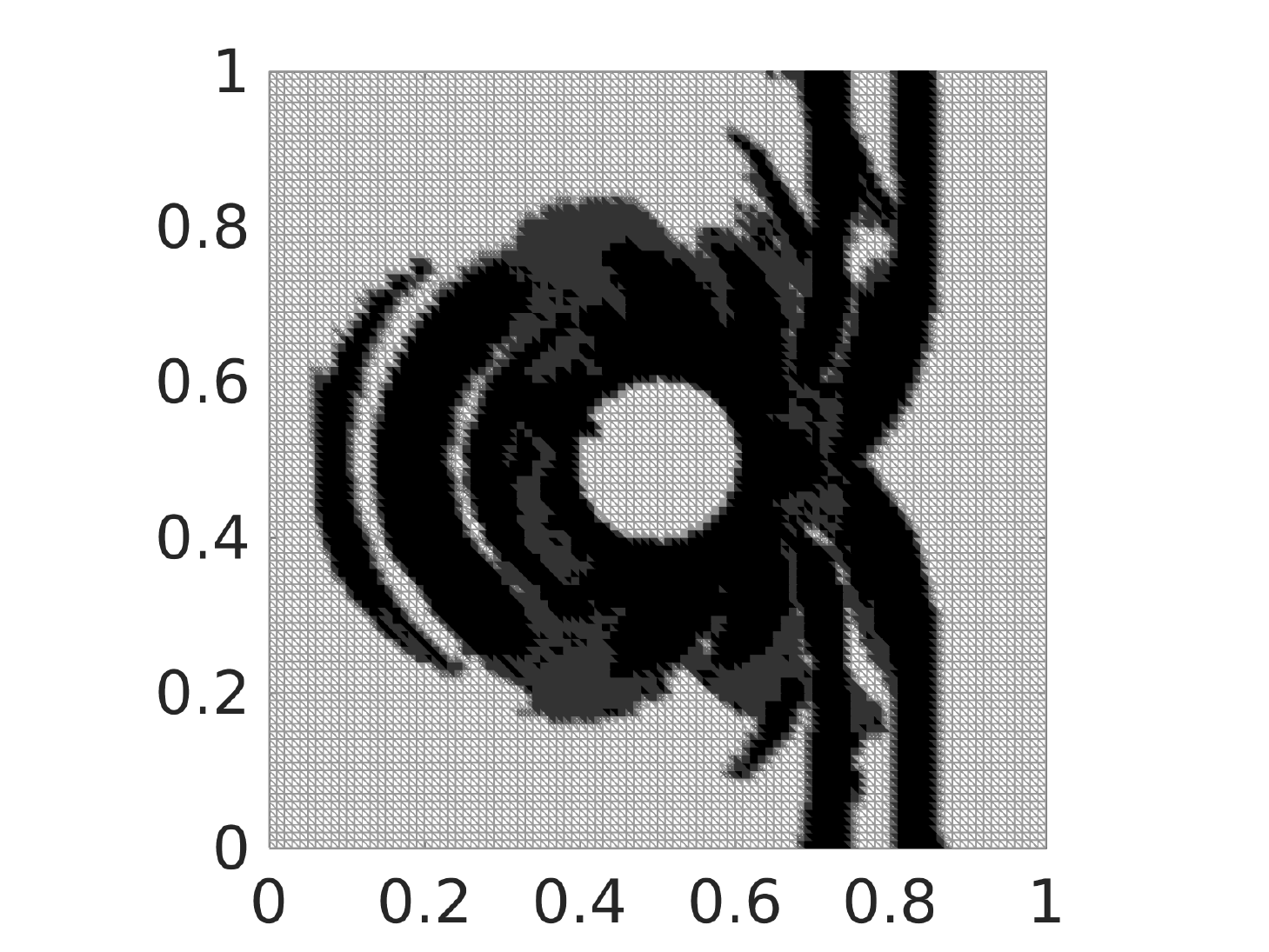}\\
		\vspace*{5mm}
		\caption{Example 2: $w$ component of the solution of
					\eref{eq:ex2cB}-\eref{eq:ex2cI}
                                        at $t =0.14$ (top) and $t=0.2$
                                        (bottom) obtained by the
                                        adaptive central-upwind scheme
                                        using WLR error indicator (left
			column). The corresponding adaptive
                        meshes reconstructed from uniform mesh
                        $2\times 100\times 100$ by using the WLR error
                        indicator (middle column) and the gradient indicator (right column).
		} \label{fig:ex2un}
	\end{figure}

	  Finally,  we illustrate the advantages
                  of WLR error as the error indicator,  and hence
                  compare it with another example of the error
                  indicator which uses the unlimited gradients of the
                  water surface $(w_x,w_y)$, see
                  e.g. \citep{GHAZIZADEH2020104633}. We continue to
                  consider the IVP
                  \eref{eq:swe1}-\eref{eq:swe3},\eref{eq:ex2B}-\eref{eq:ex2I}. In
                  \fref{fig:ex2un} (left column), we first present the
                  contour plots of water surface $w$ computed at
                  $t=0.14$ and $t=0.2$ using the meshes obtained by
                  WLR error indicator,   \fref{fig:ex2un} (middle column). The adaptive meshes in
                  \fref{fig:ex2un} (middle and right columns) are
                  reconstructed from the uniform mesh $2\times
                  100\times 100$, respectively,  by using the WLR error
                  indicator and the gradient indicator. For WLR error,
                  we apply the threshold $\omega=0.01\max_j(e_j)$ and
                  refine a cell $T_j$ if local WLR error
                  $e_j>\omega\times2^{m-1}$ for $m\leq\mathcal{M}=3$
                  levels. For the gradient indicator,  a cell $T_j$ is
                  flagged for $m, m\leq3,$ levels of refinement if it
                  has $(w_x)_j^2+(w_y)_j^2>0.0005\times2^{m-1}$. From
                  \fref{fig:ex2un}, we can observe that the WLR error indicator
                  captures more subtle features of the solution than
                  the error indicator based on the gradient. 
		
	\subsection{Example 3---Dam Break Test}\label{sec4_3}
	
	In example 3 taken from \citep{LAEK}, we simulate the
	propagation of the dam-break flood wave which produces a
	moving wet/dry fronts over an irregular dry bed with three
	obstacles. The test allows to verify capability of the proposed
	adaptive algorithm to handle wet/dry interfaces. The bottom topography is defined by, 	
	\begin{equation}
	B(x,y)=\max\left[0.5e^{-8(x-2)^2-10(y-3)^2},0.2e^{-3(x-4)^2-4(y-4.8)^2},0.2e^{-3(x-4)^2-4(y-1.2)^2}\right],
	\end{equation}
	in the computational domain is $[0, 6]\times [0, 6]$. At $t=0$, an upstream reservoir in the region $[0,1]\times [0,6]$ filled with water up to $w(x, y, 0)=0.5$ is suddenly released. Hence, the following initial condition is imposed,
	
	\begin{equation}
	w(x,y,0)=\begin{cases}
	0.5, \quad& 0\leq x<1,\\ B(x,y), &\mbox{otherwise,}
	\end{cases}\quad u(x,y,0)\equiv v(x,y,0)\equiv 0.
	\end{equation}
	
	In this example, we consider the friction effects by modifying the governing equation \eref{eq:swe2} and \eref{eq:swe3} with the Manning friction terms as follows,	
	\begin{equation}
	\begin{split}
	(hu)_t+\left(hu^2+\dfrac{g}{2}h^2\right)_x+(huv)_y=-ghB_x-\dfrac{gn^2_bu\sqrt{u^2+v^2}}{h^{1/3}},\\
	(hv)_t+(huv)_x+\left(hv^2+\dfrac{g}{2}h^2\right)_y=-ghB_y-\dfrac{gn^2_bv\sqrt{u^2+v^2}}{h^{1/3}},
	\end{split}
	\end{equation}	
	where $n_b=0.01$ is the Manning roughness coefficient. We have used the homogeneous Neumann boundary condition for the left boundary $x=6$ and the solid wall boundary condition for the other boundaries. 
	
\begin{figure}[h!]
		\centering
		\centering
		\includegraphics[width=1\textwidth]{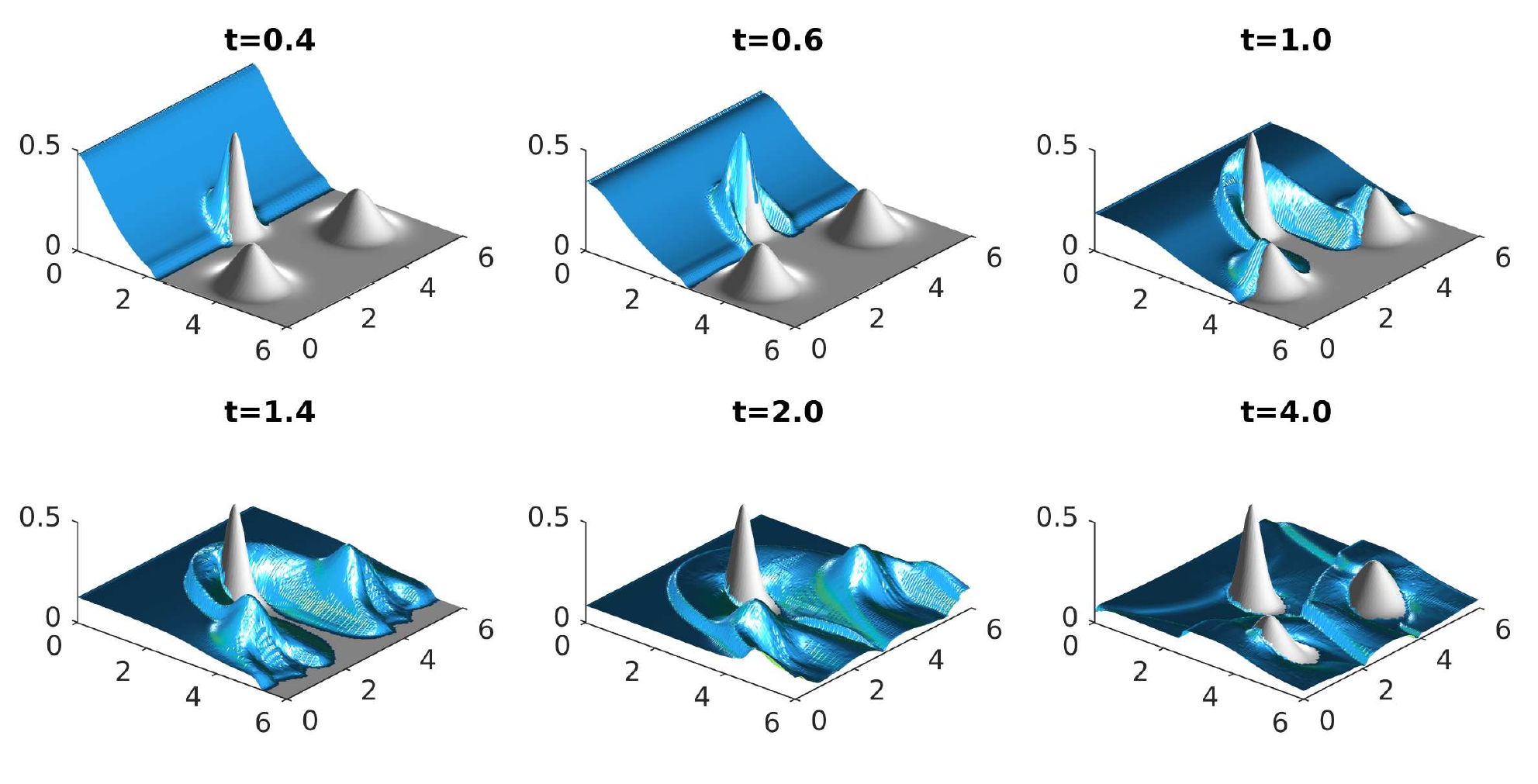}
		\caption{Example 3: Simulated water surface $w$ at
			different times on uniform grid $2\times 100\times 100$.}\label{fig:ex33du1}
		\vspace*{1.mm}
	\end{figure}

\begin{figure}[h!]
		\centering
		\centering
		\includegraphics[width=1\textwidth]{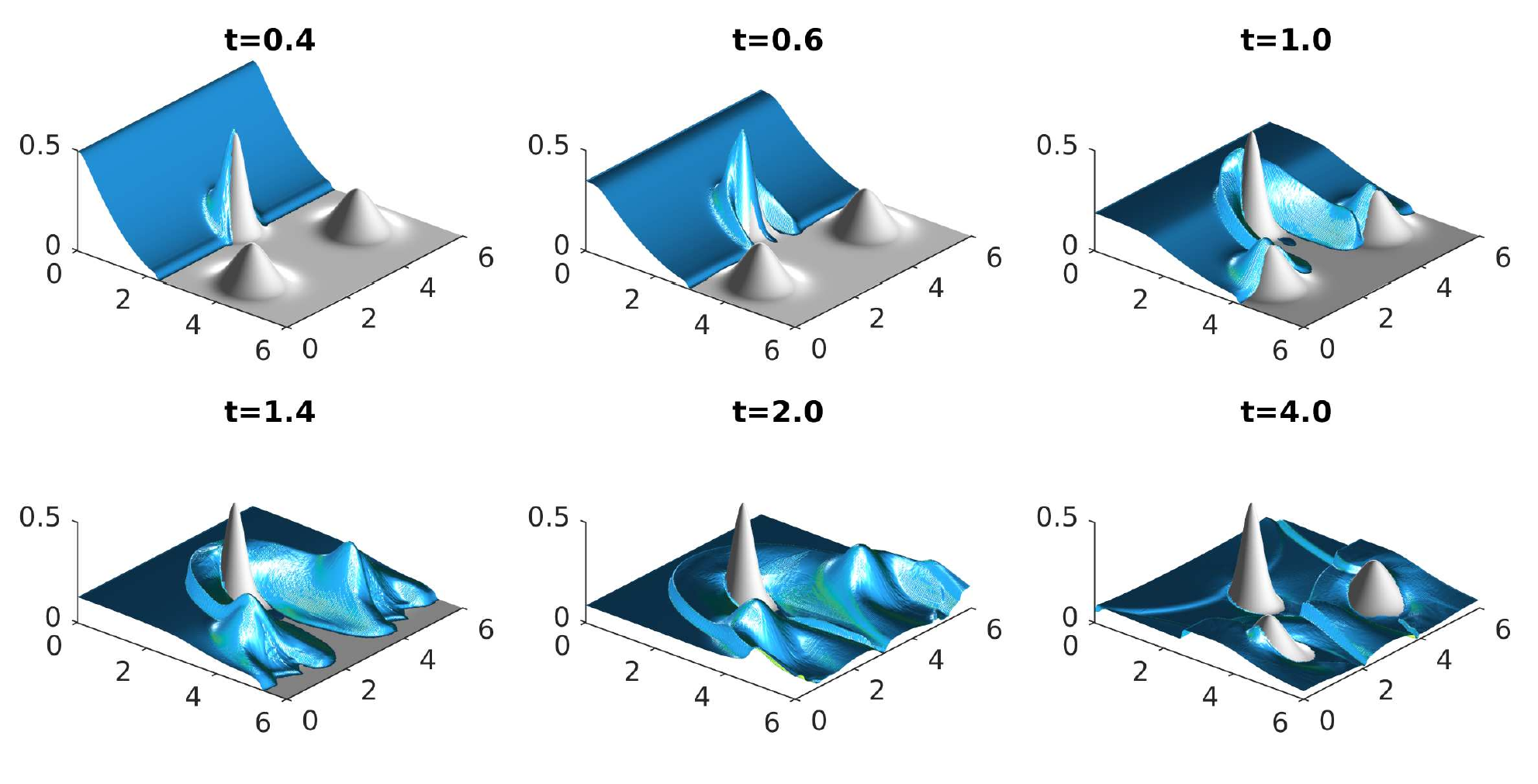}
		\caption{Example 3: Simulated water surface at
			different times on uniform grid $2\times 200\times 200$.}\label{fig:ex33du2}
	\end{figure}

\begin{figure}[h!]
		\centering
		\centering
		\includegraphics[width=1\textwidth]{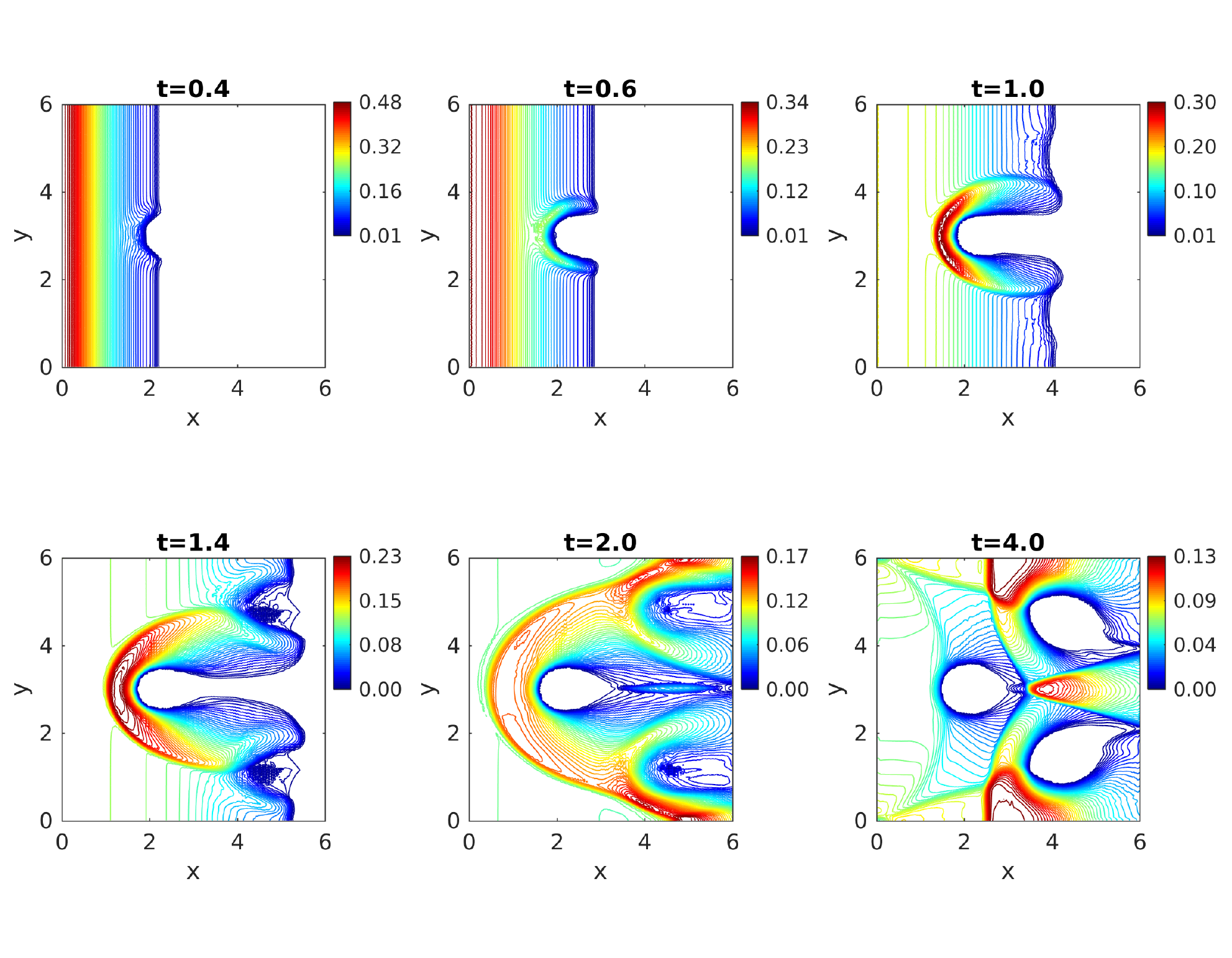}\\
		\vspace*{-5mm}
		\caption{Example 3: Contour of the water depth $h$ at
			different times on uniform grid $2\times 100\times 100$.}\label{fig:ex33du1c}
	\end{figure}

\begin{figure}[h!]
		\centering
		\centering
		\includegraphics[width=1\textwidth]{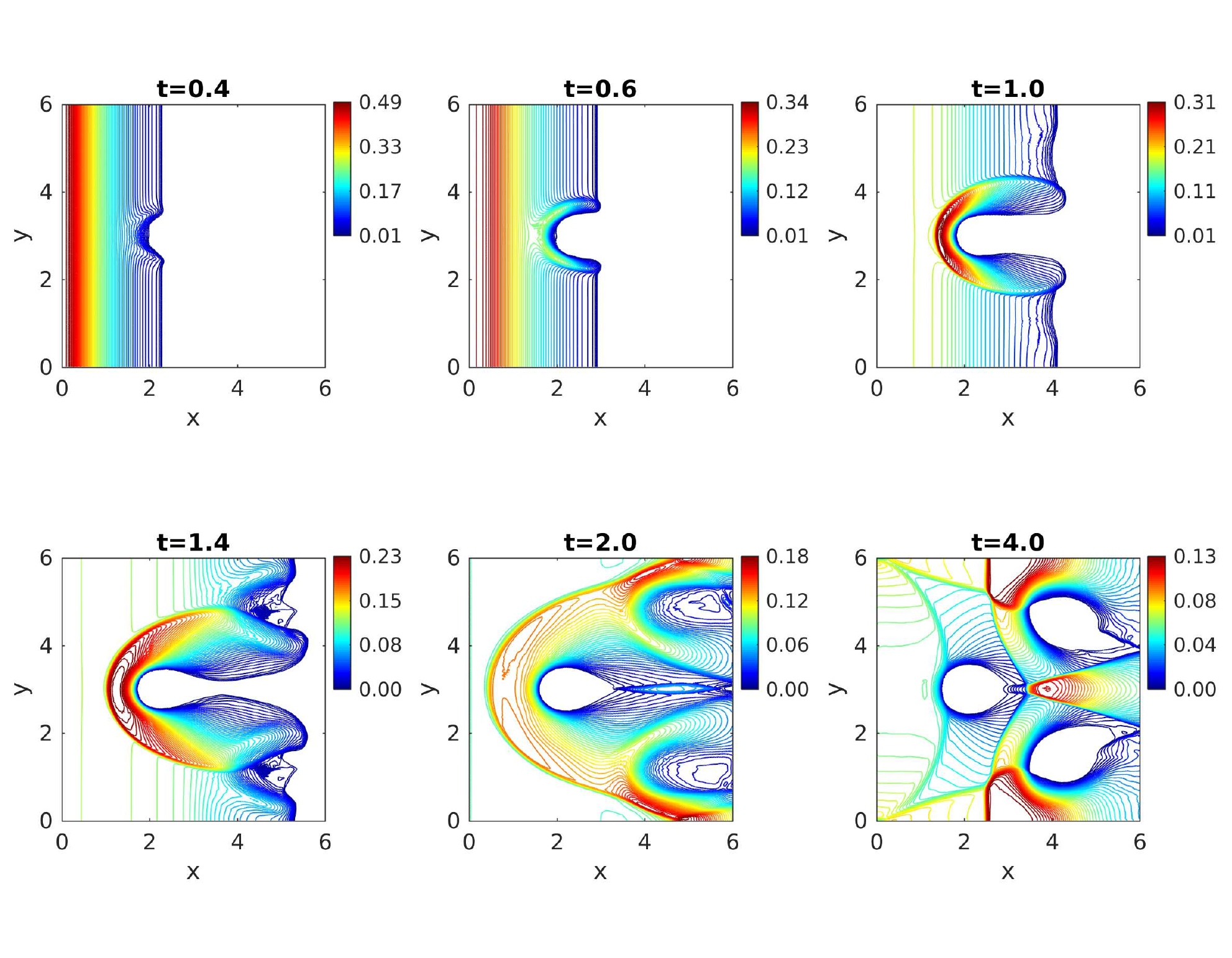}\\
		\vspace*{-5mm}
		\caption{Example 3: Contour of the water depth $h$ at
			different times on uniform grid $2\times 200\times 200$.}\label{fig:ex33du2c}
		\vspace*{2mm}
	\end{figure}

	We first present the numerical solution computed by applying
	the central-upwind scheme without the adaptivity \citep{LAEK}
	on the uniform meshes. \fref{fig:ex33du1} and
	\fref{fig:ex33du2} are the 3-D view of the dam-break wave
	propagation over the initially dry bed obtained respectively
	on $2\times 100\times 100$ and  $2\times 200\times 200$
	uniform meshes at different times $t=0.4, 0.6, 1.0, 1.4, 2.0,$
	and $4.0$. As can be observed from the figures, the water wave
	spreads from the reservoir and passes the obstacles. In
	addition, we plot the contour lines of the water depth for the
	solutions obtained on the uniform meshes $2\times 100\times
	100$ (\fref{fig:ex33du1c}) and $2\times 200\times 200$
	(\fref{fig:ex33du2c}). The contour lines clearly show the
	reflections and interactions of waves.

\begin{figure}[h!]
		\centering
		\centering
		\includegraphics[width=1\textwidth]{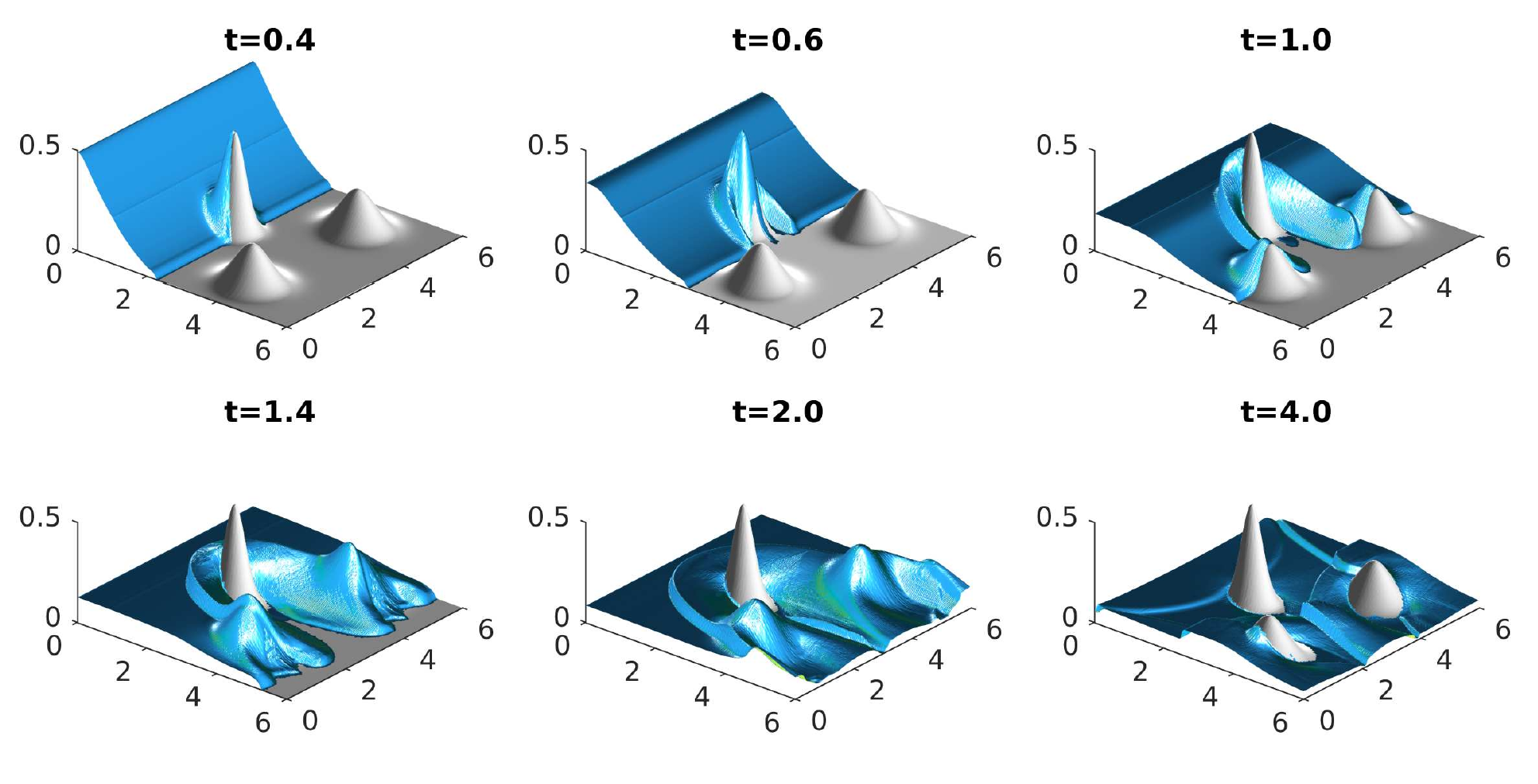}
		\caption{Example 3: Simulated water surface $w$ at different times on the adaptive mesh with $\mathcal{M}=1$.}\label{fig:ex33dr1l}
	\end{figure}

\begin{figure}[h!]
		\centering
		\vspace*{5mm}\includegraphics[width=1\textwidth]{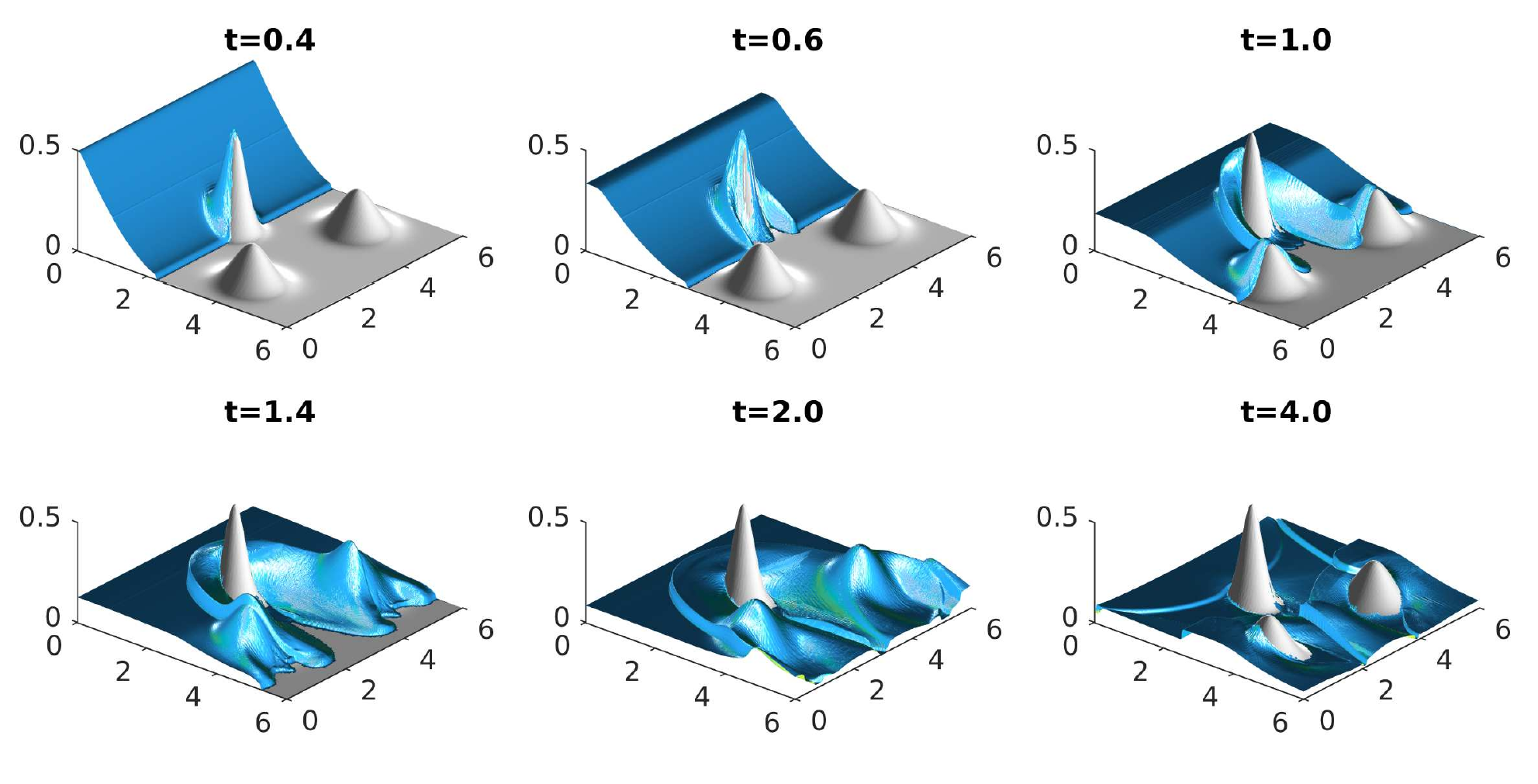}
		\caption{Example 3:  Simulated water surface $w$ at different times on the adaptive mesh with $\mathcal{M}=2$.}\label{fig:ex33dr2l}
		\vspace*{2mm}	
	\end{figure}

\begin{figure}[h!]
		\centering
		\centering
		\includegraphics[width=1\textwidth]{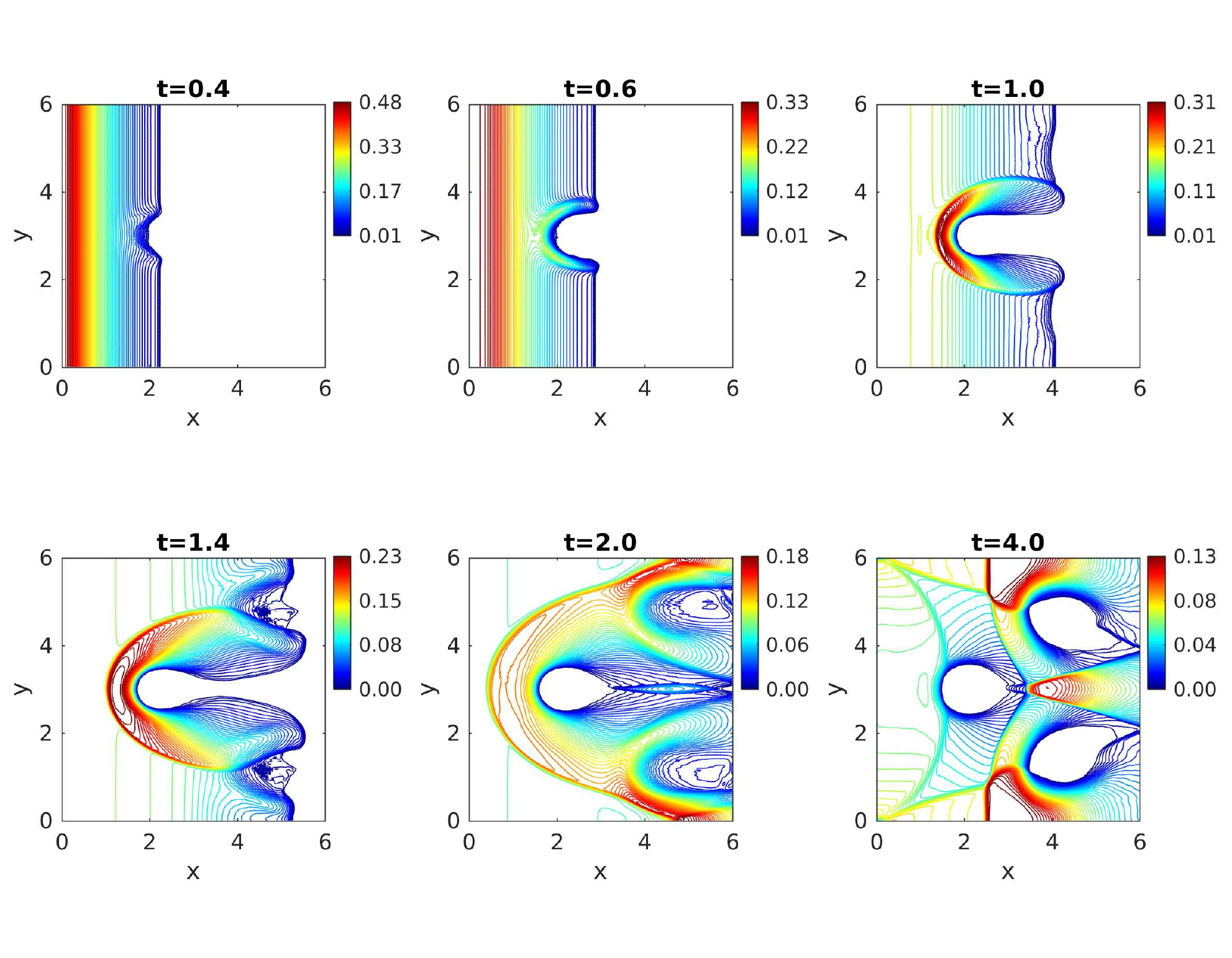}\\
		\vspace*{-5mm}
		\caption{Example 3: Contour of the water depth $h$ at different times on the adaptive mesh with $\mathcal{M}=1$.}\label{fig:ex3r1lc}
		\vspace*{2mm}	
	\end{figure}

\begin{figure}[h!]
		\centering
		\centering
		\includegraphics[width=1\textwidth]{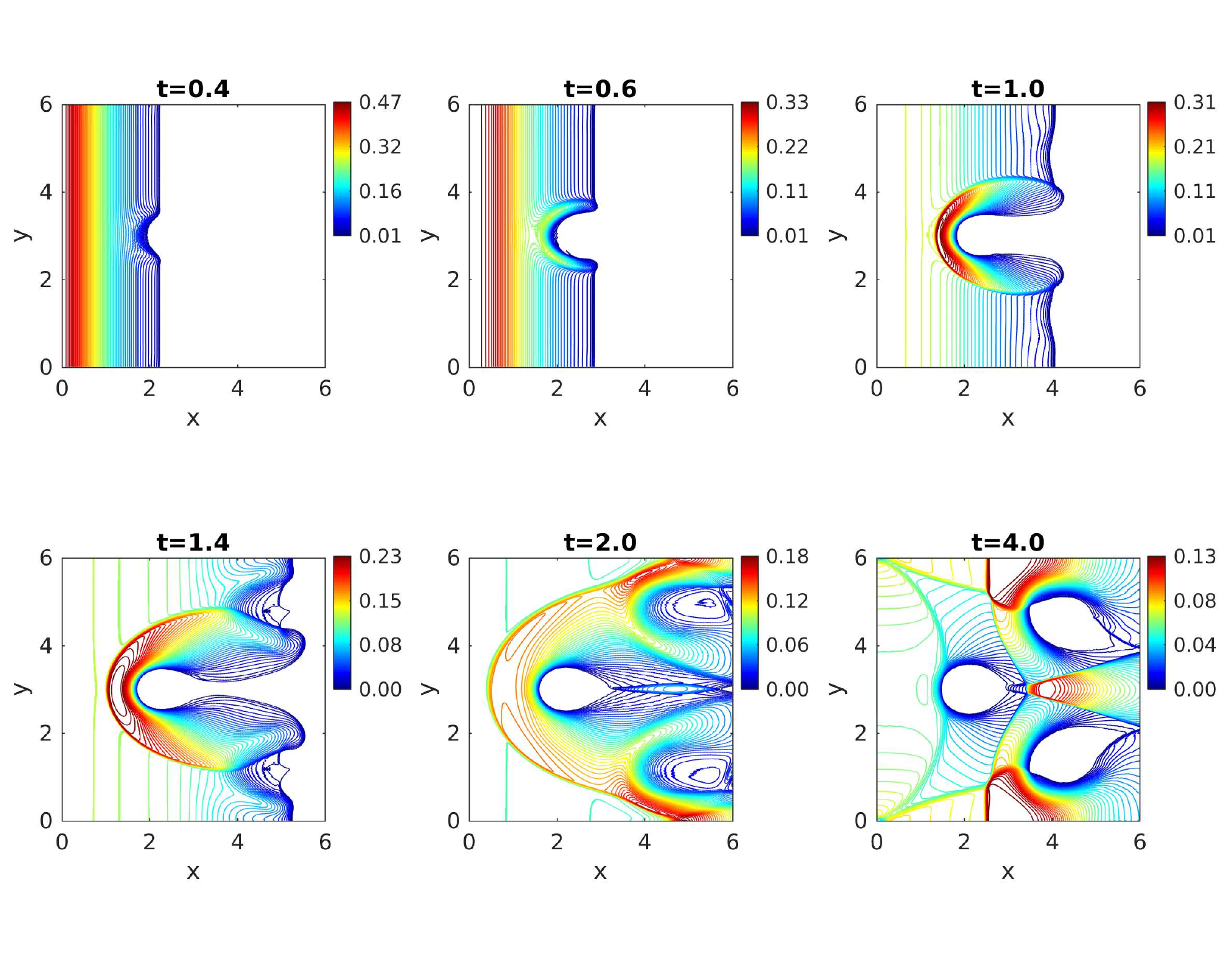}\\
		\vspace*{-5mm}
		\caption{Example 3: Contour of the water depth $h$ at different times on the adaptive mesh with $\mathcal{M}=2$.}\label{fig:ex3r2lc}
		\vspace{2mm}
	\end{figure}

\begin{figure}[h!]
		\centering
		\includegraphics[width=0.255\textwidth]{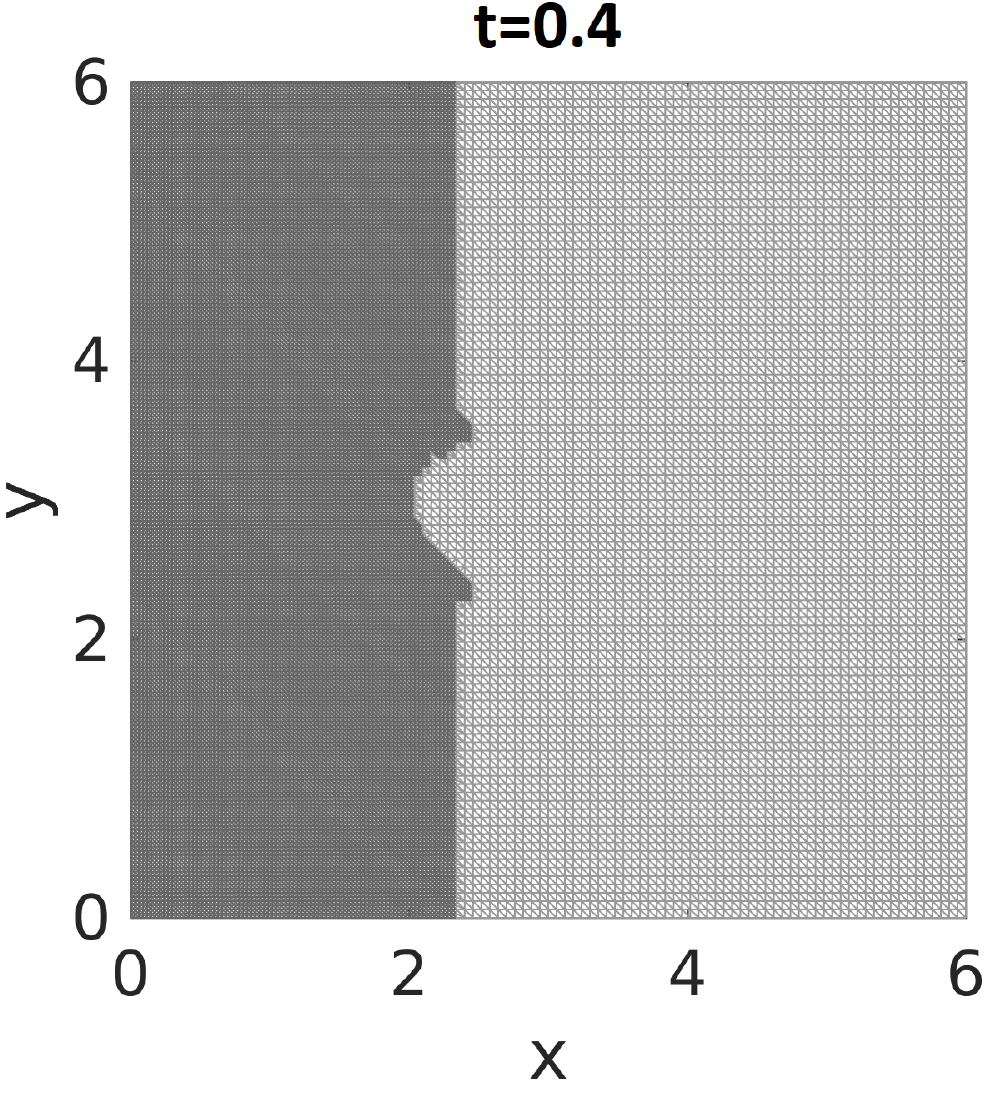}\hspace{1cm}
		\includegraphics[width=0.255\textwidth]{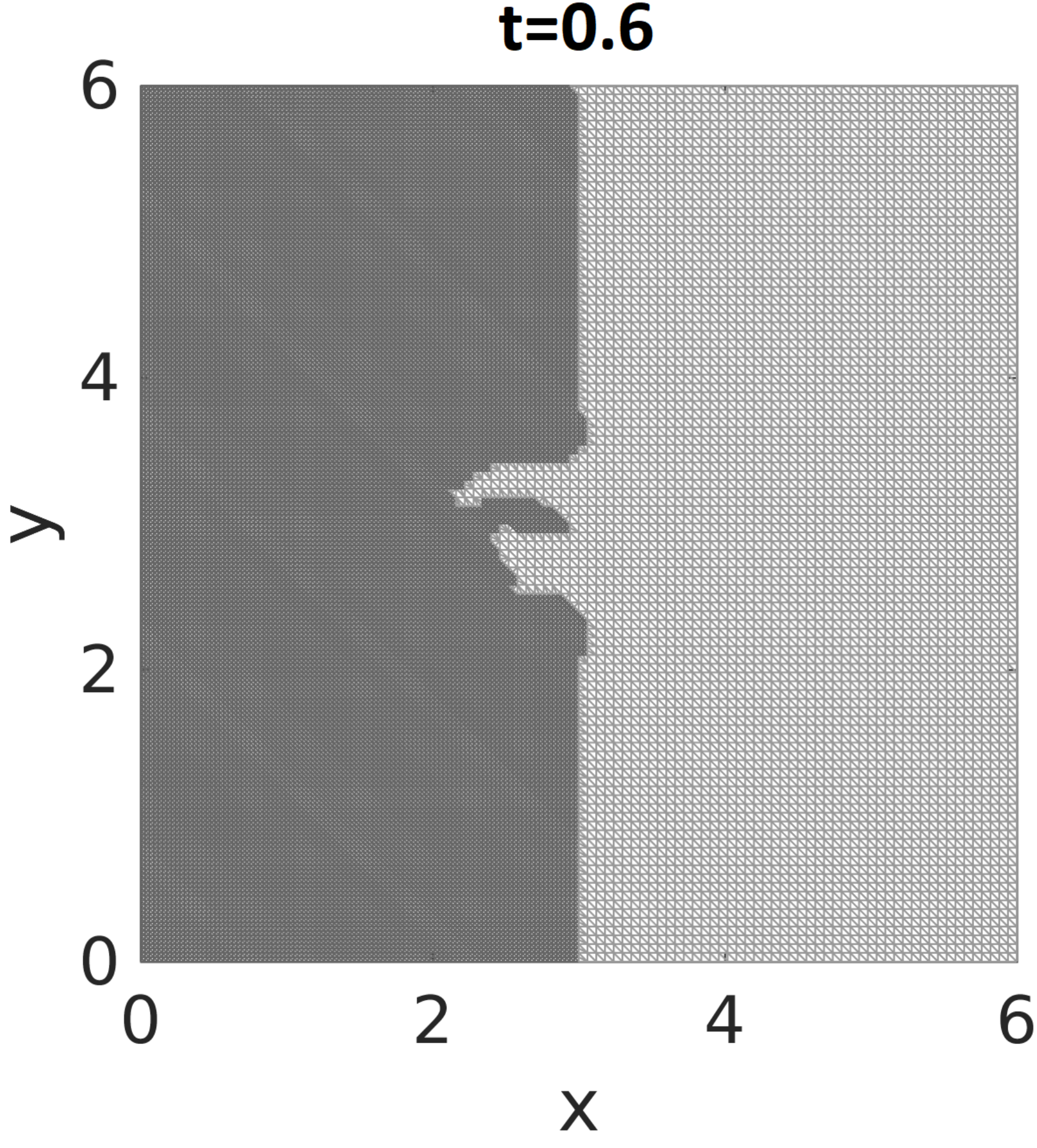}\hspace{1cm}
		\includegraphics[width=0.255\textwidth]{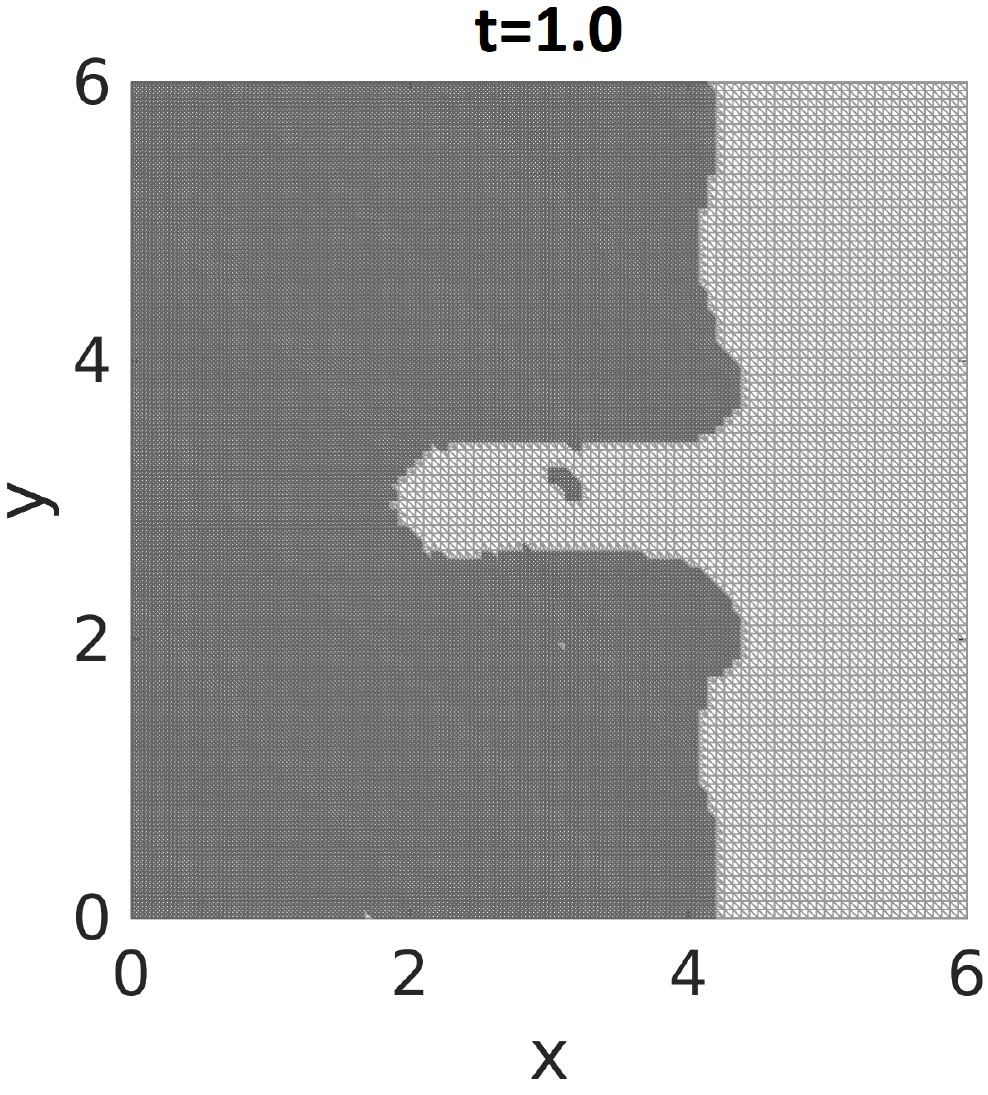}\\
		\vspace*{5mm}
		\includegraphics[width=0.255\textwidth]{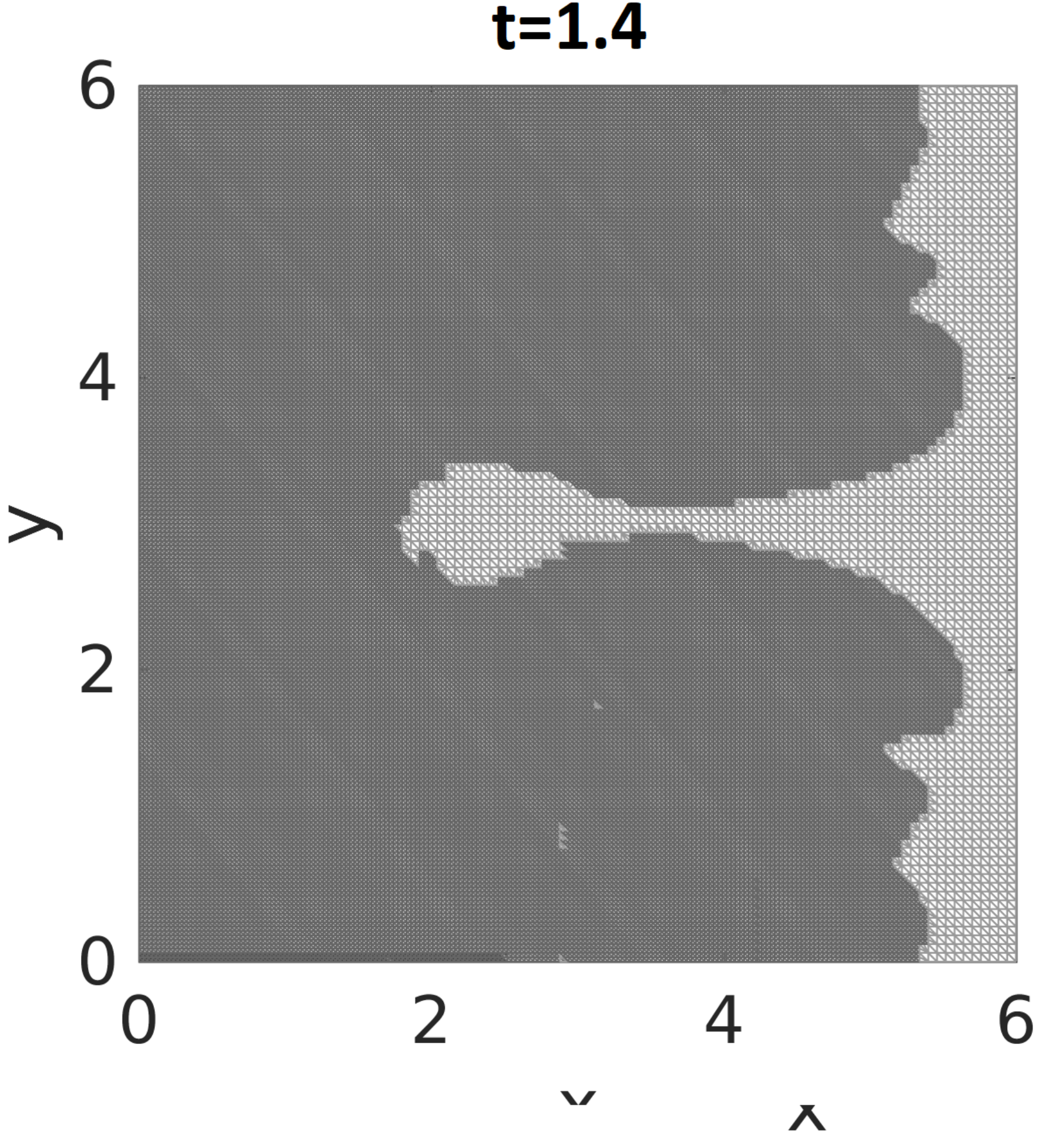}\hspace{1cm}
		\includegraphics[width=0.255\textwidth]{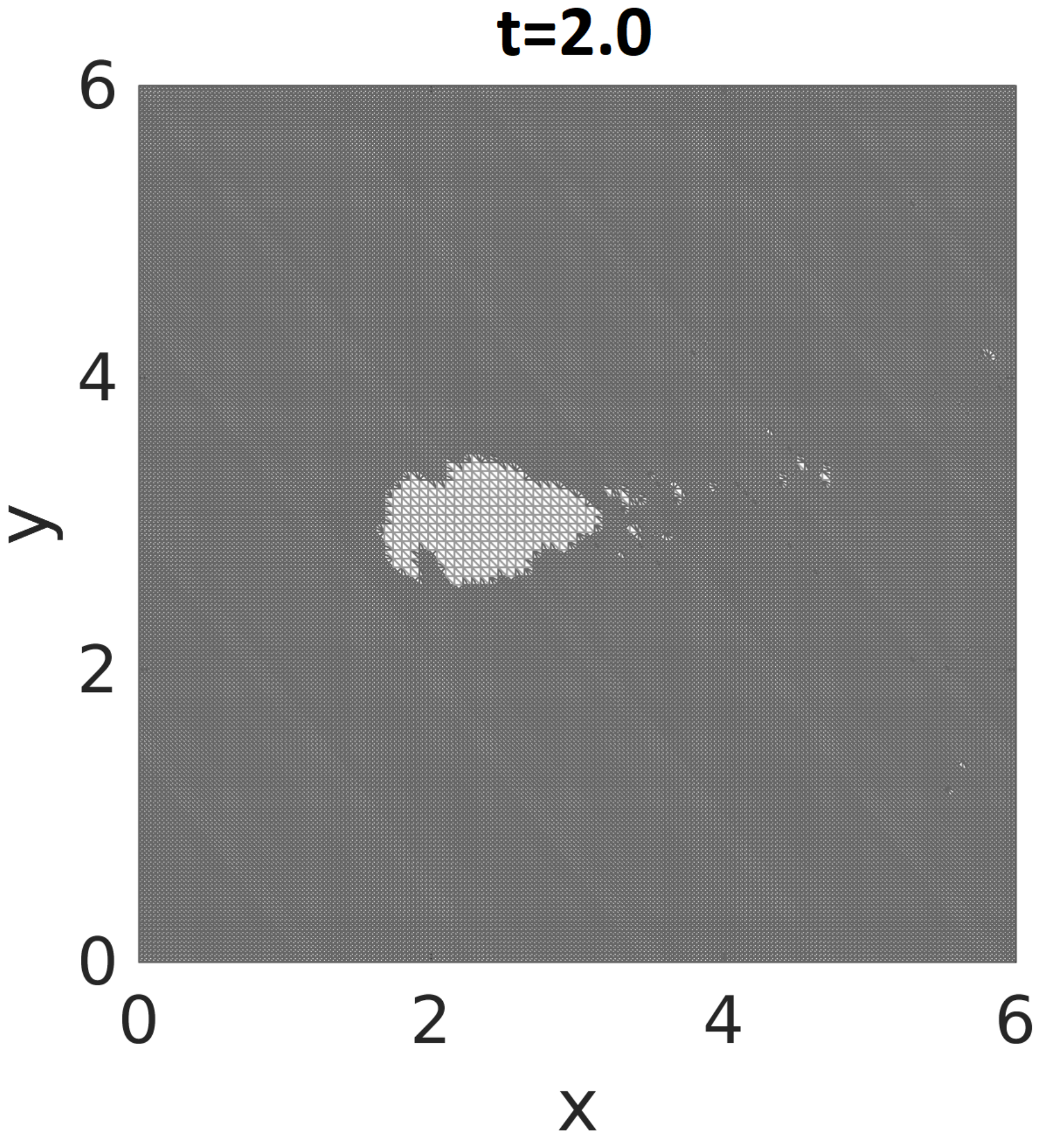}\hspace{1cm}
		\includegraphics[width=0.255\textwidth]{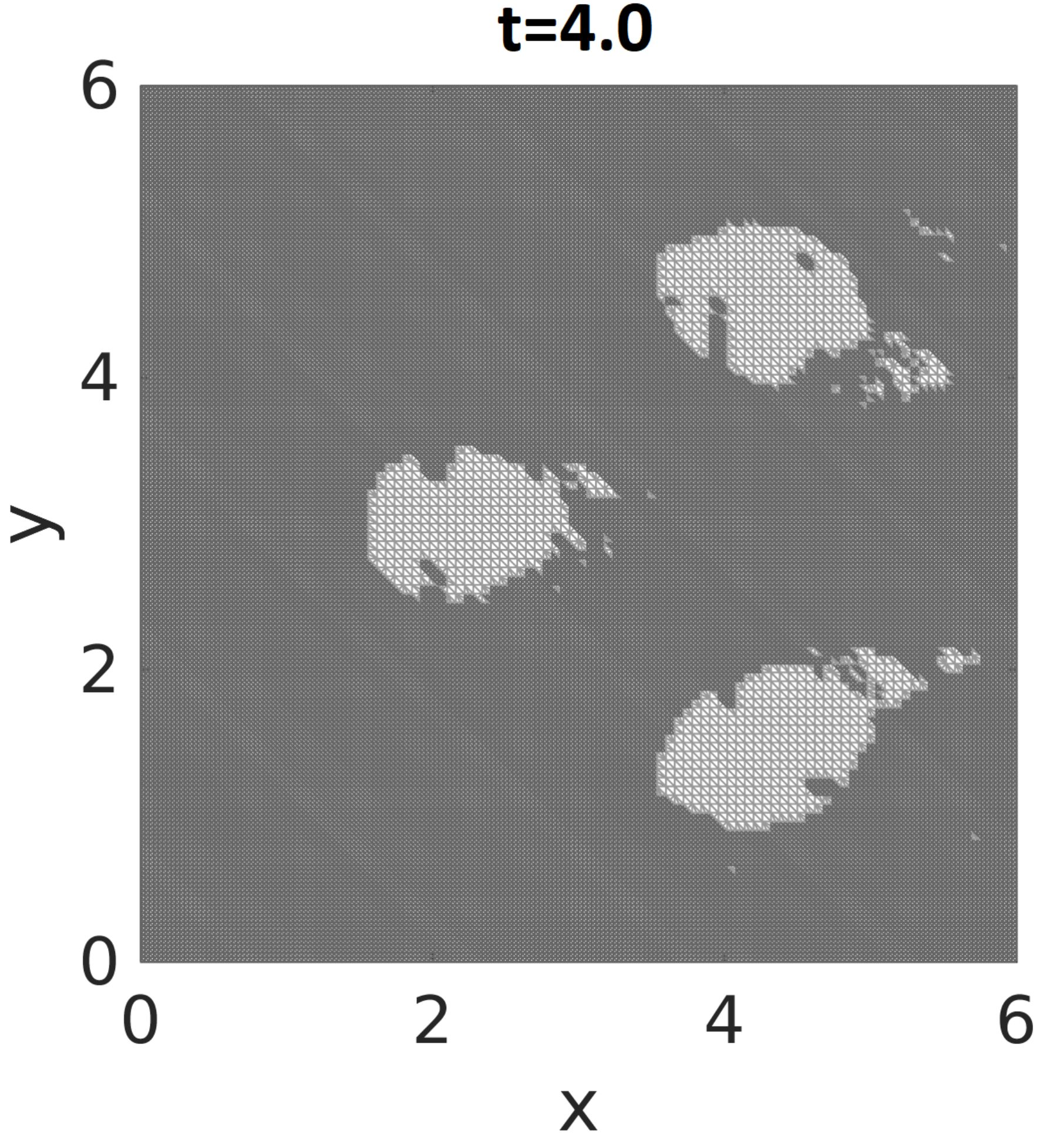}\\
		\vspace*{5mm}
		\caption{Example 3:  Adaptive mesh at different times with one level of refinement $\mathcal{M}=1$.}\label{fig:ex3r1lm}
		\vspace*{2mm}
	\end{figure}

\begin{figure}[h!]
		\centering
		\vspace*{5mm}
		\includegraphics[width=0.255\textwidth]{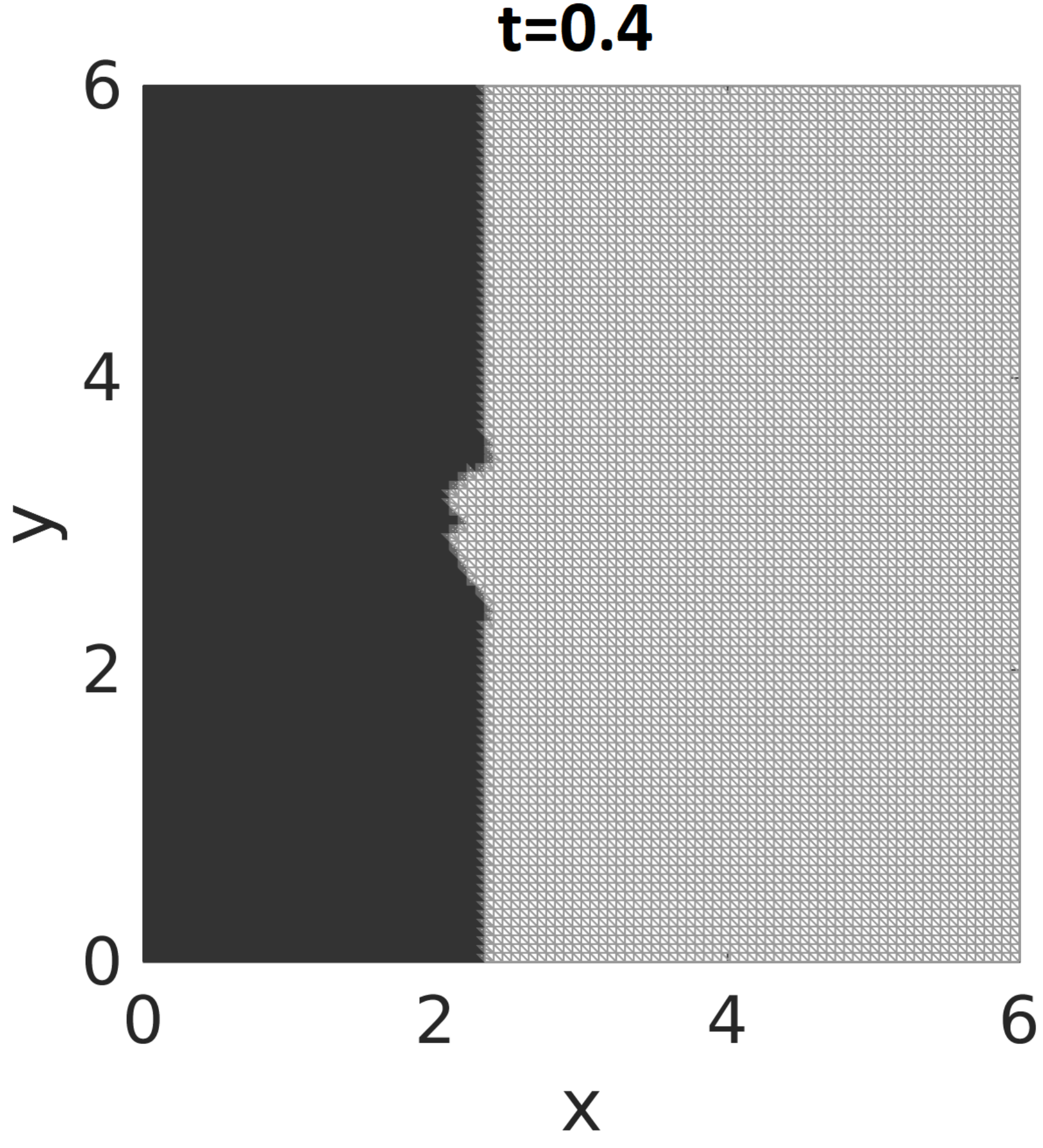}\hspace{1cm}
		\includegraphics[width=0.255\textwidth]{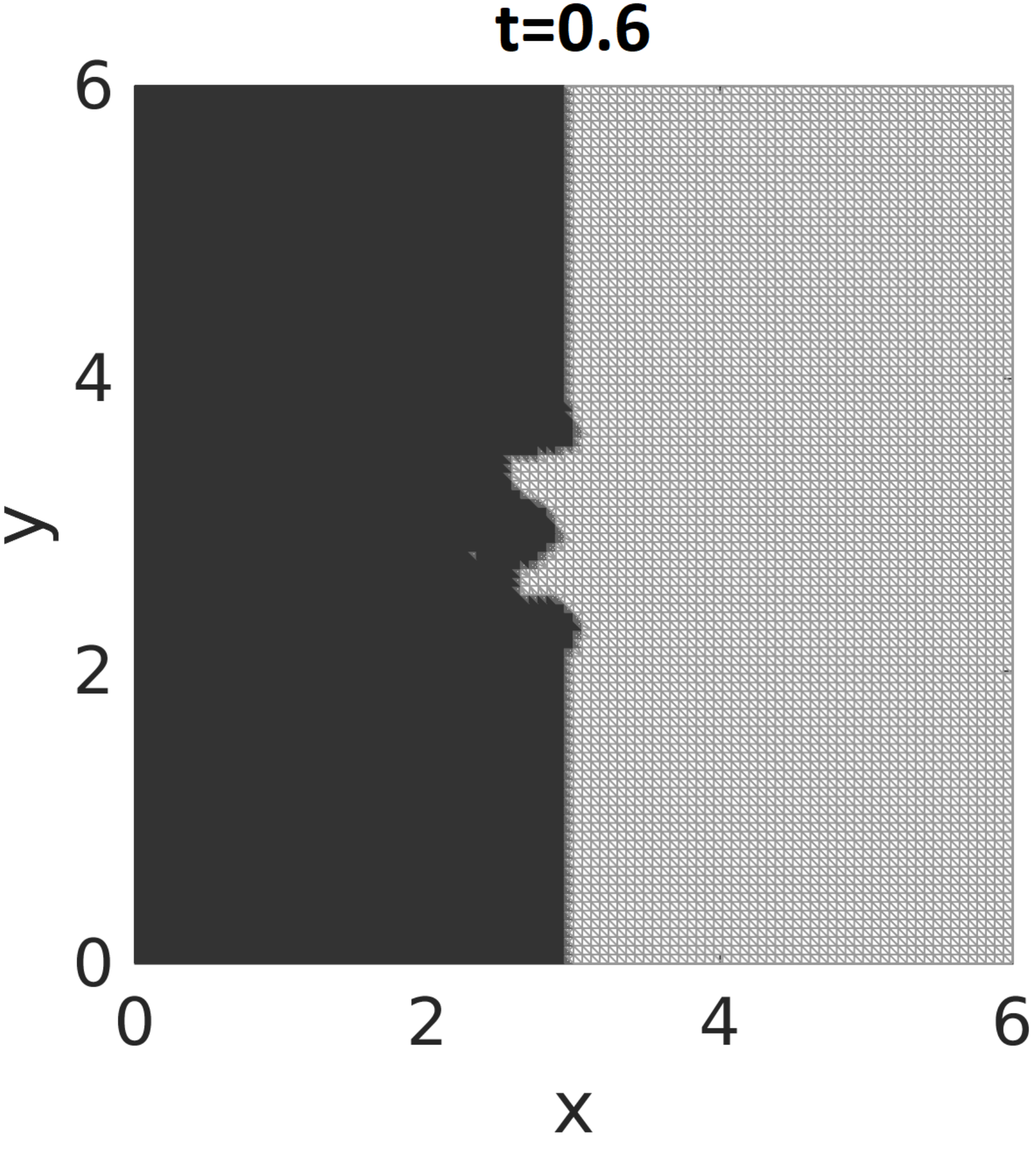}\hspace{1cm}
		\includegraphics[width=0.255\textwidth]{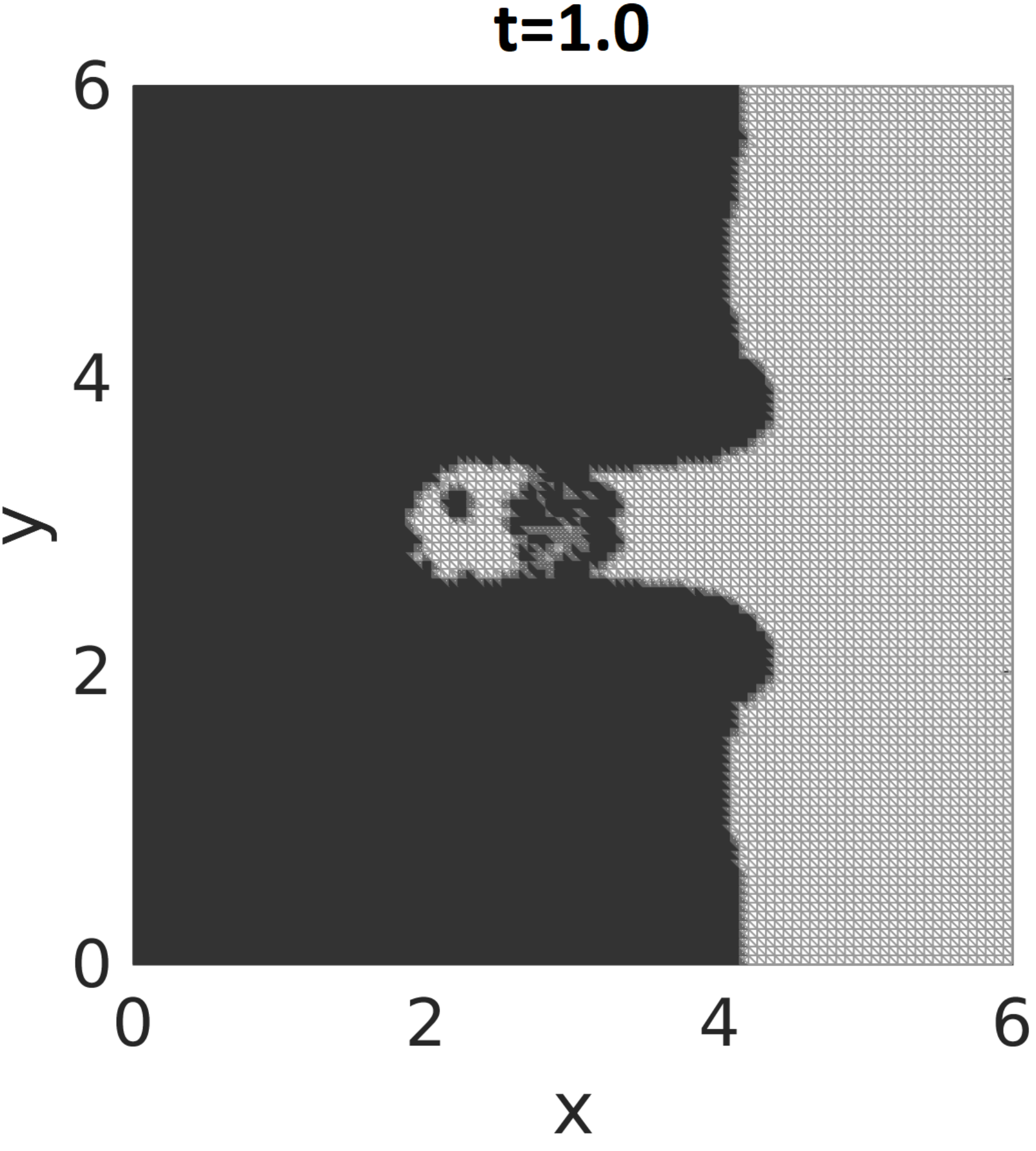}\\
		\vspace*{5mm}
		\includegraphics[width=0.255\textwidth]{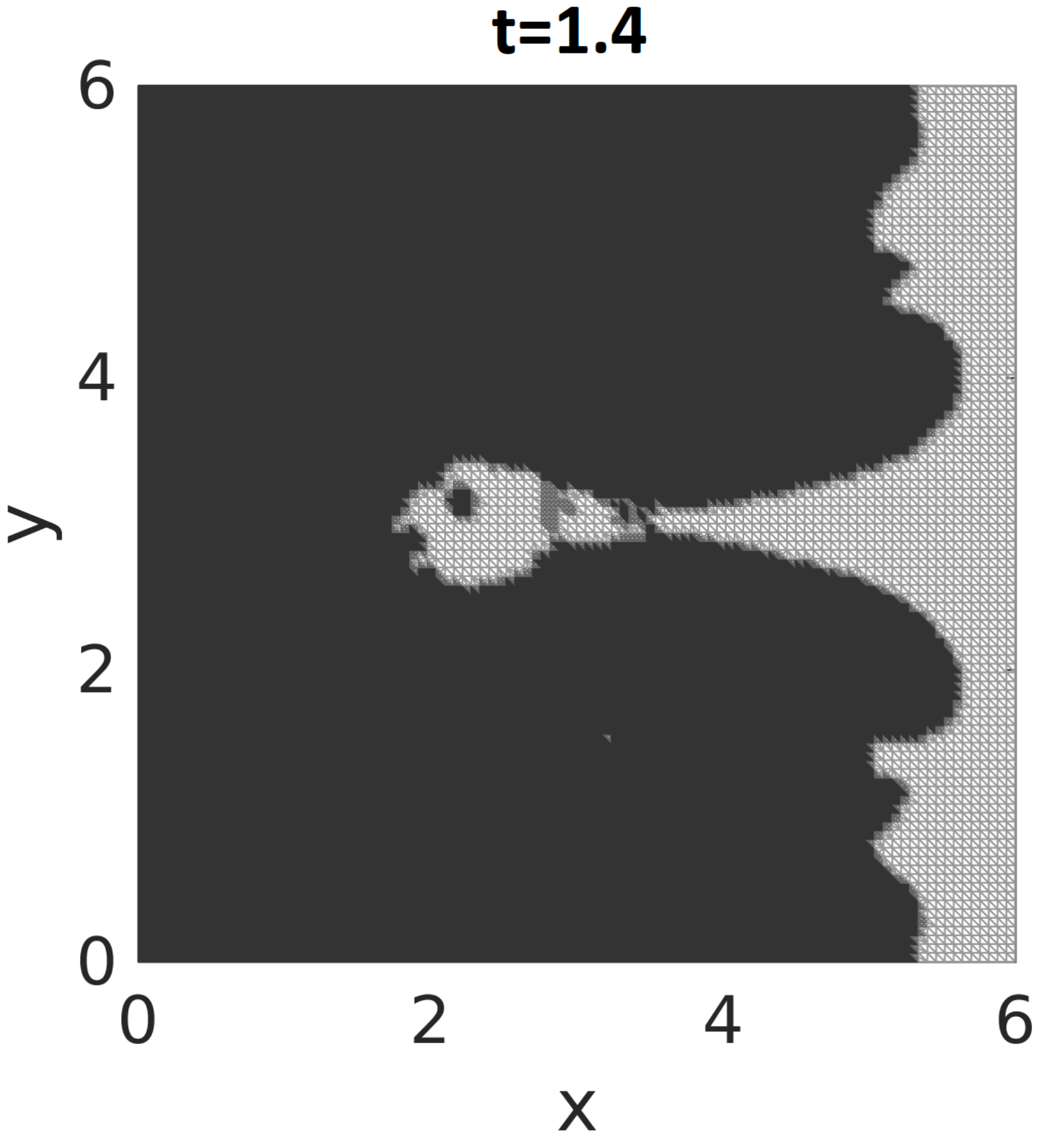}\hspace{1cm}
		\includegraphics[width=0.255\textwidth]{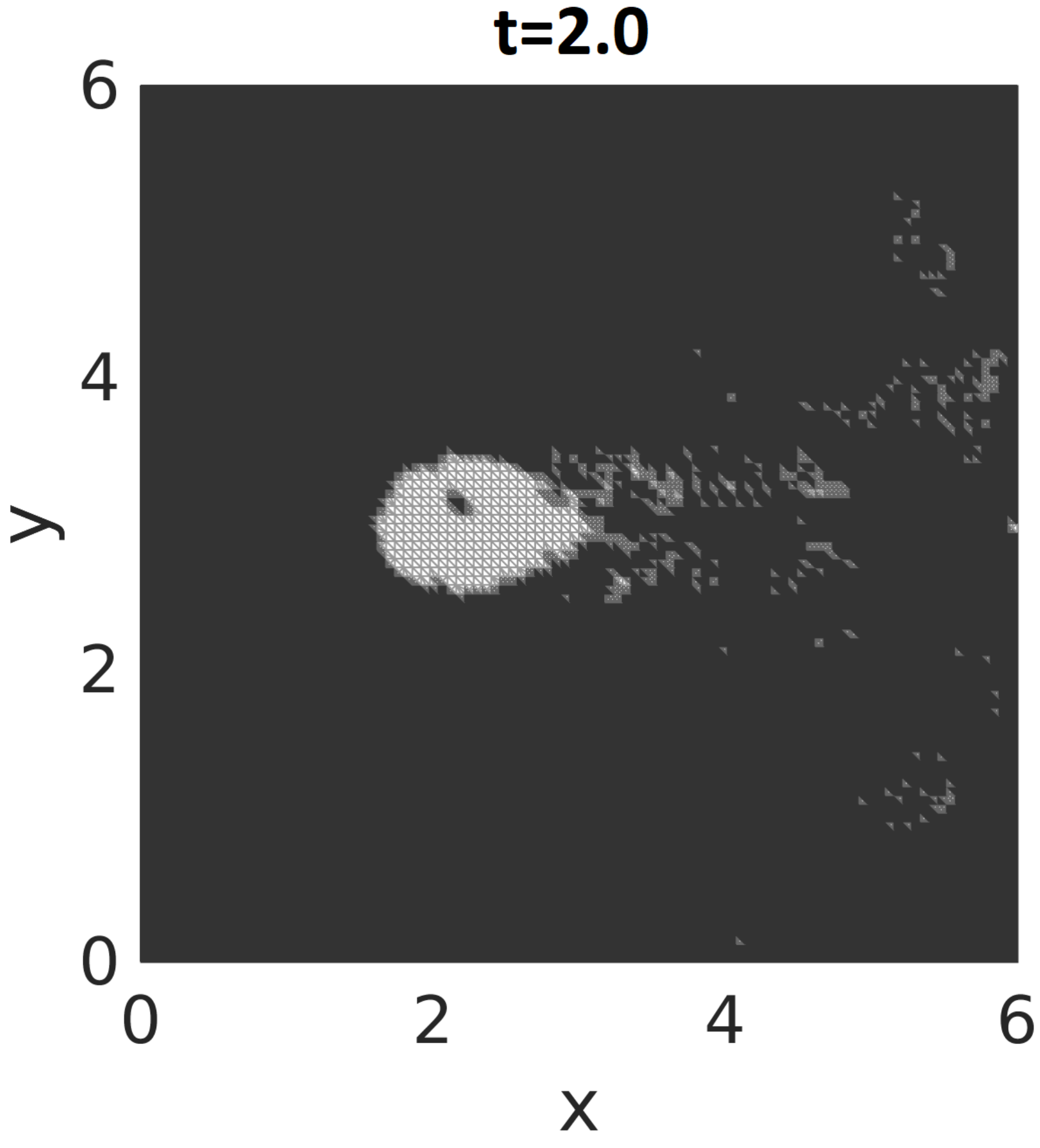}\hspace{1cm}
		\includegraphics[width=0.255\textwidth]{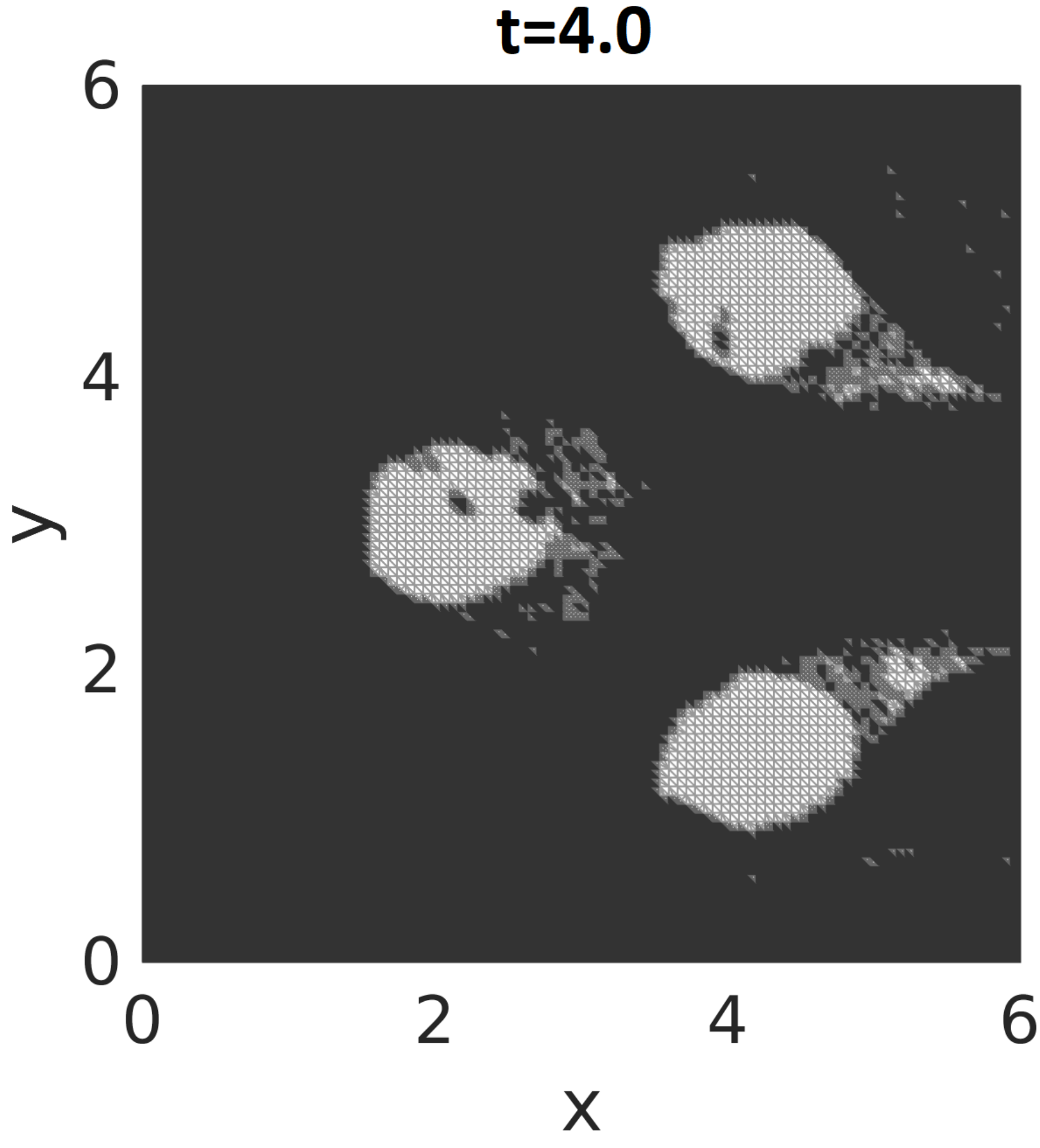}\\
		\vspace*{5mm}
		\caption{Example 3:  Adaptive mesh at different times with two levels of refinement $\mathcal{M}=2$.}\label{fig:ex3r2lm}
		\vspace*{2mm}
	\end{figure}

	Now, we continue the test for the proposed adaptive
	central-upwind method with the refined meshes generated from
	the uniform mesh $2\times 100\times 100$. The threshold for
	the grid refinement is set to $\omega=0.01 \max_j(e_j)$ in
	this example. In \fref{fig:ex33dr1l} and \fref{fig:ex33dr2l},
	we show the 3-D view of the simulated water computed by the
	adaptive scheme on the adaptive meshes with two cases of the
	highest refinement level $\mathcal{M}=1$ and $\mathcal{M}=2$. The
	behavior of the wave is similar to the result obtained by the
	central-upwind scheme without adaptivity,  see
	\fref{fig:ex33du1} and \fref{fig:ex33du2}. This means that the
	adaptive scheme performs well in simulating the wetting/drying
	processes. There are no non-physical spurious waves appear as
	a result of the simulation. We also present the contour lines
	of the water depth obtained on the adaptive meshes with
	$\mathcal{M}=1$ (\fref{fig:ex3r1lc}) and  $\mathcal{M}=2$
	(\fref{fig:ex3r2lc}). Clearly, the simulated solution captures
	correctly the reflections and interactions of the waves with
	no oscillations or disturbances showing up at the wet/dry
	interfaces. The
	considered adaptive meshes with $\mathcal{M}=1$ are plotted in \fref{fig:ex3r1lm} and with $\mathcal{M}=2$ in \fref{fig:ex3r2lm}. As we expected, the moving refined/de-refined regions match with the wetting and drying processes in the propagation of the flow.
	
	Finally, we present the $R_{CPU}$ ratios at time $t=1.0$ in
	\tref{tab:ex3ccpu}. The result in \tref{tab:ex3ccpu} shows
	that the average cost for the adaptive central-upwind method
	is about half of the cost for the central-upwind method
	without adaptivity. Note that at $t=1.0$, the refined region
	is larger than half of the computational domain, see
	\fref{fig:ex3r1lm} and \fref{fig:ex3r2lm}. Therefore, the
	numerical cost is not as remarkably reduced with the adaptive
	grid as the cost in example 2, see \tref{tab:ex2ccpu}. As illustrated, the  adaptive
	central-upwind method preserves the advantages of the
	well-balanced positivity preserving central-upwind scheme
	proposed in \citep{LAEK}, but at a less computational cost.
	
	\begin{table}[h!]
		\vspace*{5mm}
		\centering
		\begin{tabular}{ |c|c| c|c|c |}
			\hline
			\makecell{uniform mesh\\(cells)} &\makecell{adaptive mesh\\ $\mathcal{M}=1$\\(cells)}&\makecell{$\mathcal{R}_{CPU}$ \\with $\mathcal{M}=1$} &\makecell{adaptive mesh\\$\mathcal{M}=2$\\(cells)}&\makecell{$\mathcal{R}_{CPU}$ \\with $\mathcal{M}=2$}\\
			\hline
			$2\times100\times100$&15,064 &2.47&14,299   & 2.64\\
			$2\times200\times200$&59,252 &1.62&54,518 &  2.13\\
			$2\times400\times400$&238,485&1.53&217,075&1.69 \\
			\hline
			\multicolumn{2}{|c}{$\mathcal{R}_{CPU}$ average:}&
				\multicolumn{1}{c}{1.87} &\multicolumn{1}{c}{} &\multicolumn{1}{c|}{2.02}\\
			\hline
		\end{tabular}
		\caption{Example 3: The $R_{CPU}$ ratios at $t=1.0$.}\label{tab:ex3ccpu}
	\end{table}
	
	\section{Conclusion}\label{8sect05}
We have developed a new adaptive well-balanced and positivity
preserving central-upwind scheme on unstructured traingular meshes for
shallow water equations. The designed scheme is an extension and
improvement of the scheme in \citep{LAEK}. In addition, as a part of the adaptive
algorithm, we obained a robust local error indicator for the efficient
mesh refinement strategy. We conducted several challenging numerical
tests for shallow water equations and we demonstrated that the new
adaptive central-upwind scheme maintains important stability
properties (i.e.,  well-balanced and positivity-preserving properties) and delivers high-accuracy at a reduced computational cost.
		
	\section*{References} 

	\bibliography{ref}

\def\cprime{$'$} \def\cprime{$'$} \def\cprime{$'$} \def\cprime{$'$}
  \def\cprime{$'$} \def\cprime{$'$} \def\cprime{$'$} \def\cprime{$'$}
  \def\cprime{$'$} \def\cprime{$'$} \def\cprime{$'$} \def\cprime{$'$}
  \def\cprime{$'$} \def\cprime{$'$}
\begin{thebibliography}{10}

\bibitem{ABBLP}
{\sc E.~Audusse, F.~Bouchut, M.-O. Bristeau, R.~Klein, and B.~Perthame}, {\em A
  fast and stable well-balanced scheme with hydrostatic reconstruction for
  shallow water flows}, SIAM J. Sci. Comput., 25 (2004), pp.~2050--2065
  (electronic).

\bibitem{MR3738315}
{\sc P.~Azerad, J.-L. Guermond, and B.~Popov}, {\em Well-balanced second-order
  approximation of the shallow water equation with continuous finite elements},
  SIAM J. Numer. Anal., 55 (2017), pp.~3203--3224.

\bibitem{MR3071176}
{\sc A.~Bollermann, G.~Chen, A.~Kurganov, and S.~Noelle}, {\em A well-balanced
  reconstruction of wet/dry fronts for the shallow water equations}, J. Sci.
  Comput., 56 (2013), pp.~267--290.

\bibitem{MR2799646}
{\sc A.~Bollermann, S.~Noelle, and M.~Luk\'a{\v c}ov\'a-Medvi{\v d}ov\'a}, {\em
  Finite volume evolution {G}alerkin methods for the shallow water equations
  with dry beds}, Commun. Comput. Phys., 10 (2011), pp.~371--404.

\bibitem{AIPBEKP}
{\sc S.~Bryson, Y.~Epshteyn, A.~Kurganov, and G.~Petrova}, {\em Central
  {U}pwind {S}cheme on {T}riangular {G}rids for the {S}aint {V}enant {S}ystem
  of {S}hallow {W}ater {E}quations}, AIP Conference Proceedings, 1389 (2011),
  pp.~686--689.

\bibitem{MR2804645}
{\sc S.~Bryson, Y.~Epshteyn, A.~Kurganov, and G.~Petrova}, {\em Well-balanced
  positivity preserving central-upwind scheme on triangular grids for the
  {S}aint-{V}enant system}, ESAIM Math. Model. Numer. Anal., 45 (2011),
  pp.~423--446.

\bibitem{SDBL}
{\sc S.~Bryson and D.~Levy}, {\em Balanced central schemes for the shallow
  water equations on unstructured grids}, SIAM J. Sci. Comput., 27 (2005),
  pp.~532--552 (electronic).

\bibitem{MR2513368}
{\sc S.~Bunya, E.~J. Kubatko, J.~J. Westerink, and C.~Dawson}, {\em A wetting
  and drying treatment for the {R}unge-{K}utta discontinuous {G}alerkin
  solution to the shallow water equations}, Comput. Methods Appl. Mech. Engrg.,
  198 (2009), pp.~1548--1562.

\bibitem{CEHK}
{\sc A.~Chertock, Y.~Epshteyn, H.~Hu, and A.~Kurganov}, {\em High-order
  positivity-preserving hybrid finite-volume-finite-difference methods for
  chemotaxis systems}, Adv. Comput. Math., 44 (2018), pp.~327--350.

\bibitem{KM}
{\sc A.~{Chertock}, K.~{Fellner}, A.~{Kurganov}, A.~{Lorz}, and P.~A.
  {Markowich}}, {\em {Sinking, merging and stationary plumes in a coupled
  chemotaxis-fluid model: a high-resolution numerical approach}}, Journal of
  Fluid Mechanics, 694 (2012), pp.~155--190.

\bibitem{MR3229988}
{\sc A.~Chertock, A.~Kurganov, and Y.~Liu}, {\em Central-upwind schemes for the
  system of shallow water equations with horizontal temperature gradients},
  Numer. Math., 127 (2014), pp.~595--639.

\bibitem{Sai}
{\sc A.~J.~C. de~Saint-Venant}, {\em Th\`{e}orie du mouvement non-permanent des
  eaux, avec application aux crues des rivi\`{e}re at \`{a} l'introduction des
  mar\`{e}es dans leur lit.}, C.R. Acad. Sci. Paris, 73 (1871), pp.~147--154.

\bibitem{domingues2008adaptive}
{\sc M.~O. Domingues, S.~M. Gomes, O.~Roussel, and K.~Schneider}, {\em An
  adaptive multiresolution scheme with local time stepping for evolutionary
  pdes}, Journal of Computational Physics, 227 (2008), pp.~3758--3780.

\bibitem{Donat}
{\sc R.~Donat, M.~Mart\'{i}, A.~Mart\'{i}nez-Gavara, and P.~Mulet}, {\em
  Well-balanced adaptive mesh refinement for shallow water flows}, Journal of
  Computational Physics, 257 (2014), pp.~937--953.

\bibitem{MR2562498}
{\sc U.~S. Fjordholm, S.~Mishra, and E.~Tadmor}, {\em Energy preserving and
  energy stable schemes for the shallow water equations}, in Foundations of
  computational mathematics, {H}ong {K}ong 2008, vol.~363 of London Math. Soc.
  Lecture Note Ser., Cambridge Univ. Press, Cambridge, 2009, pp.~93--139.

\bibitem{MR2799526}
\leavevmode\vrule height 2pt depth -1.6pt width 23pt, {\em Well-balanced and
  energy stable schemes for the shallow water equations with discontinuous
  topography}, J. Comput. Phys., 230 (2011), pp.~5587--5609.

\bibitem{GHS}
{\sc T.~Gallou{\"e}t, J.-M. H{\'e}rard, and N.~Seguin}, {\em Some approximate
  {G}odunov schemes to compute shallow-water equations with topography},
  Comput. \& Fluids, 32 (2003), pp.~479--513.

\bibitem{DLGeorge}
{\sc D.~L. George}, {\em Finite volume methods and adaptive refinement for
  tsunami propagation and inundation}, {P}h.{D}, University of Washington,
  (2006).

\bibitem{GHAZIZADEH2020104633}
{\sc M.~A. Ghazizadeh, A.~Mohammadian, and A.~Kurganov}, {\em An adaptive
  well-balanced positivity preserving central-upwind scheme on quadtree grids
  for shallow water equations}, Computers $\&$ Fluids, 208 (2020), p.~104633.

\bibitem{GST}
{\sc S.~Gottlieb, C.-W. Shu, and E.~Tadmor}, {\em Strong stability-preserving
  high-order time discretization methods}, SIAM Rev., 43 (2001), pp.~89--112.

\bibitem{NHA}
{\sc N.~Hannoun and V.~Alexiades}, {\em Issues in adaptive mesh refinement
  implementation}, Electronic Journal of Differential Equations, 15 (2007),
  pp.~141--151.

\bibitem{Jin}
{\sc S.~Jin}, {\em A steady-state capturing method for hyperbolic systems with
  geometrical source terms}, M2AN Math. Model. Numer. Anal., 35 (2001),
  pp.~631--645.

\bibitem{JW}
{\sc S.~Jin and X.~Wen}, {\em Two interface-type numerical methods for
  computing hyperbolic systems with geometrical source terms having
  concentrations}, SIAM J. Sci. Comput., 26 (2005), pp.~2079--2101.

\bibitem{MR2126235}
{\sc S.~Karni and A.~Kurganov}, {\em Local error analysis for approximate
  solutions of hyperbolic conservation laws}, Adv. Comput. Math., 22 (2005),
  pp.~79--99.

\bibitem{KLab}
{\sc A.~Kurganov and D.~Levy}, {\em Central-upwind schemes for the
  {S}aint-{V}enant system}, M2AN Math. Model. Numer. Anal., 36 (2002),
  pp.~397--425.

\bibitem{KurLiu}
{\sc A.~Kurganov and Y.~Liu}, {\em New adaptive artificial viscosity method for
  hyperbolic systems of conservation laws}, J. Comput. Phys., 231 (2012),
  pp.~8114--8132.

\bibitem{MR3187922}
{\sc A.~Kurganov and J.~Miller}, {\em Central-upwind scheme for
  {S}avage-{H}utter type model of submarine landslides and generated tsunami
  waves}, Comput. Methods Appl. Math., 14 (2014), pp.~177--201.

\bibitem{KNP}
{\sc A.~Kurganov, S.~Noelle, and G.~Petrova}, {\em Semi-discrete central-upwind
  schemes for hyperbolic conservation laws and {H}amilton-{J}acobi equations},
  SIAM J. Sci. Comput., 23 (2001), pp.~707--740.

\bibitem{KP05}
{\sc A.~Kurganov and G.~Petrova}, {\em Central-upwind schemes on triangular
  grids for hyperbolic systems of conservation laws}, Numer. Methods Partial
  Differential Equations, 21 (2005), pp.~536--552.

\bibitem{KP1}
{\sc A.~Kurganov and G.~Petrova}, {\em A second-order well-balanced positivity
  preserving scheme for the {S}aint-{V}enant system}, Commun. Math. Sci., 5
  (2007), pp.~133--160.

\bibitem{MR1756766}
{\sc A.~Kurganov and E.~Tadmor}, {\em New high-resolution central schemes for
  nonlinear conservation laws and convection-diffusion equations}, J. Comput.
  Phys., 160 (2000), pp.~241--282.

\bibitem{MR1763829}
\leavevmode\vrule height 2pt depth -1.6pt width 23pt, {\em New high-resolution
  semi-discrete central schemes for {H}amilton-{J}acobi equations}, J. Comput.
  Phys., 160 (2000), pp.~720--742.

\bibitem{LeV}
{\sc R.~LeVeque}, {\em Balancing source terms and flux gradients in
  high-resolution {G}odunov methods: the quasi-steady wave-propagation
  algorithm}, J. Comput. Phys., 146 (1998), pp.~346--365.

\bibitem{LeVbook}
\leavevmode\vrule height 2pt depth -1.6pt width 23pt, {\em Finite volume
  methods for hyperbolic problems}, Cambridge Texts in Applied Mathematics,
  Cambridge University Press, Cambridge, 2002.

\bibitem{MR3640671}
{\sc M.~Li, P.~Guyenne, F.~Li, and L.~Xu}, {\em A positivity-preserving
  well-balanced central discontinuous {G}alerkin method for the nonlinear
  shallow water equations}, J. Sci. Comput., 71 (2017), pp.~994--1034.

\bibitem{LAEK}
{\sc X.~Liu, J.~Albright, Y.~Epshteyn, and A.~Kurganov}, {\em Well-balanced
  positivity preserving central-upwind scheme with a novel wet/dry
  reconstruction on triangular grids for the {S}aint-{V}enant system}, Journal
  of Computational Physics, 374 (2018), pp.~213 -- 236.

\bibitem{MORE}
{\sc M.~{Moreira Lopes}, M.~O. Domingues, K.~Schneider, and O.~Mendes}, {\em
  Local time-stepping for adaptive multiresolution using natural extension of
  runge–kutta methods}, Journal of Computational Physics, 382 (2019), pp.~291
  -- 318.

\bibitem{MR1047564}
{\sc H.~Nessyahu and E.~Tadmor}, {\em Nonoscillatory central differencing for
  hyperbolic conservation laws}, J. Comput. Phys., 87 (1990), pp.~408--463.

\bibitem{NPPN}
{\sc S.~Noelle, N.~Pankratz, G.~Puppo, and J.~Natvig}, {\em Well-balanced
  finite volume schemes of arbitrary order of accuracy for shallow water
  flows}, J. Comput. Phys., 213 (2006), pp.~474--499.

\bibitem{PS}
{\sc B.~Perthame and C.~Simeoni}, {\em A kinetic scheme for the
  {S}aint-{V}enant system with a source term}, Calcolo, 38 (2001),
  pp.~201--231.

\bibitem{Rus1}
{\sc G.~Russo}, {\em Central schemes for balance laws}, in Hyperbolic problems:
  theory, numerics, applications, {V}ol. {I}, {II} ({M}agdeburg, 2000),
  vol.~141 of Internat. Ser. Numer. Math., 140, Birkh\"auser, Basel, 2001,
  pp.~821--829.

\bibitem{Rus2}
{\sc G.~Russo}, {\em Central schemes for conservation laws with application to
  shallow water equations}, in {T}rends and applications of mathematics to
  mechanics: {S}{T}{A}{M}{M} 2002, {S}. {R}ionero and {G}. {R}omano (eds.),
  (2005), pp.~225--246.

\bibitem{MR3315267}
{\sc M.~L. S\ae~tra, A.~R. Brodtkorb, and K.-A. Lie}, {\em Efficient
  {GPU}-implementation of adaptive mesh refinement for the shallow-water
  equations}, J. Sci. Comput., 63 (2015), pp.~23--48.

\bibitem{MR3440159}
{\sc H.~Shirkhani, A.~Mohammadian, O.~Seidou, and A.~Kurganov}, {\em A
  well-balanced positivity-preserving central-upwind scheme for shallow water
  equations on unstructured quadrilateral grids}, Comput. \& Fluids, 126
  (2016), pp.~25--40.

\bibitem{MR1111445}
{\sc E.~Tadmor}, {\em Local error estimates for discontinuous solutions of
  nonlinear hyperbolic equations}, SIAM J. Numer. Anal., 28 (1991),
  pp.~891--906.

\bibitem{XS1}
{\sc Y.~Xing and C.-W. Shu}, {\em High order finite difference {WENO} schemes
  with the exact conservation property for the shallow water equations}, J.
  Comput. Phys., 208 (2005), pp.~206--227.

\bibitem{XS2}
\leavevmode\vrule height 2pt depth -1.6pt width 23pt, {\em High order
  well-balanced finite volume {WENO} schemes and discontinuous {G}alerkin
  methods for a class of hyperbolic systems with source terms}, J. Comput.
  Phys., 214 (2006), pp.~567--598.

\end{thebibliography}
	\bibliographystyle{siam}
\end{document}